\documentclass[a4paper,12pt,dvips]{article} 
\usepackage{amssymb}
\usepackage{amscd}
\usepackage{amsmath}
\usepackage{graphicx}
\usepackage{latexsym}
\usepackage{enumerate}

\usepackage{theorem}
\theoremstyle{plain}
\theorembodyfont{\itshape}
\newtheorem{theorem}{Theorem}[section]

\newtheorem{lemma}[theorem]{Lemma}
\newtheorem{corollary}[theorem]{Corollary}
\theorembodyfont{\rmfamily}
\newtheorem{definition}[theorem]{Definition} 
\newtheorem{remark}[theorem]{Remark}
\newtheorem{problem}[theorem]{Problem} 
\newtheorem{example}[theorem]{Example} 
 
\def\Ad{\mathrm{Ad}} 

\def\Aut{\mathrm{Aut}}
\def\det{\mathrm{det}}
\def\End{\mathrm{End}}
\def\Hom{\mathrm{Hom}}
\def\IMF{\mathrm{IMF}}
\def\Ker{\mathrm{Ker}}
\def\Kl{\mathrm{Kl}} 
\def\id{\mathrm{id}}
\def\incl{\sc{\mathrm{inclusion}}}

\def\pardeg{\mathrm{pardeg}}
\def\PVI{\mathrm{P}_{\mathrm{VI}}}
\def\HVI{\mathrm{H}_{\mathrm{VI}}}
\def\r{\mathrm{r}}
\def\rank{\mathrm{rank}} 
\def\Res{\mathrm{Res}}
\def\RH{\mathrm{RH}} 
\def\rh{\mathrm{rh}}

\def\s{\mathrm{s}} 

\def\Spec{\mathrm{Spec}}
\def\Tr{\mathrm{Tr}}
\def\CC{\mathcal{C}}
\def\E{\mathcal{E}}
\def\F{\mathcal{F}}
\def\I{\mathcal{I}}
\def\K{\mathcal{K}}
\def\L{\mathcal{L}}
\def\M{\mathcal{M}}
\def\O{\mathcal{O}}
\def\R{\mathcal{R}}
\def\bR{\mbox{\bm $\mathcal{R}$}}

\def\C{\mathbb{C}}
\def\D{\mathbb{D}}
\def\P{\mathbb{P}}

\def\Z{\mathbb{Z}} 
\def\bIMF{\mathbf{IMF}}
\def\Wall{\mathbf{Wall}} 
\def\a{\alpha} 
\def\b{\beta}
\def\ga{\gamma} 
\def\l{\lambda}
\def\G{\varGamma}
\def\Sig{\varSigma} 
\def\k{\kappa}
\def\fraksl{\mathfrak{sl}}
\def\Om{\Omega}
\def\si{\sigma}
\def\th{\theta} 
\def\Th{\Theta}
\def\ve{\varepsilon}  
\def\RR{\sf{R}}
\def\car{\curvearrowright}

\def\lra{\longrightarrow}
  
\def\ra{\rightarrow}
\def\mt{\mapsto}
\def\lmt{\longmapsto}  
\def\ol{\overline}
\def\ul{\underline}
\def\aa{$\phan{aaaaaaaaaaaaaa}$}
\def\ci{\circ}
\def\dfrac#1#2{{\displaystyle\frac{#1}{#2}}}
\def\ds{\displaystyle}
\def\sc{\scriptstyle}
\def\bm{\boldmath}

\def\phan{\phantom}   
\def\la{\langle}
\def\ra{\rangle}
\def\op{\oplus}
\def\ot{\otimes}
\def\rot{\rotatebox} 
\def\W{W(D_4^{(1)})}  
\title{\bf Dynamics of the Sixth Painlev\'e Equation}  
\author{Michi-aki Inaba, Katsunori Iwasaki, and Masa-Hiko Saito 
\\ \\
Kyushu University and Kobe University}   
\date{}  
\begin{document}
\maketitle
\begin{abstract} 
The sixth Painlev\'e equation is hiding extremely rich geometric 
structures behind its outward appearance.  
This article tries to give as a total picture as possible of its 
dynamical nature, based on the Riemann-Hilbert approach recently 
developed by the authors, using various techniques from algebraic 
geometry. 
A good part of the contents is extended to Garnier systems, 
while this article is restricted to the original sixth Painlev\'e 
equation.   
\end{abstract} 
\section{Introduction}  \label{sec:intro} 
The sixth Painlev\'e equation $\PVI = \PVI(\k)$ is 
among the six kinds of differential equations that were discovered 
by Painlev\'e \cite{Painleve} and his student 
Gambier \cite{Gambier} around the turn of the twentieth century. 
It is a second-order nonlinear ordinary differential equation 
with an independent variable $x \in \P^1-\{0,1,\infty\}$ 
and an unknown function $q = q(x)$,
\begin{equation} \label{eqn:PVI}
\begin{array}{rcl}
q_{xx} &=& \dfrac{1}{2} 
\left(\dfrac{1}{q}+\dfrac{1}{q-1}+\dfrac{1}{q-x} \right) q_x^2 
- \left(\dfrac{1}{x}+\dfrac{1}{x-1}+\dfrac{1}{q-x} \right) q_x \\[6mm]
&+& \dfrac{q(q-1)(q-x)}{2x^2(x-1)^2} 
\left\{\k_4^2 - \k_1^2 \dfrac{x}{q^2} + 
\k_2^2 \dfrac{x-1}{(q-1)^2} + 
(1-\k_3^2) \dfrac{x(x-1)}{(q-x)^2} \right\}, 
\end{array} 
\end{equation}
depending on complex parameters 
$\k = (\k_0,\k_1,\k_2,\k_3,\k_4)$ in a 
$4$-dimensional affine space 
\begin{equation} \label{eqn:K}
\mathcal{K} = \{\k = 
(\k_0,\k_1,\k_2,\k_3,\k_4) \in \C^5 \,:\, 
2 \k_0 + \k_1 + \k_2 + \k_3 +\k_4 = 1\}. 
\end{equation} 
\par
This highly nonlinear and seemingly rather ugly equation is 
only a small visible part of a more substantial entity. 
The large invisible part has extremely rich geometric 
structures that are related to symplectic geometry, moduli 
spaces of stable parabolic connections, moduli spaces of 
representations of a fundamental group, Riemann-Hilbert 
correspondence, geometry of cubic surfaces, braid and modular 
groups, simple isolated singularities and their resolutions 
of singularities, affine Weyl groups, discrete dynamical 
systems, and so on. 
The aim of this survey article is to discuss various aspects 
of these illuminating structures, trying to give as a total 
picture as possible of the sixth Painlev\'e equation. 
\par
Among other features, Painlev\'e equation is primarily 
a dynamical system and a dynamical system in general is 
characterized by two aspects: {\it laws} and {\it phenomena}. 
Mathematically, laws refer to equations that govern the 
dynamics, symmetries of the system, etc., while 
phenomena refer to solutions of the equations, 
behaviors of trajectories, etc. 
These two aspects often show a sharp contrast.  
For example, in classical mechanics, the simple laws of 
Newton create extremely rich and complicated phenomena. 
The Painlev\'e dynamics is also in this case, being 
algebraic in its laws and transcendental in its phenomena 
(see Table \ref{tab:aspects}).    
\begin{table}[t] 
\begin{center}
\begin{tabular}{|l|l|l|}
\hline 
aspect & contents & nature \\ 
\hline 
laws  & equations, symmetry, \dots & algebraic \\ 
\hline 
phenomena & solutions, trajectories, \dots & transcendental \\ 
\hline
\end{tabular}
\end{center} 
\caption{Two aspects of Painlev\'e equation} 
\label{tab:aspects}
\end{table}
For comparison, we should remark that there exists an 
interesting dynamics whose laws are already transcendental, 
like a dynamics on a K3 surface recently explored by 
McMullen \cite{McMullen}, who showed that the existence of 
Siegel disks implies the transcendence of the K3 surface.
\par
Generally speaking, the two pricipal approaches to dynamical 
systems are perhaps: 
\par\medskip\noindent
\centerline{(L) Lyapunov's methods, \qquad (C) 
conjugacy methods.}  
\par\medskip\noindent
In Lyapunov's methods (L), we examine, control, or confine the 
behaviors of trajectories by estimating suitable functions 
called ``Lyapunov functions''.  
Main tools of the methods are estimations by inequalities. 
On the other hand, in the conjugacy methods (C),   
we try to find a ``conjugacy map'' that converts the 
difficult dynamical system we want to study to a more 
tractable one, to extract informations from the latter, 
and to send feedback to the former 
(see \S\ref{subsec:conjugacy} for more details).   
Our approach to the Painlev\'e equation, which we call the 
{\it Riemann-Hilbert approach}, falls into this category (C),  
making use of Riemann-Hilbert correspondence as a  
conjugacy map between Painlev\'e flow and isomonodromic flow.  
\par  
Of course, the Riemann-Hilbert approach is closely related 
to the isomonodromic approach represented by the classical 
works of Fuchs \cite{Fuchs}, Schlesinger \cite{Schlesinger}, 
Garnier \cite{Garnier}, Jimbo, Miwa and Ueno \cite{JMU,JM} and 
others, but differs from the latter in its definitive 
employment of the method of conjugacy maps and in its extensive 
use of a complete solution to the Riemann-Hilbert problem. 
The Riemann-Hilbert approach {\sl a priori} has a global 
nature once Riemann-Hilbert correspondence is 
formulated appropriately, while the 
isomonodromic approach mostly stands on the  
infinitesimal point of view and pays little attention  
to the target space of Riemann-Hilbert correspondence, 
namely, moduli space of monodromy representations. 
In the Riemann-Hilbert appraoch, we consciously 
distinguish the Painlev\'e flow on the moduli space 
of stable parabolic connections and the isomonodromic 
flow on the moduli space of monodromy representations, 
and build a bridge between them via the Riemann-Hilbert 
correspondence. 
\par    
This approach has been explored by 
Iwasaki \cite{Iwasaki1,Iwasaki2,Iwasaki3,Iwasaki4}, 
Hitchin \cite{Hitchin}, Kawai \cite{Kawai1,Kawai2}, 
Boalch \cite{Boalch1,Boalch2,Boalch3}, 
Dubrovin and Mazzocco \cite{DM} and others. 
Recently it was thoroughly developed by Inaba, Iwasaki 
and Saito \cite{IIS1,IIS2}. 
The exposition of this article is largely based on the 
contents of the last papers.  
We focus our attension on the original case of $\PVI$   
in the merit of presenting, for the most basic model, 
those materials which can be expected to be universal 
throughout various generalizations.   
A good part of the contents is extended to Garnier 
systems, a several-variable version of $\PVI$; 
see \cite{IIS2}.   
\par  
In addtion to the general methods represented by 
approaches (L) and (C), which are conceivable in general 
situations, there are also various particular methods 
applicable to various particular dynamical systems. 
For example, for the class of dynamical systems 
that are called completely integrable systems, 
there exist  
\par\medskip\noindent
\centerline{(CI) Completely integrating methods,}    
\par\medskip\noindent
which are characterized by such keywords as 
$\tau$-functions, bilnear equations, Lax pairs, 
separations of variables, combinatorics and representation 
theory, etc. 
Painlev\'e equations are usually thought of as 
a member of this class and many works have been 
done from this point of view.  
See Conte \cite{Conte}, Noumi \cite{Noumi} and 
the references therein. 
But we shall not touch on this aspect 
in this article. 
Among other things, we wish to lay a sound foundation on 
the sixth Painlev\'e equation to such an extent that it 
can be a basis for the investigations into the 
{\it transcendental} natures of $\PVI$.   
To do so, many things should be done, even within 
the general framework of dynamical systems,  
before entering into those subjects which 
are particular to integrable systems. 
Therefore the integrable aspects should be 
discussed later and elsewhere. 
\par 
Lyapunov-type approaches to Painlev\'e equations 
will not completely be discussed in this article. 
There have been long traditions as well as recent 
developments of establishing Painlev\'e property 
by these methods.  
We refer to Painlev\'e \cite{Painleve}, 
Hukuhara \cite{Hukuhara} (see Okamoto and 
Takano \cite{OT} for a part of these unpublished lectures), 
Joshi and Kruskal \cite{JK}, Steinmetz \cite{Steinmetz}, 
Iwasaki, Kimura, Shimomura and Yoshida \cite{IKSY}, 
Gromak, Laine and Shimomura \cite{GLS} and the references therein. 
\par\vspace{0.3cm}\noindent 
\begin{center}
\begin{minipage}{.70\linewidth}
\centerline{\bf Table of Contents} 
\vspace{0.1cm}
\begin{enumerate}
\item[1.] Introduction 
\vspace{-0.2cm}
\item[2.] Painlev\'e Equation as a Dynamical System 
\vspace{-0.2cm} 
\item[3.] Moduli Spaces of Parabolic Connections (Phase Spaces)
\vspace{-0.2cm}
\item[4.] Riemann-Hilbert Correspondence (Conjugacy Map)
\vspace{-0.2cm} 
\item[5.] Isomonodromic flow and Painlev\'e flow
\vspace{-0.2cm}
\item[6.] Family of Affine Cubic Surfaces
\vspace{-0.2cm} 
\item[7.] B\"acklund Transformations (Symmetry)
\vspace{-0.2cm} 
\item[8.] Nonlinear Monodromy (Poincar\'e Return Maps)
\vspace{-0.2cm}
\item[9.] Singularities and Riccati Solutions (Classical Trajectories) 
\vspace{-0.2cm}
\item[10.] Canonical Coordinates
\vspace{-0.2cm}
\item[11.] Summary  
\end{enumerate} 
\end{minipage}
\end{center} 
\par\vspace{0.3cm} 
The plan of this article is as follows: 
In Section 2 we introduce a general formalism of dynamical 
systems and push $\PVI$ into this framework. 
We present the Guiding Diagram that encodes major 
dynamical natures of $\PVI$.     
Section 3 is devoted to the construction of moduli spaces 
of stable parabolic connections, which, in the dynamical 
context, means the construction of phase spaces of $\PVI$. 
In Section 4, after setting up moduli spaces of monodromy 
representations, we formulate Riemann-Hilbert 
correspondence, $\RH$, and settle Riemann-Hilbert problems  
in suitable ways.  
In the dynamical context, theses parts correspond to the 
construction of conjugacy maps. 
In Section 5 we formulate isomonodromic flows $\F_{\IMF}$ 
and Painlev\'e flows $\F_{\PVI}$ in such a manner that 
$\RH$ yields analytic conjugacy from $\F_{\PVI}$ to 
$\F_{\IMF}$. 
Section 6 is devoted to the construction of a family of 
affine cubic surfaces, which enables us to describe all 
the previous constructions more explicitly. 
In Section 7 we give a characterization of B\"acklund 
transformations, namely, the symmetries of $\PVI$, in 
terms of Riemann-Hilbert correspondence. 
Section 8 describes the nonlinear monodromy or the  
Poincar\'e return map of $\PVI$ that extracts the 
global natures of trajectories of $\PVI$. 
In Section 9 we characterize the classical components 
of $\PVI$, called the Riccati flows, in terms of  
singularities on cubic surfaces and their 
resolutions of singularities. 
In Section 10 we construct canonical coordinate systems 
of moduli spaces (phase spaces) which make it possible 
to write down the Painlev\'e dynamics explicitly. 
This article is closed with a brief summary, 
especially with the Closing Diagram, in Section 11. 
\section{Painlev\'e Equation as a Dynamical System} 
\label{sec:dynamics}
A total picture of the sixth Painlev\'e equation is 
most clearly described in the framework of dynamical 
systems, or, more specifically as a time-dependent 
Hamiltonian system with Painlev\'e property. 
So we begin by establishing a general formalism of 
dynamical systems, based on which we shall develop 
our whole story. 
\subsection{General Formalism of Dynamical Systems} 
\label{subsec:formalism}
\begin{definition}[Time-Dependent Dynamical System] 
\label{def:dynamics}  
A {\it time-dependent dynamical system} $(M,\F)$ is a 
fibration $\pi : M \to T$ of complex manifolds together with 
a complex foliation $\F$ on $M$ that is transverse to each   
fiber $M_t = \pi^{-1}(t)$, $t \in T$. 
The total space $M$ is referred to as the {\it phase space}, 
while the base space $T$ is called the {\it space of 
time-variables}. 
Moreover, the fiber $M_t$ is called the 
{\it space of initial conditions} at time $t$.  
\end{definition} 
\par 
The space of initial conditions becomes a meaningful concept 
if the dynamical system enjoys Painlev\'e property. 
It is this property that makes it possible to think of 
Poincar\'e return maps or the nonlinear monodromy, 
which is the discrete dynamical system on a space of initial 
conditions that represents the global natures of a continuous 
dynamical system.  
\begin{figure}[t] 
\begin{center}
\unitlength 0.1in
\begin{picture}(36.60,31.80)(1.60,-32.10)
%
\special{pn 20}%
\special{pa 2840 1300}%
\special{pa 2960 1310}%
\special{fp}%
\special{sh 1}%
\special{pa 2960 1310}%
\special{pa 2895 1285}%
\special{pa 2907 1306}%
\special{pa 2892 1324}%
\special{pa 2960 1310}%
\special{fp}%
%
\special{pn 20}%
\special{pa 610 3010}%
\special{pa 3810 3010}%
\special{fp}%
%
\special{pn 20}%
\special{pa 1800 420}%
\special{pa 1800 2220}%
\special{fp}%
%
\special{pn 20}%
\special{pa 3000 410}%
\special{pa 3000 2220}%
\special{fp}%
%
\special{pn 13}%
\special{pa 3000 2210}%
\special{pa 3000 3020}%
\special{da 0.070}%
%
\special{pn 13}%
\special{pa 1800 2220}%
\special{pa 1800 3020}%
\special{da 0.070}%
%
\special{pn 20}%
\special{pa 2300 2370}%
\special{pa 2300 2810}%
\special{fp}%
\special{sh 1}%
\special{pa 2300 2810}%
\special{pa 2320 2743}%
\special{pa 2300 2757}%
\special{pa 2280 2743}%
\special{pa 2300 2810}%
\special{fp}%
\put(16.8000,-3.3000){\makebox(0,0)[lb]{$M_t$}}%
\put(29.3000,-3.3000){\makebox(0,0)[lb]{$M_{t'}$}}%
\put(20.9000,-26.1000){\makebox(0,0)[lb]{$\pi$}}%
\put(1.6000,-30.8000){\makebox(0,0)[lb]{$T$}}%
\put(1.7000,-14.5000){\makebox(0,0)[lb]{$M$}}%
%
\special{pn 20}%
\special{sh 0.600}%
\special{ar 1800 3010 81 81  0.0000000 6.2831853}%
%
\special{pn 20}%
\special{sh 0.600}%
\special{ar 2990 3000 81 81  0.0000000 6.2831853}%
%
\special{pn 20}%
\special{sh 0.600}%
\special{ar 1800 1900 81 81  0.0000000 6.2831853}%
%
\special{pn 20}%
\special{sh 0.600}%
\special{ar 3000 1580 81 81  0.0000000 6.2831853}%
\put(24.5000,-28.8000){\makebox(0,0)[lb]{$\ga$}}%
%
\special{pn 20}%
\special{pa 1790 1890}%
\special{pa 1824 1890}%
\special{pa 1858 1890}%
\special{pa 1891 1889}%
\special{pa 1925 1888}%
\special{pa 1958 1887}%
\special{pa 1991 1884}%
\special{pa 2024 1881}%
\special{pa 2056 1878}%
\special{pa 2088 1873}%
\special{pa 2119 1867}%
\special{pa 2150 1859}%
\special{pa 2180 1850}%
\special{pa 2210 1840}%
\special{pa 2238 1828}%
\special{pa 2266 1814}%
\special{pa 2293 1798}%
\special{pa 2319 1780}%
\special{pa 2345 1761}%
\special{pa 2370 1740}%
\special{pa 2395 1719}%
\special{pa 2420 1699}%
\special{pa 2445 1679}%
\special{pa 2470 1659}%
\special{pa 2497 1642}%
\special{pa 2524 1627}%
\special{pa 2552 1614}%
\special{pa 2581 1602}%
\special{pa 2610 1593}%
\special{pa 2640 1585}%
\special{pa 2671 1579}%
\special{pa 2703 1574}%
\special{pa 2735 1570}%
\special{pa 2767 1568}%
\special{pa 2800 1566}%
\special{pa 2834 1566}%
\special{pa 2867 1566}%
\special{pa 2901 1566}%
\special{pa 2935 1567}%
\special{pa 2969 1569}%
\special{pa 3000 1570}%
\special{sp}%
%
\special{pn 20}%
\special{pa 1810 1630}%
\special{pa 1844 1630}%
\special{pa 1878 1630}%
\special{pa 1911 1629}%
\special{pa 1945 1628}%
\special{pa 1978 1627}%
\special{pa 2011 1624}%
\special{pa 2044 1621}%
\special{pa 2076 1618}%
\special{pa 2108 1613}%
\special{pa 2139 1607}%
\special{pa 2170 1599}%
\special{pa 2200 1590}%
\special{pa 2230 1580}%
\special{pa 2258 1568}%
\special{pa 2286 1554}%
\special{pa 2313 1538}%
\special{pa 2339 1520}%
\special{pa 2365 1501}%
\special{pa 2390 1480}%
\special{pa 2415 1459}%
\special{pa 2440 1439}%
\special{pa 2465 1419}%
\special{pa 2490 1399}%
\special{pa 2517 1382}%
\special{pa 2544 1367}%
\special{pa 2572 1354}%
\special{pa 2601 1342}%
\special{pa 2630 1333}%
\special{pa 2660 1325}%
\special{pa 2691 1319}%
\special{pa 2723 1314}%
\special{pa 2755 1310}%
\special{pa 2787 1308}%
\special{pa 2820 1306}%
\special{pa 2854 1306}%
\special{pa 2887 1306}%
\special{pa 2921 1306}%
\special{pa 2955 1307}%
\special{pa 2989 1309}%
\special{pa 3020 1310}%
\special{sp}%
%
\special{pn 20}%
\special{pa 600 410}%
\special{pa 3820 410}%
\special{pa 3820 2220}%
\special{pa 600 2220}%
\special{pa 600 410}%
\special{fp}%
%
\special{pn 20}%
\special{pa 1800 1360}%
\special{pa 1834 1360}%
\special{pa 1868 1360}%
\special{pa 1901 1359}%
\special{pa 1935 1358}%
\special{pa 1968 1357}%
\special{pa 2001 1354}%
\special{pa 2034 1351}%
\special{pa 2066 1348}%
\special{pa 2098 1343}%
\special{pa 2129 1337}%
\special{pa 2160 1329}%
\special{pa 2190 1320}%
\special{pa 2220 1310}%
\special{pa 2248 1298}%
\special{pa 2276 1284}%
\special{pa 2303 1268}%
\special{pa 2329 1250}%
\special{pa 2355 1231}%
\special{pa 2380 1210}%
\special{pa 2405 1189}%
\special{pa 2430 1169}%
\special{pa 2455 1149}%
\special{pa 2480 1129}%
\special{pa 2507 1112}%
\special{pa 2534 1097}%
\special{pa 2562 1084}%
\special{pa 2591 1072}%
\special{pa 2620 1063}%
\special{pa 2650 1055}%
\special{pa 2681 1049}%
\special{pa 2713 1044}%
\special{pa 2745 1040}%
\special{pa 2777 1038}%
\special{pa 2810 1036}%
\special{pa 2844 1036}%
\special{pa 2877 1036}%
\special{pa 2911 1036}%
\special{pa 2945 1037}%
\special{pa 2979 1039}%
\special{pa 3010 1040}%
\special{sp}%
%
\special{pn 20}%
\special{pa 2780 1570}%
\special{pa 2880 1570}%
\special{fp}%
\special{sh 1}%
\special{pa 2880 1570}%
\special{pa 2813 1550}%
\special{pa 2827 1570}%
\special{pa 2813 1590}%
\special{pa 2880 1570}%
\special{fp}%
%
\special{pn 20}%
\special{pa 2850 1030}%
\special{pa 2970 1040}%
\special{fp}%
\special{sh 1}%
\special{pa 2970 1040}%
\special{pa 2905 1015}%
\special{pa 2917 1036}%
\special{pa 2902 1054}%
\special{pa 2970 1040}%
\special{fp}%
\put(17.5000,-33.7000){\makebox(0,0)[lb]{$t$}}%
\put(29.1000,-33.8000){\makebox(0,0)[lb]{$t'$}}%
\put(15.1000,-19.7000){\makebox(0,0)[lb]{$p$}}%
\put(30.8000,-18.0000){\makebox(0,0)[lb]{$p'$}}%
\put(23.8000,-19.9000){\makebox(0,0)[lb]{$\tilde{\ga}_p$}}%
%
\special{pn 20}%
\special{pa 1820 2930}%
\special{pa 2910 2930}%
\special{fp}%
\special{sh 1}%
\special{pa 2910 2930}%
\special{pa 2843 2910}%
\special{pa 2857 2930}%
\special{pa 2843 2950}%
\special{pa 2910 2930}%
\special{fp}%
\put(19.6000,-12.5000){\makebox(0,0)[lb]{$\F$}}%
%
\special{pn 8}%
\special{pa 1930 710}%
\special{pa 1840 710}%
\special{fp}%
\special{sh 1}%
\special{pa 1840 710}%
\special{pa 1907 730}%
\special{pa 1893 710}%
\special{pa 1907 690}%
\special{pa 1840 710}%
\special{fp}%
%
\special{pn 20}%
\special{sh 0.600}%
\special{ar 1810 710 81 81  0.0000000 6.2831853}%
%
\special{pn 20}%
\special{pa 3010 1590}%
\special{pa 3045 1587}%
\special{pa 3079 1583}%
\special{pa 3113 1580}%
\special{pa 3147 1575}%
\special{pa 3181 1570}%
\special{pa 3213 1564}%
\special{pa 3245 1556}%
\special{pa 3276 1548}%
\special{pa 3306 1537}%
\special{pa 3335 1525}%
\special{pa 3362 1511}%
\special{pa 3388 1494}%
\special{pa 3413 1476}%
\special{pa 3436 1454}%
\special{pa 3457 1430}%
\special{pa 3476 1404}%
\special{pa 3493 1375}%
\special{pa 3508 1345}%
\special{pa 3521 1314}%
\special{pa 3532 1281}%
\special{pa 3540 1248}%
\special{pa 3546 1214}%
\special{pa 3550 1180}%
\special{pa 3551 1147}%
\special{pa 3549 1114}%
\special{pa 3545 1082}%
\special{pa 3538 1051}%
\special{pa 3527 1022}%
\special{pa 3515 993}%
\special{pa 3499 967}%
\special{pa 3482 941}%
\special{pa 3462 917}%
\special{pa 3440 895}%
\special{pa 3416 873}%
\special{pa 3390 853}%
\special{pa 3363 833}%
\special{pa 3335 815}%
\special{pa 3305 799}%
\special{pa 3274 783}%
\special{pa 3242 768}%
\special{pa 3209 754}%
\special{pa 3176 741}%
\special{pa 3142 730}%
\special{pa 3108 719}%
\special{pa 3074 709}%
\special{pa 3040 700}%
\special{pa 3005 691}%
\special{pa 2972 684}%
\special{pa 2938 677}%
\special{pa 2905 671}%
\special{pa 2872 666}%
\special{pa 2840 661}%
\special{pa 2807 657}%
\special{pa 2775 654}%
\special{pa 2743 651}%
\special{pa 2711 649}%
\special{pa 2680 647}%
\special{pa 2648 646}%
\special{pa 2617 645}%
\special{pa 2585 644}%
\special{pa 2554 644}%
\special{pa 2523 645}%
\special{pa 2492 645}%
\special{pa 2461 646}%
\special{pa 2429 648}%
\special{pa 2398 649}%
\special{pa 2367 651}%
\special{pa 2335 653}%
\special{pa 2304 655}%
\special{pa 2272 657}%
\special{pa 2240 659}%
\special{pa 2209 661}%
\special{pa 2177 664}%
\special{pa 2145 667}%
\special{pa 2113 669}%
\special{pa 2081 672}%
\special{pa 2049 675}%
\special{pa 2017 678}%
\special{pa 1985 681}%
\special{pa 1953 684}%
\special{pa 1921 688}%
\special{pa 1889 691}%
\special{pa 1857 694}%
\special{pa 1825 697}%
\special{pa 1800 700}%
\special{sp 0.070}%
%
\special{pn 20}%
\special{pa 2040 680}%
\special{pa 1920 680}%
\special{da 0.070}%
\special{sh 1}%
\special{pa 1920 680}%
\special{pa 1987 700}%
\special{pa 1973 680}%
\special{pa 1987 660}%
\special{pa 1920 680}%
\special{fp}%
%
\special{pn 20}%
\special{pa 3000 2930}%
\special{pa 3490 2930}%
\special{da 0.070}%
%
\special{pn 20}%
\special{pa 3480 3100}%
\special{pa 1880 3100}%
\special{da 0.070}%
\special{sh 1}%
\special{pa 1880 3100}%
\special{pa 1947 3120}%
\special{pa 1933 3100}%
\special{pa 1947 3080}%
\special{pa 1880 3100}%
\special{fp}%
%
\special{pn 20}%
\special{pa 3480 2930}%
\special{pa 3511 2956}%
\special{pa 3538 2981}%
\special{pa 3556 3006}%
\special{pa 3560 3030}%
\special{pa 3548 3053}%
\special{pa 3524 3076}%
\special{pa 3493 3098}%
\special{pa 3490 3100}%
\special{sp 0.070}%
\put(23.8000,-33.6000){\makebox(0,0)[lb]{$\ell$}}%
\put(23.8000,-9.0000){\makebox(0,0)[lb]{$\tilde{\ell}_p$}}%
\put(15.0000,-7.9000){\makebox(0,0)[lb]{$p''$}}%
\put(7.4000,-11.4000){\makebox(0,0)[lb]{Poincar\'e}}%
\put(7.1000,-13.9000){\makebox(0,0)[lb]{return map}}%
%
\special{pn 13}%
\special{pa 2090 270}%
\special{pa 2750 270}%
\special{fp}%
\special{sh 1}%
\special{pa 2750 270}%
\special{pa 2683 250}%
\special{pa 2697 270}%
\special{pa 2683 290}%
\special{pa 2750 270}%
\special{fp}%
\put(22.9000,-2.0000){\makebox(0,0)[lb]{$\ga_*$}}%
\end{picture}%
\end{center}
\caption{Dynamical system with Painlev\'e property}
\label{fig:lift} 
\end{figure}
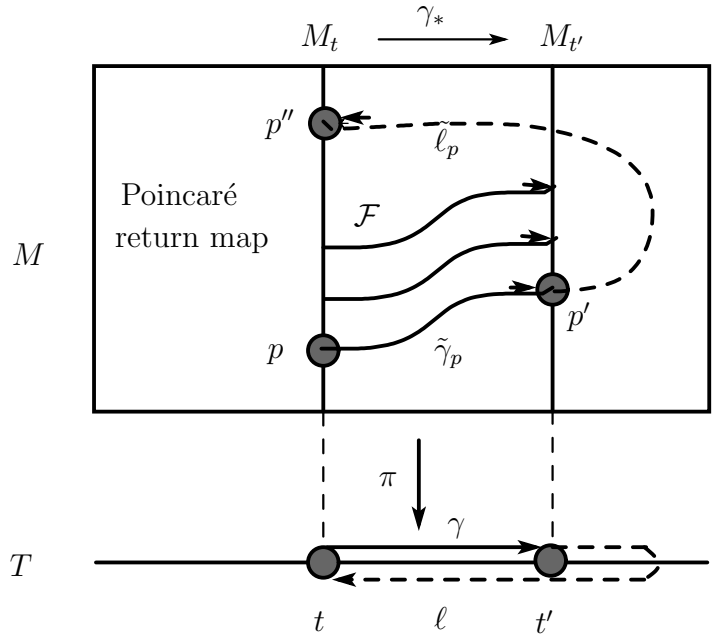       
\begin{definition}[Geometric Painlev\'e Property] 
\label{def:Painleve} 
A dynamical system $(M,\F)$ is said to have  
{\it geometric Painlev\'e property} (GPP) if for any path 
$\ga$ in $T$ and any point $p \in M_t$, where $t$ is the 
initial point of $\ga$, there exists a unique $\F$-horizontal 
lift $\tilde{\ga}_p$ of $\ga$ with initial point at $p$ (see 
Figure \ref{fig:lift}).  
Here a curve in $M$ is said to be $\F$-horizontal if it 
lies in a leaf of $\F$. 
\end{definition} 
\begin{remark}[Uniqueness] \label{rem:unique} 
Under the {\it transversality} assumption, the lifting 
problem is reduced to solving a Cauchy problem for a 
regular ODE along the curve $\ga$. 
Hence the local existence and uniquness of the lift 
$\tilde{\ga}_p$ always hold, due to the classical Cauchy 
theorem on ODE's. 
The question in Definition \ref{def:Painleve} is the 
existence of the {\it global} lift $\tilde{\ga}_p$ for 
any path $\ga$ in $T$.      
\end{remark}  
\begin{definition}[Poincar\'e Return Map] \label{def:PRM} 
If $(M,\F)$ has geometric Painlev\'e property, then any 
path $\ga$ in $T$ with initial point $t$ and terminal point 
$t'$ induces an isomorphism 
\begin{equation} \label{eqn:isom} 
\ga_* : M_t \rightarrow M_{t'}, \qquad p \mt p', 
\end{equation}
where $p'$ is the terminal point of the lift $\tilde{\ga}_p$. 
When $\ga$ is a loop in $T$ with base point at $t$, 
we have an automorphism $\ga_*$ of $M_t$, which is called   
the {\it Poincare return map} along the loop $\ga$. 
Since $\ga_*$ depends only on the homotopy class of $\ga$, 
we have the group homomorphism 
\[
\pi_1(T,t) \to \Aut\,M_t, \qquad \ga \mt \ga_*,
\]   
which we call the {\it nonlinear monodromy} of the dynamical 
system $(M,\F)$. 
\end{definition}
\begin{definition}[Hamiltonian System] 
\label{def:Hamilton}
A dynamical system $(M, \F)$ with Painlev\'e property is 
said to be {\it Hamiltonian} if $\pi : M \to T$ is a 
fibration of symplectic manifolds and the map 
(\ref{eqn:isom}) is a symplectic isomorphism for any 
path $\ga$ in $T$. 
If $(M,\F)$ is Hamiltonian, then there exists a unique closed 
$2$-form $\Om$ on $M$, called the {\it fundamental $2$-form} 
for $(M,\F)$, such that 
\begin{enumerate}
\item the form $\Om$ restricted to each fiber $M_t$ yields 
the symplectic structure $\Om_t$ on $M_t$, 
\vspace{-0.2cm}  
\item the form $\Om$ vanishes along the foliation $\F$, 
that is, for any $\F$-horizontal vector field $v$,  
\begin{equation} \label{eqn:inner}
\iota_v \Om = 0, 
\end{equation} 
where $\iota_v \Om = \Om(v, \cdot)$ stands for the inner 
product of $\Om$ by $v$. 
\end{enumerate}
\end{definition} 
\begin{remark}[Differential Equations]  \label{rem:hamiltonian}  
Condition (\ref{eqn:inner}), when expressed in terms of canonical 
local coordinates on $M$, leads to a Hamiltonian system of 
differential equations. 
\end{remark} 
\par
There are two definitions of Painlev\'e property; one is the 
geometric definition given in Definition \ref{def:Painleve}, 
which is addressed to a foliation, and the other is the 
analytic one addressed to a {\it nonlinear 
differential equation}. 
As for the latter, a differential equation is said to have 
Painlev\'e property if any solution has no movable 
singularities other than movable poles. 
This traditional but rather ambiguous definition can be 
made rigorous by the following definition.    
\begin{definition}[Analytic Painlev\'e Property] 
\label{def:Painleve2}
A nonlinear differential equation in a domain $X$ is said to have 
{\it analytic Painlev\'e property} (APP) if any meromorphic 
solution germ at any point $x \in X$ has an analytic continuation 
as a meromorphic function, along any path $\ga$ in $X$ emanating 
from $x$. 
\end{definition}
\par 
Here is a relation between geometric and analytic Painlev\'e 
properties.    
\begin{remark}[GPP Versus APP] \label{rem:relation}   
Given a dynamical system (foliation), assume that its phase space 
is an {\it algebraic variety}.  
Then, in terms of affine algebraic coordinates, the foliation is 
expressed as a differential equation and the geometric 
Painlev\'e property for the foliation implies the analytic 
Painlev\'e property for the differential equation.  
\end{remark}
\par 
In this sense the algebraicity of phase space is an 
important issue. 
Remark \ref{rem:relation} will be applied to the Hamiltonian 
system of differential equaitons in Remark \ref{rem:hamiltonian} 
and in particular to the Painlev\'e 
equation (see Theorem \ref{thm:APPHVI}). 
\subsection{Conjugacy Maps} \label{subsec:conjugacy}
As is mentioned in the Introduction, 
one of the major approaches in dynamical system theory is 
to find out a conjugacy map that converts a ``difficult'' 
dynamical system to an ``easy'' one; to extract as much 
information as possible from the latter; and to send feedback 
to obtain nontrivial (hopefully striking) results on the former. 
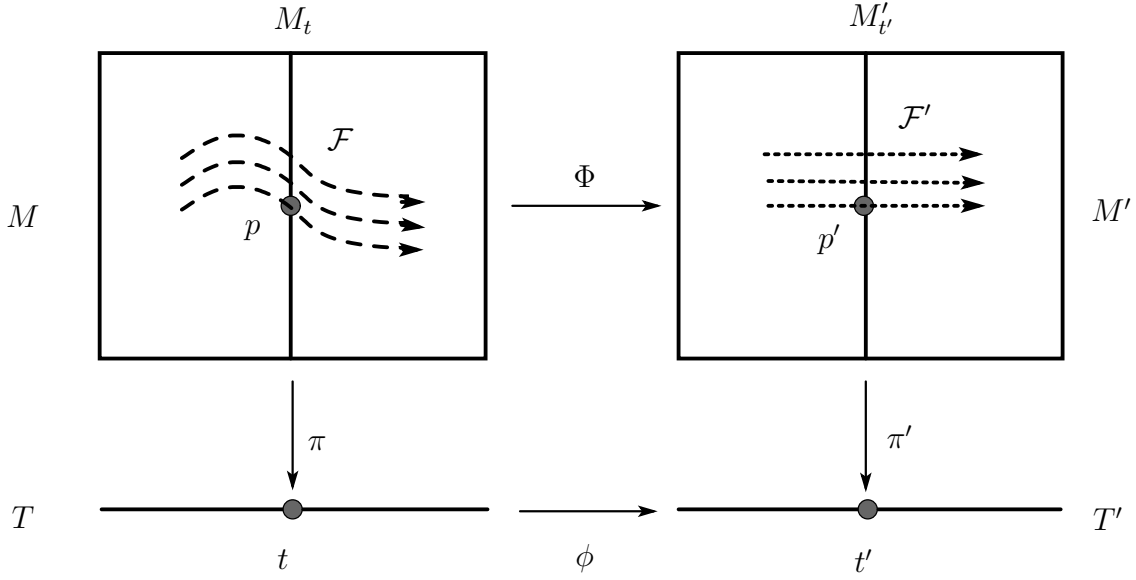
\begin{figure}[t]
\begin{center}
\unitlength 0.1in
\begin{picture}(56.90,28.20)(5.20,-33.60)
%
\special{pn 20}%
\special{pa 1010 810}%
\special{pa 3030 810}%
\special{pa 3030 2410}%
\special{pa 1010 2410}%
\special{pa 1010 810}%
\special{fp}%
%
\special{pn 20}%
\special{pa 1020 3200}%
\special{pa 3040 3200}%
\special{fp}%
%
\special{pn 13}%
\special{pa 3170 1610}%
\special{pa 3930 1610}%
\special{fp}%
\special{sh 1}%
\special{pa 3930 1610}%
\special{pa 3863 1590}%
\special{pa 3877 1610}%
\special{pa 3863 1630}%
\special{pa 3930 1610}%
\special{fp}%
%
\special{pn 20}%
\special{pa 4040 810}%
\special{pa 6050 810}%
\special{pa 6050 2410}%
\special{pa 4040 2410}%
\special{pa 4040 810}%
\special{fp}%
%
\special{pn 20}%
\special{pa 4040 3200}%
\special{pa 6040 3200}%
\special{fp}%
%
\special{pn 20}%
\special{pa 2010 810}%
\special{pa 2010 2410}%
\special{fp}%
%
\special{pn 20}%
\special{pa 5020 810}%
\special{pa 5020 2410}%
\special{fp}%
%
\special{pn 13}%
\special{pa 2020 2530}%
\special{pa 2020 3080}%
\special{fp}%
\special{sh 1}%
\special{pa 2020 3080}%
\special{pa 2040 3013}%
\special{pa 2020 3027}%
\special{pa 2000 3013}%
\special{pa 2020 3080}%
\special{fp}%
%
\special{pn 13}%
\special{pa 5020 2530}%
\special{pa 5020 3090}%
\special{fp}%
\special{sh 1}%
\special{pa 5020 3090}%
\special{pa 5040 3023}%
\special{pa 5020 3037}%
\special{pa 5000 3023}%
\special{pa 5020 3090}%
\special{fp}%
%
\special{pn 8}%
\special{sh 0.600}%
\special{ar 2010 1610 50 50  0.0000000 6.2831853}%
%
\special{pn 8}%
\special{sh 0.600}%
\special{ar 2020 3200 50 50  0.0000000 6.2831853}%
%
\special{pn 8}%
\special{sh 0.600}%
\special{ar 5010 1610 50 50  0.0000000 6.2831853}%
%
\special{pn 8}%
\special{sh 0.600}%
\special{ar 5030 3200 50 50  0.0000000 6.2831853}%
\put(19.2000,-7.1000){\makebox(0,0)[lb]{$M_t$}}%
\put(19.4000,-35.1000){\makebox(0,0)[lb]{$t$}}%
\put(34.9000,-15.1000){\makebox(0,0)[lb]{$\Phi$}}%
\put(5.2000,-17.2000){\makebox(0,0)[lb]{$M$}}%
\put(62.0000,-16.9000){\makebox(0,0)[lb]{$M'$}}%
\put(49.3000,-7.1000){\makebox(0,0)[lb]{$M'_{t'}$}}%
\put(5.5000,-33.0000){\makebox(0,0)[lb]{$T$}}%
\put(62.1000,-33.1000){\makebox(0,0)[lb]{$T'$}}%
\put(21.0000,-28.9000){\makebox(0,0)[lb]{$\pi$}}%
\put(51.3000,-28.8000){\makebox(0,0)[lb]{$\pi'$}}%
\put(17.7000,-17.9000){\makebox(0,0)[lb]{$p$}}%
\put(47.7000,-18.6000){\makebox(0,0)[lb]{$p'$}}%
\put(49.6000,-35.3000){\makebox(0,0)[lb]{$t'$}}%
%
\special{pn 20}%
\special{pa 4510 1610}%
\special{pa 5610 1610}%
\special{dt 0.054}%
\special{sh 1}%
\special{pa 5610 1610}%
\special{pa 5543 1590}%
\special{pa 5557 1610}%
\special{pa 5543 1630}%
\special{pa 5610 1610}%
\special{fp}%
%
\special{pn 20}%
\special{pa 4510 1480}%
\special{pa 5610 1490}%
\special{dt 0.054}%
\special{sh 1}%
\special{pa 5610 1490}%
\special{pa 5544 1469}%
\special{pa 5557 1490}%
\special{pa 5543 1509}%
\special{pa 5610 1490}%
\special{fp}%
\special{pa 4490 1340}%
\special{pa 5600 1340}%
\special{dt 0.054}%
\special{sh 1}%
\special{pa 5600 1340}%
\special{pa 5533 1320}%
\special{pa 5547 1340}%
\special{pa 5533 1360}%
\special{pa 5600 1340}%
\special{fp}%
%
\special{pn 20}%
\special{pa 1440 1630}%
\special{pa 1468 1611}%
\special{pa 1496 1592}%
\special{pa 1524 1574}%
\special{pa 1552 1557}%
\special{pa 1581 1542}%
\special{pa 1610 1528}%
\special{pa 1639 1517}%
\special{pa 1669 1508}%
\special{pa 1700 1502}%
\special{pa 1731 1500}%
\special{pa 1763 1500}%
\special{pa 1795 1505}%
\special{pa 1828 1512}%
\special{pa 1860 1521}%
\special{pa 1891 1533}%
\special{pa 1921 1547}%
\special{pa 1951 1563}%
\special{pa 1979 1580}%
\special{pa 2005 1598}%
\special{pa 2029 1617}%
\special{pa 2051 1637}%
\special{pa 2071 1657}%
\special{pa 2092 1677}%
\special{pa 2112 1696}%
\special{pa 2134 1716}%
\special{pa 2158 1734}%
\special{pa 2185 1752}%
\special{pa 2216 1769}%
\special{pa 2251 1784}%
\special{pa 2291 1797}%
\special{pa 2337 1809}%
\special{pa 2387 1819}%
\special{pa 2440 1827}%
\special{pa 2492 1834}%
\special{pa 2542 1839}%
\special{pa 2586 1843}%
\special{pa 2624 1846}%
\special{pa 2652 1848}%
\special{pa 2668 1849}%
\special{pa 2671 1850}%
\special{pa 2657 1850}%
\special{pa 2630 1850}%
\special{sp 0.070}%
%
\special{pn 20}%
\special{pa 2600 1840}%
\special{pa 2670 1840}%
\special{fp}%
\special{sh 1}%
\special{pa 2670 1840}%
\special{pa 2603 1820}%
\special{pa 2617 1840}%
\special{pa 2603 1860}%
\special{pa 2670 1840}%
\special{fp}%
%
\special{pn 20}%
\special{pa 1440 1500}%
\special{pa 1468 1481}%
\special{pa 1496 1462}%
\special{pa 1524 1444}%
\special{pa 1552 1427}%
\special{pa 1581 1412}%
\special{pa 1610 1398}%
\special{pa 1639 1387}%
\special{pa 1669 1378}%
\special{pa 1700 1372}%
\special{pa 1731 1370}%
\special{pa 1763 1370}%
\special{pa 1795 1375}%
\special{pa 1828 1382}%
\special{pa 1860 1391}%
\special{pa 1891 1403}%
\special{pa 1921 1417}%
\special{pa 1951 1433}%
\special{pa 1979 1450}%
\special{pa 2005 1468}%
\special{pa 2029 1487}%
\special{pa 2051 1507}%
\special{pa 2071 1527}%
\special{pa 2092 1547}%
\special{pa 2112 1566}%
\special{pa 2134 1586}%
\special{pa 2158 1604}%
\special{pa 2185 1622}%
\special{pa 2216 1639}%
\special{pa 2251 1654}%
\special{pa 2291 1667}%
\special{pa 2337 1679}%
\special{pa 2387 1689}%
\special{pa 2440 1697}%
\special{pa 2492 1704}%
\special{pa 2542 1709}%
\special{pa 2586 1713}%
\special{pa 2624 1716}%
\special{pa 2652 1718}%
\special{pa 2668 1719}%
\special{pa 2671 1720}%
\special{pa 2657 1720}%
\special{pa 2630 1720}%
\special{sp 0.070}%
%
\special{pn 20}%
\special{pa 1440 1360}%
\special{pa 1468 1341}%
\special{pa 1496 1322}%
\special{pa 1524 1304}%
\special{pa 1552 1287}%
\special{pa 1581 1272}%
\special{pa 1610 1258}%
\special{pa 1639 1247}%
\special{pa 1669 1238}%
\special{pa 1700 1232}%
\special{pa 1731 1230}%
\special{pa 1763 1230}%
\special{pa 1795 1235}%
\special{pa 1828 1242}%
\special{pa 1860 1251}%
\special{pa 1891 1263}%
\special{pa 1921 1277}%
\special{pa 1951 1293}%
\special{pa 1979 1310}%
\special{pa 2005 1328}%
\special{pa 2029 1347}%
\special{pa 2051 1367}%
\special{pa 2071 1387}%
\special{pa 2092 1407}%
\special{pa 2112 1426}%
\special{pa 2134 1446}%
\special{pa 2158 1464}%
\special{pa 2185 1482}%
\special{pa 2216 1499}%
\special{pa 2251 1514}%
\special{pa 2291 1527}%
\special{pa 2337 1539}%
\special{pa 2387 1549}%
\special{pa 2440 1557}%
\special{pa 2492 1564}%
\special{pa 2542 1569}%
\special{pa 2586 1573}%
\special{pa 2624 1576}%
\special{pa 2652 1578}%
\special{pa 2668 1579}%
\special{pa 2671 1580}%
\special{pa 2657 1580}%
\special{pa 2630 1580}%
\special{sp 0.070}%
%
\special{pn 20}%
\special{pa 2590 1710}%
\special{pa 2680 1720}%
\special{fp}%
\special{sh 1}%
\special{pa 2680 1720}%
\special{pa 2616 1693}%
\special{pa 2627 1714}%
\special{pa 2612 1733}%
\special{pa 2680 1720}%
\special{fp}%
%
\special{pn 20}%
\special{pa 2590 1590}%
\special{pa 2680 1590}%
\special{fp}%
\special{sh 1}%
\special{pa 2680 1590}%
\special{pa 2613 1570}%
\special{pa 2627 1590}%
\special{pa 2613 1610}%
\special{pa 2680 1590}%
\special{fp}%
%
\special{pn 13}%
\special{pa 3210 3210}%
\special{pa 3910 3210}%
\special{fp}%
\special{sh 1}%
\special{pa 3910 3210}%
\special{pa 3843 3190}%
\special{pa 3857 3210}%
\special{pa 3843 3230}%
\special{pa 3910 3210}%
\special{fp}%
\put(35.0000,-35.3000){\makebox(0,0)[lb]{$\phi$}}%
\put(22.0000,-13.2000){\makebox(0,0)[lb]{$\F$}}%
\put(51.9000,-12.0000){\makebox(0,0)[lb]{$\F'$}}%
\end{picture}%
\end{center}
\caption{Conjugacy map}
\label{fig:conjugacy} 
\end{figure}
\begin{definition}[Conjugacy] \label{def:conjugacy} 
A {\it conjugacy map} $\Phi$ between two dynamical systems 
$(M,\F)$ and $(M',\F')$ is a commutative diagram as in 
Figure \ref{fig:conjugacy} such that 
\begin{enumerate} 
\item the map $\Phi$ is an isomorphism that preserves geometric 
structures under consideration, e.g., measure, topology,
analytic structure, Hamiltonian structure, etc.,  
\item the foliation $\F$ is the pull-back of $\F'$ by $\Phi$, 
that is, $\Phi^* \F' = \F$.  
\end{enumerate} 
\end{definition}
\par 
It is expected that a good conjugacy map should be highly 
transcendental, reflecting the difficulty of the ``difficult'' 
dynamical system. 
This transcendental nature would make the problem not so 
tractable but at the same time so attractive. 
A few examples of conjugacy maps are presented in 
Table \ref{tab:conjugacy}, with some explanations below. 
\begin{example}[Examples of Conjugacy] \label{ex:conjugacy} 
$\phan{a}$ 
\begin{enumerate}
\item The KdV flow is conjugated to an isospectral flow by the 
scattering map, whose inversion is the Gel'fand-Levitan-Marchenko 
procedure; a seminal discovery by Gardner, Green, Kruskal and 
Miura \cite{GGKM} which opened up the soliton theory. 
\item The quadratic dynamics $P_c(z) = z^2 + c$ on $\C-K_c$  
is conjugated to the standard angle-doubling $P_0(z) = z^2$ 
on $\C-\bar{\D}$ by the B\"ottcher function, where $K_c$ is the 
filled Julia set of $P_c$ and $\bar{\D}$ is the closed unit disk. 
Using this fact, Douady and Hubbard \cite{DH} made deep studies 
on Julia sets and the Mandelbrot set for the quadratic dynamics 
on $\C$. 
\item The Painlev\'e flow on a moduli space of stable parabolic 
connections is conjugated to an isomonodromic 
flow on a moduli space of monodromy representations 
by a Riemann-Hilbert correspondence. 
\end{enumerate}
\end{example}
\begin{table}[t]
\begin{center}
\begin{tabular}{|c|l|l|l|}
\hline 
 & ``difficult" dynamics & ``easy" dynamics & conjugacy map \\ 
\hline 
1 & KdV flow & isospectral flow & scattering map \\ 
\hline 
2 & quadratic dynamics & angle-doubling map & 
B\"ottcher function \\ 
\hline 
3 & Painlev\'e flow & isomonodromic flow & 
Riemann-Hilbert map \\ 
\hline 
\end{tabular}
\end{center}
\caption{Examples of conjugacy maps} \label{tab:conjugacy} 
\end{table} 
\par 
The third example is exactly what is focused on in this article. 
As is mentioned in the Introduction, our approach is closely 
related with the isomonodromic deformation theory. 
But the latter has been mainly concerned with infinitesimal 
deformations of linear differential equations, 
without paying serious attensions to the global structure of 
Riemann-Hilbert correspondence, $\RH$, especially to its 
target space, moduli space of monodromy representations. 
Let us repeat to say that our objective is to set up the source 
and target spaces of $\RH$ firmly; to establish a  precise 
correspondence between them via $\RH$; to interrelate two 
dynamics on both sides; and to understand the dynamics of $\PVI$ 
through all these procedures. 
A similar situation seems to have occurred with KdV: 
While the machinary of inverse scattering method had 
been known since 1967 (\cite{GGKM}), it was not so soon that 
the precise correspondence was established between a class of 
potentials (of Schr\"odinger equations) and a class of scattering data, 
as in Deift and Trubowitz \cite{DT}. 
\par 
According to Definition \ref{def:conjugacy}, a conjugacy map 
$\Phi$ must be an isomorphism, namely, a {\it biholomorphism} 
in the case of holomorphic dynamics. 
However we sometimes come across such cases where this 
condition is too restrictive, that is, where the injectivity 
of $\Phi$ slightly fails to hold.   
To cover those cases, we make the following definition.  
\begin{definition}[Quasi-Conjugacy] \label{def:quasiconjugacy}
A {\it quasi-conjugacy map} $\Phi$ between two dynamical systems 
$(M,\F)$ and $(M',\F')$ is a commutative diagram as in 
Figure \ref{fig:conjugacy} such that 
\begin{enumerate} 
\item the map $\Phi$ is a surjective, proper, holomorphic map, 
\item there exists an $\F'$-invariant,   
analytic-Zariski closed, proper subset $Z' \subset M'$ such that 
$\Phi : M - Z \to M' - Z'$ is a biholomophism that 
preserves geometric structures under consideration, 
where $Z = \Phi^{-1}(Z')$,  
\item the foliation $\F$ restricted to $M - Z$ is the 
pull-back of $\F'$ on $M' - Z'$, 
\item $\Phi$ maps $\F$-trajectories in $Z$ to 
$\F'$-trajectories in $Z'$.    
\end{enumerate} 
Here the time-correspondence  $\phi : T \to T'$  is assumed 
to be biholomorphic.    
\end{definition} 
\par 
The properness condition on the map $\Phi$ in 
Definition \ref{def:quasiconjugacy} has a significant meaning 
for the geometric Painlev\'e property. 
Indeed the following lemma follows from the properness of 
$\Phi$.    
\begin{lemma}[Properness and GPP] \label{lem:proper} 
Let $\Phi : (M,\F) \to (M',\F')$ be a quasi-conjugacy map. 
Assume that the target dynamics $(M',\F')$ has geometric 
Painlev\'e property, then so does the source dynamics 
$(M,\F)$. 
\end{lemma}
{\it Proof}. 
Let $\ga$ be any compact path in $T$ with initial point $t$, 
and let $p$ be any point on $M_t$. 
By GPP for $(M',\F')$, the 
path $\phi(\ga)$ in $T'$ is lifted to 
$K' = \widetilde{\phi(\ga)}_{\Phi(p)}$ emanating 
from $\Phi(p) \in M'_{\phi(t)}$.    
Since $K'$ is compact, the properness of $\Phi$ implies that 
$K = \Phi^{-1}(K')$ is also compact. 
By conditions (3) and (4) the Cauchy problem for constructing 
the lift $\tilde{\ga}_p$ has a solution within $K$.  
Since $K$ is compact, the solution $\tilde{\ga}_p$ exists 
globally over $\ga$. 
Hence GPP for $(M,\F)$ follows. \hfill $\Box$
\par\medskip\noindent 
Lemma \ref{lem:proper} is a guiding principle in 
establishing Painlev\'e property by the conjugacy 
method: The GPP for a difficult dynamical system follows 
from that for an easy one. 
\par  
It will turn out that the Riemann-Hilbert correspondence 
is a conjugacy map (in the strict sense) for generic values 
of $\k \in \K$, but is only a quasi-conjugacy map for 
nongeneric $\k \in \K$, due to the presence of what we call 
the Riccati locus (see Definition \ref{def:riccati}).  
This locus carries classical trajectories that can be 
linearlized in terms of Gauss hypergeometric equations.         
\subsection{Application to Painlev\'e Equation}  
\label{sec:fumdamental}
\begin{figure}[t] 
\unitlength 0.1in
\begin{picture}(58.10,28.10)(1.90,-33.10)
%
\special{pn 20}%
\special{pa 800 780}%
\special{pa 2820 780}%
\special{pa 2820 2380}%
\special{pa 800 2380}%
\special{pa 800 780}%
\special{fp}%
%
\special{pn 20}%
\special{pa 810 3170}%
\special{pa 2830 3170}%
\special{fp}%
%
\special{pn 13}%
\special{pa 2960 1580}%
\special{pa 3720 1580}%
\special{fp}%
\special{sh 1}%
\special{pa 3720 1580}%
\special{pa 3653 1560}%
\special{pa 3667 1580}%
\special{pa 3653 1600}%
\special{pa 3720 1580}%
\special{fp}%
%
\special{pn 20}%
\special{pa 3830 780}%
\special{pa 5840 780}%
\special{pa 5840 2380}%
\special{pa 3830 2380}%
\special{pa 3830 780}%
\special{fp}%
%
\special{pn 20}%
\special{pa 3830 3170}%
\special{pa 5830 3170}%
\special{fp}%
%
\special{pn 20}%
\special{pa 1800 780}%
\special{pa 1800 2030}%
\special{fp}%
%
\special{pn 13}%
\special{pa 1810 2500}%
\special{pa 1810 3050}%
\special{fp}%
\special{sh 1}%
\special{pa 1810 3050}%
\special{pa 1830 2983}%
\special{pa 1810 2997}%
\special{pa 1790 2983}%
\special{pa 1810 3050}%
\special{fp}%
%
\special{pn 13}%
\special{pa 4810 2500}%
\special{pa 4810 3060}%
\special{fp}%
\special{sh 1}%
\special{pa 4810 3060}%
\special{pa 4830 2993}%
\special{pa 4810 3007}%
\special{pa 4790 2993}%
\special{pa 4810 3060}%
\special{fp}%
%
\special{pn 8}%
\special{sh 0.600}%
\special{ar 1800 1580 50 50  0.0000000 6.2831853}%
%
\special{pn 8}%
\special{sh 0.600}%
\special{ar 1810 3170 50 50  0.0000000 6.2831853}%
%
\special{pn 8}%
\special{sh 0.600}%
\special{ar 4800 1580 50 50  0.0000000 6.2831853}%
%
\special{pn 8}%
\special{sh 0.600}%
\special{ar 4820 3170 50 50  0.0000000 6.2831853}%
\put(14.6000,-6.8000){\makebox(0,0)[lb]{$\M_t(\k)$}}%
\put(17.5000,-34.8000){\makebox(0,0)[lb]{$t$}}%
\put(31.5000,-14.9000){\makebox(0,0)[lb]{$\RH_{\k}$}}%
\put(1.9000,-16.9000){\makebox(0,0)[lb]{$\M(\k)$}}%
\put(59.9000,-16.6000){\makebox(0,0)[lb]{$\R(a)$}}%
\put(45.9000,-6.7000){\makebox(0,0)[lb]{$\R_t(a) \simeq \mathcal{S}(\th)$}}%
\put(3.9000,-32.4000){\makebox(0,0)[lb]{$T$}}%
\put(60.0000,-32.4000){\makebox(0,0)[lb]{$T$}}%
\put(18.9000,-28.6000){\makebox(0,0)[lb]{$\pi_{\k}$}}%
\put(49.2000,-28.5000){\makebox(0,0)[lb]{$\pi_a$}}%
\put(14.8000,-17.9000){\makebox(0,0)[lb]{$Q$}}%
\put(45.3000,-18.0000){\makebox(0,0)[lb]{$\rho$}}%
\put(47.8000,-34.8000){\makebox(0,0)[lb]{$t$}}%
%
\special{pn 20}%
\special{pa 4300 1580}%
\special{pa 5400 1580}%
\special{dt 0.054}%
\special{sh 1}%
\special{pa 5400 1580}%
\special{pa 5333 1560}%
\special{pa 5347 1580}%
\special{pa 5333 1600}%
\special{pa 5400 1580}%
\special{fp}%
%
\special{pn 20}%
\special{pa 4300 1450}%
\special{pa 5400 1460}%
\special{dt 0.054}%
\special{sh 1}%
\special{pa 5400 1460}%
\special{pa 5334 1439}%
\special{pa 5347 1460}%
\special{pa 5333 1479}%
\special{pa 5400 1460}%
\special{fp}%
\special{pa 4280 1310}%
\special{pa 5390 1310}%
\special{dt 0.054}%
\special{sh 1}%
\special{pa 5390 1310}%
\special{pa 5323 1290}%
\special{pa 5337 1310}%
\special{pa 5323 1330}%
\special{pa 5390 1310}%
\special{fp}%
%
\special{pn 20}%
\special{pa 1230 1600}%
\special{pa 1258 1581}%
\special{pa 1286 1562}%
\special{pa 1314 1544}%
\special{pa 1342 1527}%
\special{pa 1371 1512}%
\special{pa 1400 1498}%
\special{pa 1429 1487}%
\special{pa 1459 1478}%
\special{pa 1490 1472}%
\special{pa 1521 1470}%
\special{pa 1553 1470}%
\special{pa 1585 1475}%
\special{pa 1618 1482}%
\special{pa 1650 1491}%
\special{pa 1681 1503}%
\special{pa 1711 1517}%
\special{pa 1741 1533}%
\special{pa 1769 1550}%
\special{pa 1795 1568}%
\special{pa 1819 1587}%
\special{pa 1841 1607}%
\special{pa 1861 1627}%
\special{pa 1882 1647}%
\special{pa 1902 1666}%
\special{pa 1924 1686}%
\special{pa 1948 1704}%
\special{pa 1975 1722}%
\special{pa 2006 1739}%
\special{pa 2041 1754}%
\special{pa 2081 1767}%
\special{pa 2127 1779}%
\special{pa 2177 1789}%
\special{pa 2230 1797}%
\special{pa 2282 1804}%
\special{pa 2332 1809}%
\special{pa 2376 1813}%
\special{pa 2414 1816}%
\special{pa 2442 1818}%
\special{pa 2458 1819}%
\special{pa 2461 1820}%
\special{pa 2447 1820}%
\special{pa 2420 1820}%
\special{sp 0.070}%
%
\special{pn 20}%
\special{pa 2390 1810}%
\special{pa 2460 1810}%
\special{fp}%
\special{sh 1}%
\special{pa 2460 1810}%
\special{pa 2393 1790}%
\special{pa 2407 1810}%
\special{pa 2393 1830}%
\special{pa 2460 1810}%
\special{fp}%
%
\special{pn 20}%
\special{pa 1230 1470}%
\special{pa 1258 1451}%
\special{pa 1286 1432}%
\special{pa 1314 1414}%
\special{pa 1342 1397}%
\special{pa 1371 1382}%
\special{pa 1400 1368}%
\special{pa 1429 1357}%
\special{pa 1459 1348}%
\special{pa 1490 1342}%
\special{pa 1521 1340}%
\special{pa 1553 1340}%
\special{pa 1585 1345}%
\special{pa 1618 1352}%
\special{pa 1650 1361}%
\special{pa 1681 1373}%
\special{pa 1711 1387}%
\special{pa 1741 1403}%
\special{pa 1769 1420}%
\special{pa 1795 1438}%
\special{pa 1819 1457}%
\special{pa 1841 1477}%
\special{pa 1861 1497}%
\special{pa 1882 1517}%
\special{pa 1902 1536}%
\special{pa 1924 1556}%
\special{pa 1948 1574}%
\special{pa 1975 1592}%
\special{pa 2006 1609}%
\special{pa 2041 1624}%
\special{pa 2081 1637}%
\special{pa 2127 1649}%
\special{pa 2177 1659}%
\special{pa 2230 1667}%
\special{pa 2282 1674}%
\special{pa 2332 1679}%
\special{pa 2376 1683}%
\special{pa 2414 1686}%
\special{pa 2442 1688}%
\special{pa 2458 1689}%
\special{pa 2461 1690}%
\special{pa 2447 1690}%
\special{pa 2420 1690}%
\special{sp 0.070}%
%
\special{pn 20}%
\special{pa 1230 1330}%
\special{pa 1258 1311}%
\special{pa 1286 1292}%
\special{pa 1314 1274}%
\special{pa 1342 1257}%
\special{pa 1371 1242}%
\special{pa 1400 1228}%
\special{pa 1429 1217}%
\special{pa 1459 1208}%
\special{pa 1490 1202}%
\special{pa 1521 1200}%
\special{pa 1553 1200}%
\special{pa 1585 1205}%
\special{pa 1618 1212}%
\special{pa 1650 1221}%
\special{pa 1681 1233}%
\special{pa 1711 1247}%
\special{pa 1741 1263}%
\special{pa 1769 1280}%
\special{pa 1795 1298}%
\special{pa 1819 1317}%
\special{pa 1841 1337}%
\special{pa 1861 1357}%
\special{pa 1882 1377}%
\special{pa 1902 1396}%
\special{pa 1924 1416}%
\special{pa 1948 1434}%
\special{pa 1975 1452}%
\special{pa 2006 1469}%
\special{pa 2041 1484}%
\special{pa 2081 1497}%
\special{pa 2127 1509}%
\special{pa 2177 1519}%
\special{pa 2230 1527}%
\special{pa 2282 1534}%
\special{pa 2332 1539}%
\special{pa 2376 1543}%
\special{pa 2414 1546}%
\special{pa 2442 1548}%
\special{pa 2458 1549}%
\special{pa 2461 1550}%
\special{pa 2447 1550}%
\special{pa 2420 1550}%
\special{sp 0.070}%
%
\special{pn 20}%
\special{pa 2380 1680}%
\special{pa 2470 1690}%
\special{fp}%
\special{sh 1}%
\special{pa 2470 1690}%
\special{pa 2406 1663}%
\special{pa 2417 1684}%
\special{pa 2402 1703}%
\special{pa 2470 1690}%
\special{fp}%
%
\special{pn 20}%
\special{pa 2380 1560}%
\special{pa 2470 1560}%
\special{fp}%
\special{sh 1}%
\special{pa 2470 1560}%
\special{pa 2403 1540}%
\special{pa 2417 1560}%
\special{pa 2403 1580}%
\special{pa 2470 1560}%
\special{fp}%
\put(13.1000,-22.0000){\makebox(0,0)[lb]{Painlev\'e flow}}%
\put(41.0000,-22.1000){\makebox(0,0)[lb]{Isomonodromic flow}}%
\put(31.8000,-32.4000){\makebox(0,0)[lb]{$=$}}%
%
\special{pn 20}%
\special{pa 4800 780}%
\special{pa 4800 2030}%
\special{fp}%
%
\special{pn 20}%
\special{pa 4800 2250}%
\special{pa 4800 2380}%
\special{fp}%
%
\special{pn 20}%
\special{pa 1800 2230}%
\special{pa 1800 2380}%
\special{fp}%
\end{picture}%
\par\vspace{0.5cm}\noindent
\caption{Guiding Diagram for $\PVI$} 
\label{fig:guiding}
\end{figure}
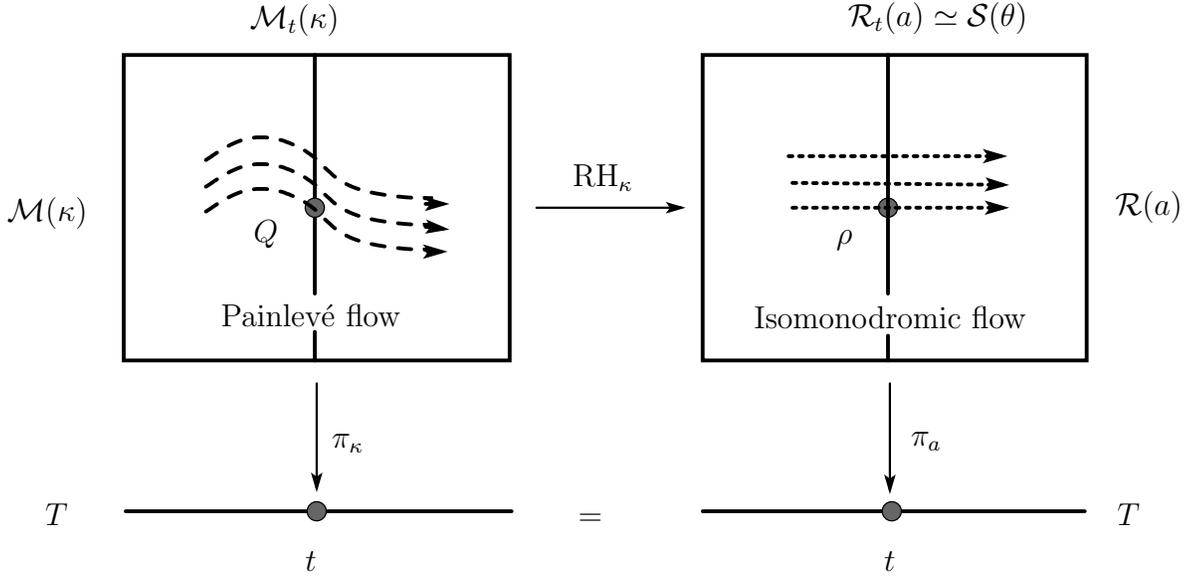

We apply the general formalism developed in the previous 
subsections to the Painlev\'e equation $\PVI(\k)$. 
In Figure \ref{fig:guiding} we present the Guiding Diagram 
that will serve as guidelines on what we shall develop in 
the sequel. 
The main ingredients of the diagram are stated as follows.
\par\vspace{0.5cm}\noindent 
{\bf Ingredients of Guiding Diagram}.   
(Figure \ref{fig:guiding})  
\par\noindent 
\begin{itemize} 
\item $T$ is the configuration space of mutually distinct 
ordered four points in $\P^1$, 
\[
T = \{\, t = (t_1,t_2,t_3,t_4) \in \P^1 \times \P^1 \times \P^1  
\times \P^1 \,:\, t_i \neq t_j \quad \mathrm{for} \quad 
i \neq j \,\}, 
\]
which serves as the space of time-variables 
(see also Remark \ref{rem:reduction}). 
\item $\M(\k)$ is the moduli space of rank-two stable 
parabolic connections on $\P^1$ having four regular singular 
points with fixed local exponents $\k \in \K$ 
(see Definitions \ref{def:parabolic1} and \ref{def:stability1}). 
It serves as the phase space of $\PVI(\k)$ as a dynamical system 
(the Painlev\'e flow). 
\item $\M_t(\k)$ is the moduli space of stable parabolic 
connections with singular points fixed at $t \in T$. 
It serves as the space of initial conditions at time $t$ for the 
Painlev\'e flow. 
\item $\R_t(a)$ is the moduli space of monodromy 
representations (up to Jordan equivalence) 
\[
\rho: \pi_1(\P^1-\{t_1,t_2,t_3,t_4\},*) \to SL_2(\C)
\]
with a fixed local monodromy data $a \in A$ (see 
Definition \ref{def:lmd}). 
It serves as the space of initial conditions at time $t$ 
for the isomonodromic flow (see Definition \ref{def:IMF}). 
\item $\R(a)$ is the disjoint union of $\R_t(a)$ over $t \in T$, 
that is, 
\begin{equation}  \label{eqn:R(a)} 
\pi_a : \R(a) = \ds \coprod_{t \in T} \R_t(a) \rightarrow T,
\end{equation}
which serves as the phase space of the isomonodromic flow 
(see Definition \ref{def:IMF}). 
\item $\RH_{\k}$ is the Riemann-Hilbert correspondence associating  
to each stable parabolic connection its monodromy representation 
(see Definition \ref{def:RH2}). 
It plays the role of a (quasi-) conjugacy map between 
the Painlev\'e flow and the isomonodromic flow. 
\item $\mathcal{S}(\th)$ is an affine cubic surface, which is 
a concrete realization of the moduli space $\R_t(a)$ of 
monodromy representations (see Theorem \ref{thm:identify}). 
\end{itemize} 
\par 
Some further explanations should be added about the objects on 
the moduli space $\M(\k)$. 
Each point $Q \in \mathcal{M}(\k)$ is representing a rank-two 
stable parabolic connection on $\P^1$ with four singular points, 
which is a refined notion of a Fuchsian system with four regular 
singular points on $\P^1$, consisting of a data 
on an algebraic vector bundle, a logarithmic connection on it,  
a prescribed determinantal structure, and a parabolic structure 
at the singular points (see Definition \ref{def:parabolic1}). 
The map $\pi_{\k} : \M(\k) \to T$ is the canonical projection 
associaitng to each connection $Q$ its ordered regular singular 
points $t = (t_1,t_2,t_3,t_4)$. 
\par 
As for the space of time variables, the following remark should 
be in order at this stage. 
\begin{remark}[Reduction of Time Variables] \label{rem:reduction} 
The original Painlev\'e equation (\ref{eqn:PVI}) has only one 
time-variable $x$, while our dynamical system has four 
time-variables $t = (t_1,t_2,t_3,t_4)$.  
The transition from $t$ to $x$ is achieved by a symplectic 
reduction that is explained as follows. 
The group of M\"obius transformations $PSL_2(\C)$ acts on $T$ 
diagonally and this action can be lifted symplectically to the 
phase space $\M(\k)$ in such a manner that the lifted action 
is commutative with the Painlev\'e flow.  
So the space of time-variables $T$ can be reduced to the 
quotient space 
\[
T/PSL_2(\C) \cong \P^1-\{0,1,\infty\}.
\]
It is well known that a natural coordinate of the quotient space 
is given by the cross ratio 
\begin{equation} \label{eqn:x} 
x = \dfrac{(t_1-t_3)(t_2-t_4)}{(t_1-t_2)(t_3-t_4)}, 
\end{equation}
which gives the independent variable of the original Painlev\'e 
equation (\ref{eqn:PVI}). 
This reduction amounts to just taking the normalization 
$t_1 = 0$, $t_2 = 1$, $t_3 = x$, $t_4 = \infty$. 
The transition from $t$ to $x$ brings slightly larger symmetry 
to the Painlev\'e equation (see Remark \ref{rem:backlund2}). 
\end{remark} 
\par 
For the most part, we shall work with four time-variables 
$t = (t_1,t_2,t_3,t_4)$, but occasionally 
we shall make use of three time-variables $t = (t_1,t_2,t_3)$ 
upon putting $t_4 = \infty$, when such a convention is more 
convenient. 
\section{Moduli Spaces of Parabolic Connections} \label{sec:phase} 
In our dynamical approach to Painlev\'e equation, 
first of all, we have to set up an appropriate phase space of 
$\PVI$ as a dynamical system.  
It is realized as the moduli space of certain stable parabolic 
connections on $\P^1$. 
Following \cite{IIS2} we shall briefly sketch its construction.  
\par 
Before entering into the subject, we should remark that there 
exist related works by Arinkin and Lysenko \cite{Arinkin, AL1}, 
who introduced moduli spaces of $SL(2)$-bundles with 
connections on $\P^1$ in the context of Painlev\'e equation. 
Unfortunately, they treated the moduli spaces mostly as stacks 
and restricted themselves to generic parameters in $\K$ to 
avoid reducible or resonant connections.  
For a total understanding of the Painlev\'e equation, however, 
the locus of nongeneric parameters often plays a significant 
part. 
In order to cover all parameters, we should take parabolic 
structures into account (see Remark \ref{rem:advantage} for 
this and for another reason).  
Moreover, in order to develop a good moduli theory in the 
framework of geometric invariant theory \cite{Mumford}, 
we need the concept of stability. 
These demands lead us to consider stable parabolic connections.             
\subsection{Parabolic Connections}  
\label{subsec:parabolic}
In what follows, a vector bundle will be identified with the 
locally free sheaf associated to it.   
For a vector bundle $E$ on $\P^1$ and a point $x \in \P^1$, 
we denote by $E|_{x}$ the fiber of $E$ over $x$ (not the stalk 
at $x$), namely we have $E|_x = E/E(-x)$ with  
$E(-x) = E \ot \O_{\P^1}(-x)$. 
\begin{definition}[Parabolic Connection] \label{def:parabolic1} 
Given $(t,\k) \in T \times \K$, a $(t,\k)$-{\it parabolic 
connection} is a quadruple $Q = (E,\nabla,\psi,l)$ such that 
the following conditions are satisfied: 
\begin{enumerate} 
\item $E$ is a rank-two vector bundle over $\P^1$. 
\item $\nabla : E \rightarrow E \ot 
\Om^1_{\P^1}(D_t)$ is a connection, where 
$D_t$ is the divisor
\[
D_t = t_1 + t_2 + t_3 + t_4.
\]
\item $\psi : \det E \to \O_{\P^1}(-t_4)$ is a horizontal 
isomorphism, where $\O_{\P^1}(-t_4) \subset \O_{\P^1}$ is 
equipped with the connection $d_{t_4}$ induced from the exterior 
differentiation $d : \O_{\P^1} \rightarrow \Om^1_{\P^1}$. 
\item $l = (l_1,l_2,l_3,l_4)$, where $l_i$ is a $1$-dimensional 
subspace of the fiber $E|_{t_i}$ over $t_i$ such that 
\[ 
(\Res_{t_i}(\nabla) - \l_i \, \id_{E|_{t_i}})|_{l_i} = 0, 
\] 
namely, $l_i$ is an eigenline of $\Res_{t_i}(\nabla)$ with 
eigenvalue $\l_i$, where $\Res_{t_i}(\nabla) \in 
\End(E|_{t_i})$ is the residue of $\nabla$ at $t_i$ and 
the parameter $\l_i$ is defined so that   
\begin{equation} \label{eqn:lambda}
\k_i = \left\{
\begin{array}{ll}
2 \l_i     \qquad & (i = 1,2,3), \\[2mm]
2 \l_4-1   \qquad & (i = 4). 
\end{array}
\right.  
\end{equation}
The data $l$ is called the {\it parabolic structure} of the 
parabolic connection $Q$.  
\end{enumerate}  
\end{definition}
\begin{table}[t]  
\begin{center} 
\begin{tabular}{|c||c|c|c|c|}
\hline 
singularity & $t_1$ & $t_2$ & $t_3$ & $t_4$ \\ 
\hline 
first exponent & $-\l_1$ & $-\l_2$ & $-\l_3$ & $-\l_4$ \\ 
\hline 
second exponent & $\l_1$ & $\l_2$  & $\l_3$ & $\l_4-1$ \\ 
\hline 
difference & $\k_1$ & $\k_2$ & $\k_3$ & $\k_4$ \\
\hline 
\end{tabular}
\end{center}   
\caption{Riemann Scheme: first exponents correspond to 
parabolic structures}
\label{tab:riemann}
\end{table} 
\begin{remark}[Riemann Scheme] \label{rem:exponent} 
An eigenvalue of $-\Res_{t_i}(\nabla)$ is called a  
{\it local exponent} of $\nabla$ at $t_i$. 
By condition (4) of Definition \ref{def:parabolic1}, 
one exponent at $t_i$ is $-\l_i$ corresponding to 
the parabolic structure $l_i$ (the first exponent).  
By condition (3) the eigenvalues of $\Res_{t_i}(\nabla)$ 
are sumed up to the residue 
$\Res_{t_i}(d_{t_4})$ of the connection $d_{t_4}$ on the line 
bundle $\O_{\P^1}(-t_4)$.  
Since 
\[
\Res_{t_i}(d_{t_4}) = \left\{\begin{array}{ll} 
0  \qquad & (i = 1,2,3), \\[2mm] 
1  \qquad & (i = 4), 
\end{array}\right.
\]
the exponents of $\nabla$ at each singular point are 
given as in Table \ref{tab:riemann}. 
Formula (\ref{eqn:lambda}) means that $\k_i$ is the 
{\it difference} of the second exponent from the first 
exponent at $t_i$. 
For brevity $\k_i$ is often referred to as the 
local exponent at $t_i$.  
We remark that the fourth singular point $t_4$ is somewhat 
distinguished from the others.   
\end{remark} 
\begin{remark}[Determinantal Structure] 
In condition (3) of Definition \ref{def:parabolic1}, the 
horizontal isomorphism $\psi : \det E \to \O_{\P^1}(-t_4)$ 
is referred to as the {\it determinantal structure} of the 
parabolic connection $Q$. 
Here the choice of $\O_{\P^1}(-t_4)$ as the target line bundle 
of $\psi$ is just for convenience. 
More generally, a determinantal structure relative to $L$ 
is conceivable for any line bundle $L$ with a connection $d_L$, 
meaning a horizontal isomorphism $\psi : \det \, E \to L$. 
\end{remark}
\par
There are at least two advantages of taking parabolic structures 
into account. 
\begin{remark}[Advantages of Parabolic Structures] 
\label{rem:advantage} $\phan{a}$ 
\begin{enumerate}
\item The Reimann-Hilbert and isomonodromic approaches become 
feasible for {\it all} parameters $\k \in \K$. 
Without parabolic structures, people usually avoid nongeneric 
parameters for ``technical" reasons, but they cannot be ruled 
out becuase many interesting phenomena occur at nongeneric 
parameters both moduli-theoretically and 
special-function-theoretically. 
Moreover, what is called the technical difficulty is in fact 
an essential difficulty. 
\item The technique of elementary transformations becomes 
available. 
Here elementary transformations are certain kinds of gauge 
transformations canonically associated to parabolic structures 
(see Definition \ref{def:elementary}). 
They, together with the powerful technique of 
Langton \cite{Langton}, play an important part in solving 
the Riemann-Hilbert problem (see Remark \ref{rem:howtosolve}). 
They also serve as some portions of the B\"acklund 
transformations (see Definition \ref{def:backlund}). 
\end{enumerate}
\end{remark}
\begin{figure}[t]
\[
\begin{CD}
  @.     @.      0    @.    0      @.            @. \\
  @.     @.           @VVV         @VVV          @. \\
0 @>>> E(-t_i) @>>> \tilde{E} @>>> l_i      @>>> 0  \\
  @.     @|           @VVV         @VVV             \\
0 @>>> E(-t_i) @>>>     E     @>>> E|_{t_i} @>>> 0  \\
  @.     @.           @VVV         @VVV             \\
  @.     @.         E|_{t_i}/l_i @= E|_{t_i}/l_i @. \\
  @.     @.           @VVV         @VVV          @. \\
  @.     @.      0    @.    0      @.            @. \\
\end{CD}
\]
\caption{Diagram related to an elementary transformation}
\label{fig:elementary}
\end{figure}
\begin{definition}[Elementary Transformation] 
\label{def:elementary}
Let $Q = (E,\nabla,\psi,l)$ be a parabolic connection with 
a determinantal structure $\psi : \det E \to L$ and a 
parabolic structure $l = (l_1,l_2,l_3,l_4)$. 
The {\it elementary transform} of $Q$ at $t_i$ is the parabolic 
connection $\tilde{Q} = 
(\tilde{E},\tilde{\nabla},\tilde{\psi},\tilde{l})$ defined in 
the following manner (see also Figure \ref{fig:elementary}). 
\begin{enumerate}
\item The bundle $\tilde{E}$ is the subsheaf  
$\tilde{E} = \Ker [\, E \to E|_{t_i}/l_i \,]$, where 
$E \to E|_{t_i}/l_i$ is the composite of the canonical 
projections  
$E \to E/E(-t_i) = E|_{t_i}$ and $E|_{t_i} \to E|_{t_i}/l_i$.
\item The connection $\tilde{\nabla} = \nabla|_{\tilde{E}} : 
\tilde{E} \to \tilde{E} \ot \Om^1_{\P^1}(D_t)$ is the 
restriction of $\nabla$ to the subsheaf $\tilde{E} \subset E$. 
This is well-defined because condition (4) 
of Definition \ref{def:parabolic1} implies that $\nabla$ maps 
$\tilde{E}$ into $\tilde{E} \ot \Om^1_{\P^1}(D_t)$.  
\item The determinantal structure $\tilde{\psi} = 
\psi|_{\det \, \tilde{E}} : \det \, \tilde{E} \to 
\tilde{L} := L(-t_i)$ is the restriction of $\psi$ to 
the subsheaf $\det \, \tilde{E} \subset \det \, E$, 
where $L(-t_i) = L \ot \O_{\P^1}(-t_i)$ is equipped 
with the connection $d_L \ot d_{t_i}$ with 
$d_{t_i}$ being the connection on $\O_{\P^1}(-t_i) 
\subset \O_{\P^1}$ induced from the exterior 
differentiation $d : \O_{\P^1} \to \Om^1_{\P^1}$. 
This is well-defined becuase one has  
$\det \, \tilde{E} = (\det \, E)(-t_i)$ and $\psi$ maps 
$(\det \, E)(-t_i)$ to $L(-t_i)$ by condition (4) 
of Definition \ref{def:parabolic1}. 
\item The parabolic structure $\tilde{l}_j$ at $t_j$ 
is defined by
\[
\tilde{l}_j = \left\{\begin{array}{cl}
E(-t_i)/\tilde{E}(-t_i) \qquad & (j = i), \\[2mm]
l_j \qquad & (j \neq i),  
\end{array}\right.
\]
where $\tilde{l}_i$ is well-defined, since 
$\tilde{l}_i = E(-t_i)/\tilde{E}(-t_i) 
\subset \tilde{E}/\tilde{E}(-t_i) = \tilde{E}|_{t_i}$. 
\end{enumerate} 
\end{definition}
\par
There are other types of elementary transformations  
defined in similar manners; see \cite{IIS2}.   
Elementary transformations were intensively studied by 
Maruyama \cite{Maruyama} and others. 
For parabolic structures appearing in various moduli 
problems, we refer to  Maruyama and Yokogawa \cite{MY}, 
Nakajima \cite{Nakajima}, Inaba \cite{Inaba} and the 
references therein. 
See also Huybrechts and Lehn \cite{HL}, 
Nitsure \cite{Nitsure} Simpson \cite{Simpson} for 
related moduli problems.   
\subsection{Stability} \label{subsec:stability1}
To obtain a good moduli space, namely, to avoid non-Hausdorff 
phenomena, we require a concept of stability for parabolic 
connections.  
\begin{definition}[Stability]  \label{def:stability1} 
A {\it weight} is a sequence of mutually distinct rational numbers
\[
\a = (\a_1,\a_1',\dots,\a_4,\a_4') \qquad 
\mbox{such that} \qquad 0 < \a_i < \a_i' < 1. 
\]
Given a weight $\a$, a parabolic connection $Q = (E,\nabla,\psi,l)$ 
is said to be $\a$-{\it stable} if for any proper subbundle  
$F \subset E$ such that $\nabla(F) \subset F \ot 
\Om^1_{\P^1}(D_t)$, one has 
\begin{equation} \label{eqn:stability1}  
\dfrac{\pardeg\, F}{\rank\, F} < \dfrac{\pardeg\, E}{\rank\, E},  
\end{equation}
where $\pardeg\, E$ and $\pardeg\, F$, called the parabolic degrees, 
are defined by 
\begin{eqnarray*}
\pardeg \, E &=& \deg E + \sum_{i=1}^4 
\left\{
\a_i \dim (E|_{t_i}/l_i) + \a_i' \dim l_i \right\}
= \deg E + \sum_{i=1}^4 (\a_i + \a_i'),  \\
\pardeg \, F &=& \deg F + \sum_{i=1}^4 
\left\{
\a_i \dim(F|_{t_i}/l_i \cap F|_{t_i}) + 
\a_i' \dim(l_i \cap F|_{t_i}), 
\right\} 
\end{eqnarray*} 
The concept of $\a$-{\it semistability} is defined in a similar 
manner by weakening the condition (\ref{eqn:stability1}) so that 
it allows equality. 
A weight $\a$ is said to be {\it generic} if every $\a$-simistable 
object is $\a$-stable. 
Hereafter the weight will be assumed to be generic. 
\end{definition} 
\subsection{Moduli Space of Stable Parabolic Connections} 
\label{subsec:moduli1} 
Based on arguments from geometric invariant theory, we can 
establish the following result \cite{IIS2}.   
\begin{theorem}[Moduli Space] \label{thm:moduli1}
Fix a generic weight $\a$.  
\begin{enumerate}
\item There exists a fine moduli scheme $\M_t(\k)$ of stable 
$(t,\k)$-parabolic connections.  
\item The moduli space $\M_t(\k)$ is a smooth, irreducible, 
quasi-projective surface. 
\item As a relative setting, there exists a family of moduli spaces 
\begin{equation}  \label{eqn:family1}
\pi : \M \rightarrow T \times \K,  
\end{equation}
such that $\pi$ is a smooth morphism whose fiber over 
$(t,\k) \in T \times \K$ is just $\M_t(\k)$. 
\item Fixing an exponent $\k \in \K$, one can also speak 
of the family 
\begin{equation} \label{eqn:family2}
\pi_{\k} : \M(\k) \rightarrow T.  
\end{equation} 
\end{enumerate} 
\end{theorem} 
\par 
We insist that the fibration (\ref{eqn:family2}) gives a precise 
phase space of $\PVI(\k)$ as a time-dependent dynamical system . 
In this regard the following remark should be in order.  
\begin{remark}[Connections on Trivial Vector Bundle]  
\label{rem:schlesinger}  
In the isomonodromic approach \linebreak to $\PVI$, 
people usually work with linear Fuchsian systems of the form    
\begin{equation} \label{eqn:schlesinger} 
\dfrac{dY}{dz} = A(z)\ Y, 
\qquad 
A(z) = \sum_{i=1}^4 \dfrac{A_i}{z-t_i},  
\end{equation}
namely, Fuchsian connections on the {\it trivial} vector 
bundle, and derive the Schlesinger system 
\begin{equation} \label{eqn:schlesinger2} 
\dfrac{\partial A_i}{\partial t_i} = 
\sum_{k \neq i} \dfrac{[A_i, A_k]}{t_k-t_i}, \qquad  
\dfrac{\partial A_i}{\partial t_j} = 
\dfrac{[A_i, A_j]}{t_j-t_i} \qquad (i \neq j),  
\end{equation}
and then recast it to the Painlev\'e equation. 
In that case they are supposing that the totality of 
the connections in (\ref{eqn:schlesinger}) forms 
a phase space of $\PVI(\k)$.        
However, it is only isomorphic to a Zariski-open proper subset 
of the true phase space, that is, our moduli space $\M(\k)$,  
and some trajectories actually escape from this open 
subset.  
Thus, with such a na\"{\i}ve setting of phase space as in 
(\ref{eqn:schlesinger}), the geometric Painlev\'e property 
is not fulfilled, (although the analytic Painlev\'e property  
for the system (\ref{eqn:schlesinger2}) holds true as was 
proved\footnote{under generic conditions on exponents} 
by Malgrange \cite{Malgrange} and Miwa \cite{Miwa}).     
This is why we had to consider connections on nontrivial  
vector bundles together with the extra data of parabolic 
structures, in order to build a complete phase space.  
In our setting, the geometric Painlev\'e property  
holds quite naturally and then the analytic Painlev\'e 
property follows from this and the algebraicity of the 
phase space 
(see Theorem \ref{thm:pp2}, Remark \ref{rem:relation} 
and Theorem \ref{thm:APPHVI}).  
\end{remark} 
\subsection{Parabolic $\mbox{\bm $\phi$}$-Connection}  
\label{subsec:parabolic2}
As the moduli space $\M_t(\k)$ is quasi-projective, 
it is natural to pose the following problem.  
\begin{problem}[Compactification] \label{prob:compactify}
Compactify the moduli space $\M_t(\k)$ in a natural manner.  
\end{problem}
\par  
This problem is settled by introducing the notion of parabolic 
$\phi$-connection, which is a generalized object of parabolic 
connections, allowing some degeneracy in the exterior 
differential part. 
This procedure reminds us of semi-classical limits of 
Schr\"odinger equations as the Plank constant tends 
to zero; we compactify the moduli space by adding 
some ``semi-classical'' objects.  
\begin{definition}[Parabolic $\mbox{\bm $\phi$}$-Connection] 
\label{def:parabolic2}
For a fixed $(t,\k) \in T \times \K$, a $(t,\k)$-{\it parabolic} 
$\phi$-{\it connection} is a sextuple 
$Q = (E_1,E_2,\phi,\nabla,\psi,l)$ such that the following 
conditions are satisfied: 
\begin{enumerate}
\item $E_1$ and $E_2$ are rank-two vector bundles over $\P^1$ 
of the same degree $\deg E_1 = \deg E_2$. 
\item $\phi : E_1 \to E_2$ is an $\O_{\P^1}$-homomorphism. 
\item $\nabla : E_1 \rightarrow E_2 \ot 
\Om^1_{\P^1}(D_t)$ is a $\C$-linear map such that  
\[
\nabla(fs) = \phi(s) \ot df + f \nabla(s) 
\qquad \mbox{for} \qquad f \in \O_{\P^1}, \, s \in E_1. 
\]
\item $\psi : \det E_2 \to \O_{\P^1}(-t_4)$ is a horizontal 
isomorphism in the sense that
\[
(\psi \ot 1)(\phi(s_1) \wedge \nabla(s_2) + 
\nabla(s_1) \wedge \phi(s_2)) = 
d_{t_4}(\psi(\phi(s_1) \wedge \phi(s_2))) 
\qquad \mbox{for} \quad s_1, s_2 \in E_1. 
\] 
\item $l = (l_1,l_2,l_3,l_4)$, where $l_i$ is a $1$-dimensional 
subspace of the fiber $E_1|_{t_i}$ over $t_i$ such that 
\[ 
(\Res_{t_i}(\nabla) - \l_i \, \phi|_{t_i})|_{l_i} = 0, 
\] 
where $\Res_{t_i}(\nabla) \in \Hom(E_1|_{t_i}, E_2|_{t_i})$ is 
the residue of $\nabla$ at $t_i$ and $\l_i$ is defined 
by (\ref{eqn:lambda}). 
\end{enumerate}
\end{definition} 
\par 
We remark that a parabolic $\phi$-connection is isomorphic 
to a parabolic connection if $\phi$ is an isomorphism, while 
it is thought of as a degenerate object if $\phi$ is not an 
isomorphism.         
\subsection{Stability} \label{subsec:stability2}
Again, to get a good moduli space, we need a concept of stability 
for parabolic $\phi$-connections. 
The following definition may be intricate at first glance, 
but works well in practice.   
\begin{definition}[Stability]  \label{def:stability2} 
A {\it weight} is a sequence $\a = (\a_1,\a_1',\dots,\a_4,\a_4')$ 
of mutually distinct rational numbers, together with positive 
integers $\b_1$, $\b_2$, $\ga$, such that
\[ 
(\b_1 + \b_2) \a_i < (\b_1 + \b_2) \a_i' < \b_1,  \qquad  
\ga \gg 0. 
\]
A $(t,\k)$-parabolic $\phi$-connection 
$Q = (E_1,E_2,\phi,\nabla,\psi,l)$  is said to be 
$(\a,\b,\ga)$-{\it stable} if for any proper 
subbundle $(F_1,F_2) \subset (E_1,E_2)$ such that 
$\phi(F_1) \subset F_2$ and 
$\nabla(F_1) \subset F_2 \ot \Om^1_{\P^1}(D_t)$, one has 
\begin{equation} \label{eqn:stability2}
\dfrac{\pardeg (F_1,F_2)}{\b_1 \rank \, F_1+ \b_2 \rank \, F_2} 
< \dfrac{\pardeg (E_1,E_2)}{\b_1 \rank\,E_1 + \b_2 \rank\,E_2},  
\end{equation}
where $\pardeg(E_1,E_2)$ and $\pardeg (F_1,F_2)$ are define by 
\begin{eqnarray*}
\pardeg (E_1,E_2) &=& \b_1 \, \deg E_1(-D_t) 
+ \b_2 \, (\deg E_2 - \ga \, \rank\,E_2) \\
&& + (\b_1+\b_2) \sum_{i=1}^4 
\left\{
\a_i \dim (E_1|_{t_i}/l_i) + \a_i' \dim l_i \right\},  \\
\pardeg (F_1,F_2) &=& \b_1 \, \deg F_1(-D_t) 
+ \b_2 \, (\deg F_2 - \ga \, \rank\,F_2) \\ 
&& + (\b_1+\b_2) \sum_{i=1}^4 
\left\{
\a_i \dim(F_1|_{t_i}/l_i \cap F_1|_{t_i}) + 
\a_i' \dim(l_i \cap F_1|_{t_i}) \right\},  
\end{eqnarray*} 
The concept of $(\a,\b,\ga)$-{\it semistability} is defined in 
a similar manner by weakening the condition (\ref{eqn:stability2}) 
so that it allows equality. 
A weight $(\a,\b,\ga)$ is said to be {\it generic} if every 
$(\a,\b,\ga)$-simistable object is $(\a,\b,\ga)$-stable. 
Hereafter the weight will be assumed to be generic. 
\end{definition} 
\subsection{Moduli Space of Stable Parabolic 
$\mbox{\bm $\phi$}$-Connections} \label{subsec:moduli2}
Again, based on arguments from geometric invariant theory, 
we have the following result \cite{IIS2}. 
\begin{theorem}[Moduli Space] \label{thm:moduli2}
Fix a generic weight $(\a,\b,\ga)$.   
\begin{enumerate} 
\item There is a coarse moduli scheme $\ol{\M}_t(\k)$ of 
stable $(t,\k)$-parabolic $\phi$-connections. 
\item The moduli space $\ol{\M}_t(\k)$ is a smooth, 
irreducible, projective surface. 
\item The moduli space $\M_t(\k)$ 
is embedded into the compactified one  
$\ol{\M}_t(\k)$ by the natural map
\[
\M_t(\k) \hookrightarrow \ol{\M}_t(\k), 
\qquad (E,\nabla,\psi,l) \mt (E,E,\id,\nabla,\psi,l), 
\]
the image of which is the open subscheme of all stable 
$(t,\k)$-parabolic $\phi$-connections 
$Q = (E_1,E_2,\phi,\nabla,\psi,l)$ such that 
$\phi : E_1 \to E_2$ is an isomorphism. 
\item As a relative setting, there exists a family of moduli 
spaces 
\[
\ol{\pi} : \ol{\M} \to T \times \K, 
\]
such that $\ol{\pi}$ is a smooth, projective morhism whose 
fiber over $(t,\k) \in T \times \K$ is just the compactified 
moduli space $\ol{\M}_t(\k)$. 
\item There exists a commutative diagram 
\[
\begin{CD}
\M @> \scriptstyle{\mathrm{embedding}} >> \ol{\M} \\
@V \pi VV @VV \ol{\pi} V \\
T \times \K @= T \times \K. 
\end{CD}
\]
\end{enumerate}    
\end{theorem} 
\subsection{Realization of Moduli Spaces} 
\label{subsec:realization}
The moduli space $\M_t(\k)$ of stable parabolic 
$(t,\k)$-connections, together with its compactification 
$\ol{\M}_t(\k)$, admits a concrete realization in terms 
of the Hirzebruch surface $\Sig_2$ of degree $2$. 
The surface $\Sig_2$ is the $\P^1$-bundle over $\P^1$ 
whose cross section at infinity, denoted by $F_0$, has 
self-intersection number $-2$. 
Moreover $\Sig_2 - F_0$ is isomorphic to the line bundle 
$\Om_{\P^1}^1(D_t)$ over $\P^1$. 
Given any $t = (t_1,t_2,t_3,t_4) \in T$ and 
$i \in \{1,2,3,4\}$, 
let $F_i$ denote the fiber over $t_i$ of the fibration 
$\Sig_2 \to \P^1$. 
Then we have the following theorem from Inaba, Iwasaki and 
Saito \cite{IIS2}. 
\begin{theorem}[Realization of Moduli Spaces] 
\label{thm:realization} 
Let $(t,\k) \in T \times \K$ be fixed. 
\begin{enumerate} 
\item  
$\ol{\M}_t(\k)$ is an $8$-point blow-up of the Hirzebruch 
surface $\Sig_2$ of degree $2$, blown up at certain 
two points on each fiber $F_i$, $i = 1,2,3,4$. 
The location of the blowing-up points, possibly infinitely 
near, is determined by the value of $\k$.  
\item $\ol{\M}_t(\k)$ has a unique effective anti-canonical 
divisor 
\[
Y_t(\k) = 2 E_0 + E_1 + E_2 + E_3 + E_4 \in 
\left|-K_{\ol{\M}_t(\k)}\right|, 
\]
where $E_i$ is the strict transform of $F_i$ for 
$i = 0,1,2,3,4$.  
Each irreducible component $E_i$ of $Y_t(\k)$ satisfies 
the condition
\[ 
K_{\ol{\M}_t(\k)} \cdot E_i = 0 \qquad (i = 0,1,2,3,4). 
\]
\item $\M_t(\k)$ is obtained from $\ol{\M}_t(\k)$ by 
removing $Y_t(\k)_{\mathrm{red}}$.    
\end{enumerate}
\end{theorem}
\begin{figure}[t]
\begin{center} 
\unitlength 0.1in
\begin{picture}(39.25,31.25)(2.20,-34.15)
%
\special{pn 20}%
\special{pa 3200 610}%
\special{pa 3200 2610}%
\special{fp}%
%
\special{pn 20}%
\special{pa 1400 600}%
\special{pa 1400 2600}%
\special{fp}%
%
\special{pn 20}%
\special{pa 2010 600}%
\special{pa 2010 2600}%
\special{fp}%
%
\special{pn 20}%
\special{pa 2600 600}%
\special{pa 2600 2590}%
\special{fp}%
%
\special{pn 20}%
\special{pa 820 600}%
\special{pa 3810 600}%
\special{pa 3810 2610}%
\special{pa 820 2610}%
\special{pa 820 600}%
\special{fp}%
%
\special{pn 20}%
\special{pa 820 3210}%
\special{pa 3810 3210}%
\special{fp}%
%
\special{pn 13}%
\special{pa 2290 2690}%
\special{pa 2290 3150}%
\special{fp}%
\special{sh 1}%
\special{pa 2290 3150}%
\special{pa 2310 3083}%
\special{pa 2290 3097}%
\special{pa 2270 3083}%
\special{pa 2290 3150}%
\special{fp}%
%
\special{pn 20}%
\special{pa 820 1020}%
\special{pa 3820 1020}%
\special{fp}%
%
\special{pn 20}%
\special{pa 1400 3160}%
\special{pa 1400 3270}%
\special{fp}%
%
\special{pn 20}%
\special{pa 3190 3150}%
\special{pa 3190 3260}%
\special{fp}%
%
\special{pn 20}%
\special{pa 2000 3160}%
\special{pa 2000 3270}%
\special{fp}%
%
\special{pn 20}%
\special{pa 2610 3160}%
\special{pa 2610 3270}%
\special{fp}%
%
\special{pn 20}%
\special{pa 1810 1210}%
\special{pa 2200 1620}%
\special{dt 0.054}%
\special{pa 2200 1620}%
\special{pa 2200 1620}%
\special{dt 0.054}%
%
\special{pn 20}%
\special{pa 2400 1200}%
\special{pa 2790 1610}%
\special{dt 0.054}%
\special{pa 2790 1610}%
\special{pa 2790 1610}%
\special{dt 0.054}%
%
\special{pn 20}%
\special{pa 3010 1220}%
\special{pa 3400 1630}%
\special{dt 0.054}%
\special{pa 3400 1630}%
\special{pa 3400 1630}%
\special{dt 0.054}%
%
\special{pn 20}%
\special{pa 1200 1210}%
\special{pa 1590 1620}%
\special{dt 0.054}%
\special{pa 1590 1620}%
\special{pa 1590 1620}%
\special{dt 0.054}%
%
\special{pn 20}%
\special{pa 1600 2010}%
\special{pa 1200 2400}%
\special{dt 0.054}%
\special{pa 1200 2400}%
\special{pa 1200 2400}%
\special{dt 0.054}%
%
\special{pn 20}%
\special{pa 2190 2010}%
\special{pa 1790 2400}%
\special{dt 0.054}%
\special{pa 1790 2400}%
\special{pa 1790 2400}%
\special{dt 0.054}%
%
\special{pn 20}%
\special{pa 2800 2000}%
\special{pa 2400 2390}%
\special{dt 0.054}%
\special{pa 2400 2390}%
\special{pa 2400 2390}%
\special{dt 0.054}%
%
\special{pn 20}%
\special{pa 3390 1990}%
\special{pa 2990 2380}%
\special{dt 0.054}%
\special{pa 2990 2380}%
\special{pa 2990 2380}%
\special{dt 0.054}%
\put(3.4000,-16.4000){\makebox(0,0)[lb]{$\varSigma_2$}}%
\put(13.9000,-3.8000){\makebox(0,0){$E_1$}}%
\put(20.3000,-3.8000){\makebox(0,0){$E_2$}}%
\put(26.1000,-3.8000){\makebox(0,0){$E_3$}}%
\put(32.1000,-3.8000){\makebox(0,0){$E_4$}}%
\put(41.3000,-10.2000){\makebox(0,0){$E_0$}}%
\put(14.0000,-34.8000){\makebox(0,0){$t_1$}}%
\put(20.0000,-34.9000){\makebox(0,0){$t_2$}}%
\put(26.1000,-34.9000){\makebox(0,0){$t_3$}}%
\put(31.9000,-35.0000){\makebox(0,0){$t_4$}}%
\put(24.5000,-28.8000){\makebox(0,0){$\pi$}}%
\put(4.9000,-32.2000){\makebox(0,0){$\P^1$}}%
\put(38.1000,-4.6000){\makebox(0,0)[lb]{$\infty$-section}}%
%
\special{pn 8}%
\special{pa 4140 520}%
\special{pa 4140 880}%
\special{fp}%
\special{sh 1}%
\special{pa 4140 880}%
\special{pa 4160 813}%
\special{pa 4140 827}%
\special{pa 4120 813}%
\special{pa 4140 880}%
\special{fp}%
%
\special{pn 20}%
\special{sh 0.600}%
\special{ar 3200 1420 35 35  0.0000000 6.2831853}%
%
\special{pn 20}%
\special{sh 0.600}%
\special{ar 3210 2180 35 35  0.0000000 6.2831853}%
%
\special{pn 20}%
\special{sh 0.600}%
\special{ar 2600 2200 35 35  0.0000000 6.2831853}%
%
\special{pn 20}%
\special{sh 0.600}%
\special{ar 2610 1430 35 35  0.0000000 6.2831853}%
%
\special{pn 20}%
\special{sh 0.600}%
\special{ar 2010 1440 35 35  0.0000000 6.2831853}%
%
\special{pn 20}%
\special{sh 0.600}%
\special{ar 2010 2200 35 35  0.0000000 6.2831853}%
%
\special{pn 20}%
\special{sh 0.600}%
\special{ar 1400 1430 35 35  0.0000000 6.2831853}%
%
\special{pn 20}%
\special{sh 0.600}%
\special{ar 1400 2200 35 35  0.0000000 6.2831853}%
\end{picture}%
\end{center}
\caption{$8$-point Blow-up of Hirzebruch surface of degree $2$}
\label{fig:hirzebruch}
\end{figure}
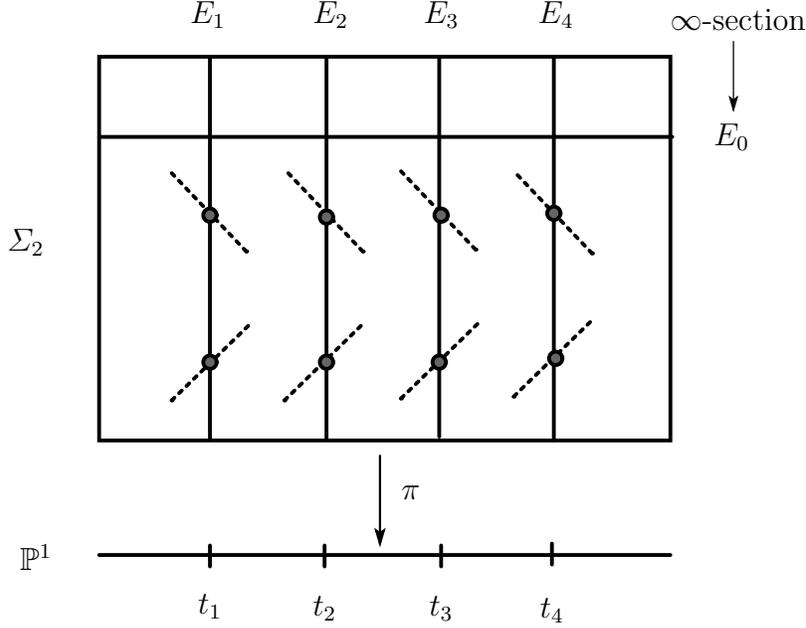
\par 
Note that $\ol{\M}_t(\k)$ is an example of generalized Halphen 
surfaces (see Definition \ref{def:Halphen}), which were 
introduced and classified by Sakai \cite{Sakai}; namely, 
a surface of type $D_4^{(1)}$ in his classification. 
\begin{definition}[Generalized Halphen Surface] 
\label{def:Halphen} 
A smooth, projective, rational surface $S$ is called a 
{\it generalized Halphen surface} if $S$ has an 
effective anti-canonical divisor 
\begin{equation} \label{eqn:Halphen} 
Y \in |-K_{S}| \qquad \mbox{such that} \qquad  
K_{S} \cdot Y_i = 0 \quad (i = 1,\dots,r), 
\end{equation}
where $Y_1,\dots,Y_r$ are the irreducible components of $Y$. 
\end{definition}  
This notion was introduced to construct discrete Painlev\'e 
equations as Cremona transformations of generalized Halphen 
surfaces and to obtain continuous Painlev\'e equations 
as their continuous limits (Cremona approach in 
Remark \ref{rem:approach}). 
\par 
Furthermore, the pair $(\ol{\M}_t(\k), Y_t(\k))$ is an instance 
of Okamoto-Painlev\'e pairs (see Definition \ref{def:OPpair}), 
which were introduced and classified by Saito, Takebe and 
Terajima \cite{STT, STa}; namely, a pair of type 
$\tilde{D}_4$ (or of type $I_0^*$ in Kodaira's notation) in 
their classification. 
\begin{definition}[Okamoto-Painlev\'e Pair] 
\label{def:OPpair}
A pair $(S,Y)$ is called a 
{\it generalized Okamoto-Painlev\'e pair} if 
$S$ is a smooth, projective surface and $Y \in |-K_S|$ 
is an effective anti-canonical divisor satisfying the 
condition (\ref{eqn:Halphen}).  
It is called an {\it Okamoto-Painlev\'e pair} if moreover  
$S - Y_{\mathrm{red}}$ contains an affine plane $\C^2$ as a 
Zariski open subset and $F := S - \C^2$ is a (reduced) divisor 
with normal crossings. 
\end{definition}
This notion was introduced to construct continous Painlev\'e 
equations as Kodaira-Spencer deformations of Okamoto-Painlev\'e 
pairs (Kodaira-Spencer approach in Remark \ref{rem:approach}). 
\par 
Definitions \ref{def:Halphen} and \ref{def:OPpair} were   
invented by speculating on the meanings of the spaces 
constructed by Okamoto \cite{Okamoto1}. 
Here is a comparison of our moduli spaces with his spaces. 
\begin{remark}[Comparison with Okamoto's space] 
\label{rem:okamoto} 
Theorem \ref{thm:realization} implies that our phase 
space $\M_t(\k)$ is isomorphic to the space constructed 
by Okamoto \cite{Okamoto1}. 
He constructed it by hand, chasing trajectories of 
differential equation (\ref{eqn:PVI})\footnote{to be more 
precise, a Hamitoninan system associated to 
equation (\ref{eqn:PVI})}, 
blowing up the points where distinct trajectories meet 
together and removing the vertical leaves.  
Our construction is more theoretical and 
intrinsic\footnote{Painlev\'e property 
follows from our construction, while 
it was presupposed in his construction.}.  
More importantly, our moduli-theoretical construction 
immediately allows us to consider the Riemann-Hilbert 
correspondence from the constructed space  
(to a moduli space of monodromy representations), since 
each point of which represents a parabolic connection.  
This means that we are in a happy situation that 
the construction of the phase space immediately results 
in the construction of a natural conjugacy map. 
\end{remark} 
\par 
Digressively, we take this opportunity to collect the 
major approaches to Painlev\'e equations we have ever 
encountered.  
Gathering those mentioned in the Introduction and 
those remarked after Definitions \ref{def:Halphen} 
and \ref{def:OPpair}, we have (at least) five approaches.  
\begin{remark}[Approaches to Painlev\'e Equations] 
\label{rem:approach} 
$\phan{a}$ 
\begin{center}
\begin{tabular}{llcll}
(1) & Isomonodromic (Fuchs) approach & & 
(2) & Lyapunov approach \\[3mm]  
(3) & Cremona approach               & & 
(4) & Kodaira-Spencer approach \\[3mm]
(5) & Riemann-Hilbert approach       & &  
    & 
\end{tabular}
\end{center} 
As is mentioned in the Introduction, the isomonodromic 
approach and the Riemann-Hilbert approach are close 
relatives. 
In this context, the meaning of our moduli-theoretical 
construction is that we were able to match Okamoto's 
spaces with the isomonodromic picture, which had hitherto    
existed independently, in the framework of 
Riemann-Hilbert approach.   
On the other hand, his spaces have {\sl a priori} had 
their raison d'\^{e}tre in the Cremona and Kodaira-Spencer 
approaches, since these approaches originate from 
searches for their intrinsic meanings.    
\end{remark}  
\section{Riemann-Hilbert Correspondence} \label{sec:RH}  
In the Riemann-Hilbert approach, undoubtedly,  
the Riemann-Hilbert correspondence plays a central part, 
as a (quasi-)conjugacy map between the Painlev\'e flow and 
the isomonodromic flow. 
We start with some basic notions concerning monodromy 
representations.  
\subsection{Monodromy Representations} \label{subsec:monodromy} 
Given $t \in T$, we consider representations of the fundamental 
group $\pi_1(\P^1-D_t,*)$ into $SL_2(\C)$, where the divisor 
$D_t$ is identified with the $4$-point set $\{t_1,t_2,t_3,t_4\}$.  
Recall that two representations $\rho_1$ and $\rho_2$ are said 
to be {\it isomorphic} if there exists a matrix $P \in SL_2(\C)$ 
such that 
\[
\rho_2(\ga) = P \rho_1(\ga) P^{-1} \qquad \mbox{for any} 
\quad \ga \in \pi_1(\P^1-D_t,*). 
\]
For a precise formulation of the Riemann-Hilbert 
correspondence, we need the concept of Jordan equivalence of 
representations, which is closely related to the 
categorical-quotient construction in algebraic geometry.  
We insist that the usual equivalence up to isomorphisms is 
not appropriate, because the set of all representations 
up to isomorphisms is {\it not} an algebraic variety. 
A more substantial reason will gradually be clear in 
the course of discussions: by a categorical-quotient 
formulation, the Riemann-Hilbert correspondence will 
become a resolution of singularities.    
\begin{definition}[Jordan Equivalence] \label{def:jordan} 
A {\it semisimplification} of a representation $\rho$ is 
the associated graded of a composition series of $\rho$. 
Two representations $\rho_1$ and $\rho_2$ are said to be 
{\it Jordan equivalent} if they have isomorphic 
simisimplifications, that is, if either 
\begin{enumerate}
\item they are both irreducible and isomorphic, or 
\item they are both reducible and their 
semisimplifications  $\rho_1'\op \rho_1/\rho_1'$ 
and $\rho_2'\op \rho_2/\rho_2'$ are isomorphic, 
where $\rho_1'$ and $\rho_2'$ are $1$-dimensional 
subrepresentations of $\rho_1$ and $\rho_2$. 
\end{enumerate}
If there is no danger of confusion, a representation and its 
Jordan equivalence class will be denoted by the same symbol. 
For each $t \in T$ let $\R_t$ denote the set of all Jordan 
equivalence classes of $SL_2(\C)$-representarions of 
$\pi_1(\P^1-D_t,*)$. 
We can also speak of the family  
\begin{equation}  \label{eqn:R}
\R = \coprod_{t \in T} R_t. 
\end{equation} 
\end{definition} 
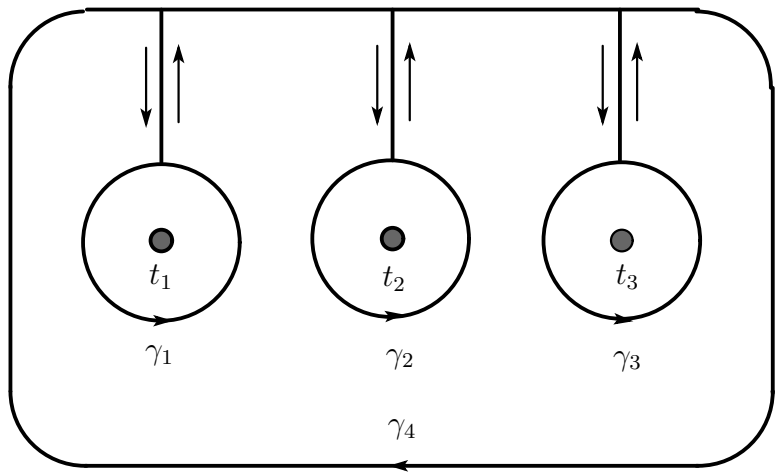
\begin{figure}[t]
\begin{center}
\unitlength 0.1in
\begin{picture}(39.90,23.95)(2.10,-25.95)
%
\special{pn 20}%
\special{ar 600 600 390 390  3.1165979 4.6873942}%
%
\special{pn 20}%
\special{pa 210 600}%
\special{pa 210 2200}%
\special{fp}%
%
\special{pn 20}%
\special{ar 600 2200 390 390  1.5458015 3.1415927}%
%
\special{pn 20}%
\special{ar 2200 1400 410 410  0.0000000 6.2831853}%
%
\special{pn 20}%
\special{ar 3410 1410 410 410  0.0000000 6.2831853}%
%
\special{pn 20}%
\special{pa 610 200}%
\special{pa 3810 200}%
\special{fp}%
%
\special{pn 20}%
\special{pa 600 2590}%
\special{pa 3810 2590}%
\special{fp}%
%
\special{pn 20}%
\special{ar 3810 2200 390 390  6.2831853 6.2831853}%
\special{ar 3810 2200 390 390  0.0000000 1.5964317}%
%
\special{pn 20}%
\special{pa 1000 200}%
\special{pa 1000 1000}%
\special{fp}%
%
\special{pn 20}%
\special{pa 2210 200}%
\special{pa 2210 990}%
\special{fp}%
%
\special{pn 20}%
\special{pa 3400 200}%
\special{pa 3400 990}%
\special{fp}%
%
\special{pn 20}%
\special{pa 4200 610}%
\special{pa 4200 2210}%
\special{fp}%
%
\special{pn 20}%
\special{ar 3800 600 391 391  4.7123890 6.2831853}%
\special{ar 3800 600 391 391  0.0000000 0.0256354}%
%
\special{pn 20}%
\special{ar 1000 1420 410 410  0.0000000 6.2831853}%
%
\special{pn 13}%
\special{sh 0.600}%
\special{ar 3410 1410 58 58  0.0000000 6.2831853}%
%
\special{pn 20}%
\special{sh 0.600}%
\special{ar 2210 1400 58 58  0.0000000 6.2831853}%
%
\special{pn 20}%
\special{sh 0.600}%
\special{ar 1000 1410 58 58  0.0000000 6.2831853}%
%
\special{pn 13}%
\special{pa 920 400}%
\special{pa 920 800}%
\special{fp}%
\special{sh 1}%
\special{pa 920 800}%
\special{pa 940 733}%
\special{pa 920 747}%
\special{pa 900 733}%
\special{pa 920 800}%
\special{fp}%
%
\special{pn 13}%
\special{pa 2130 390}%
\special{pa 2130 790}%
\special{fp}%
\special{sh 1}%
\special{pa 2130 790}%
\special{pa 2150 723}%
\special{pa 2130 737}%
\special{pa 2110 723}%
\special{pa 2130 790}%
\special{fp}%
%
\special{pn 13}%
\special{pa 3310 390}%
\special{pa 3310 790}%
\special{fp}%
\special{sh 1}%
\special{pa 3310 790}%
\special{pa 3330 723}%
\special{pa 3310 737}%
\special{pa 3290 723}%
\special{pa 3310 790}%
\special{fp}%
%
\special{pn 13}%
\special{pa 1090 790}%
\special{pa 1090 400}%
\special{fp}%
\special{sh 1}%
\special{pa 1090 400}%
\special{pa 1070 467}%
\special{pa 1090 453}%
\special{pa 1110 467}%
\special{pa 1090 400}%
\special{fp}%
%
\special{pn 13}%
\special{pa 2300 780}%
\special{pa 2300 390}%
\special{fp}%
\special{sh 1}%
\special{pa 2300 390}%
\special{pa 2280 457}%
\special{pa 2300 443}%
\special{pa 2320 457}%
\special{pa 2300 390}%
\special{fp}%
%
\special{pn 13}%
\special{pa 3490 780}%
\special{pa 3490 390}%
\special{fp}%
\special{sh 1}%
\special{pa 3490 390}%
\special{pa 3470 457}%
\special{pa 3490 443}%
\special{pa 3510 457}%
\special{pa 3490 390}%
\special{fp}%
%
\special{pn 13}%
\special{pa 950 1830}%
\special{pa 1040 1830}%
\special{fp}%
\special{sh 1}%
\special{pa 1040 1830}%
\special{pa 973 1810}%
\special{pa 987 1830}%
\special{pa 973 1850}%
\special{pa 1040 1830}%
\special{fp}%
%
\special{pn 13}%
\special{pa 2150 1810}%
\special{pa 2250 1800}%
\special{fp}%
\special{sh 1}%
\special{pa 2250 1800}%
\special{pa 2182 1787}%
\special{pa 2197 1805}%
\special{pa 2186 1827}%
\special{pa 2250 1800}%
\special{fp}%
%
\special{pn 13}%
\special{pa 3370 1820}%
\special{pa 3450 1830}%
\special{fp}%
\special{sh 1}%
\special{pa 3450 1830}%
\special{pa 3386 1802}%
\special{pa 3397 1823}%
\special{pa 3381 1842}%
\special{pa 3450 1830}%
\special{fp}%
%
\special{pn 13}%
\special{pa 2300 2590}%
\special{pa 2210 2590}%
\special{fp}%
\special{sh 1}%
\special{pa 2210 2590}%
\special{pa 2277 2610}%
\special{pa 2263 2590}%
\special{pa 2277 2570}%
\special{pa 2210 2590}%
\special{fp}%
\put(9.9000,-20.1000){\makebox(0,0){$\gamma_1$}}%
\put(22.5000,-20.3000){\makebox(0,0){$\gamma_2$}}%
\put(34.4000,-20.4000){\makebox(0,0){$\gamma_3$}}%
\put(22.6000,-24.0000){\makebox(0,0){$\gamma_4$}}%
\put(10.0000,-16.0000){\makebox(0,0){$t_1$}}%
\put(22.2000,-16.1000){\makebox(0,0){$t_2$}}%
\put(34.4000,-16.0000){\makebox(0,0){$t_3$}}%
\end{picture}%
\end{center}
\caption{The loops $\ga_i$; the fourth point $t_4$ is 
outside $\ga_4$, invisible.}  
\label{fig:loops} 
\end{figure}
\begin{definition}[Local Monodromy Data] \label{def:lmd} 
We put $A := \C^4$ and consider the map 
\begin{equation} \label{eqn:a} 
\pi_t : \R_t \rightarrow A, \quad \rho \mt a 
= (a_1,a_2,a_3,a_4), \qquad a_i = \Tr\, \rho(\ga_i).  
\end{equation} 
where $\ga_i \in \pi_1(\P^1-D_t,*)$ is a loop surrounding 
the point $t_i$ anti-clockwise, leaving the remaining three 
points outside, as in Figure \ref{fig:loops}. 
Note that $a_i$ is well-defined, that is, it depends only on 
the Jordan equivalence class of $\rho$ and does not depend 
on the choice of loop $\ga_i$.  
We call $a$ the {\it local monodromy data} of $\rho$. 
For each $a \in A$ let $\R_t(a)$ denote the fiber of the 
map (\ref{eqn:a}) over $a$. 
As the relative setting of (\ref{eqn:a}) over $T$, 
we have the family  
\[
\pi : \R \rightarrow T \times A,  
\]
where $\R$ is defined by (\ref{eqn:R}). 
For a fixed $a \in A$ we also have the family 
$\pi_a : \R(a) \to T$ as in (\ref{eqn:R(a)}). 
\end{definition} 
\subsection{Riemann-Hilbert Correspondence} 
\label{subsec:RH}
To formulate the Riemann-Hilbert correspondence, 
we first set it up in the parameter level. 
\begin{definition}[Riemann-Hilbert Correspondence in   
Parameter Level] \label{def:rh0} 
We consider \linebreak the correspondence of local exponents 
to local monodromy data   
\begin{equation} \label{eqn:rh3}
\rh : \K \to A, \qquad 
\k = (\k_0,\k_1,\k_2,\k_3,\k_4) \mt 
a = (a_1,a_2,a_3,a_4).  
\end{equation}
From Table \ref{tab:riemann} the monodromy matrix 
$\rho(\ga_i)$ along the loop $\ga_i$ has eigenvalues 
$\exp(\pm2\pi\sqrt{-1}\l_i)$ and hence has trace 
$2 \cos 2\pi\l_i$.  
Then (\ref{eqn:lambda}) and (\ref{eqn:a}) imply that  
in terms of exponents $\k \in \K$, the local monodromy 
data of $\rho$ is expressed as 
\begin{equation}  \label{eqn:rh2} 
a_i = \left\{
\begin{array}{ll}
\phan{-}2 \cos \pi \k_i \qquad & (i = 1,2,3), \\[2mm]
-2 \cos \pi \k_4 \qquad & (i = 4). 
\end{array} 
\right. 
\end{equation} 
The map (\ref{eqn:rh3}) with (\ref{eqn:rh2}) is called 
the {\it Riemann-Hilbert correspondence 
in the parameter level}. 
\end{definition}
\begin{definition}[Riemann-Hilbert Correspondence] 
\label{def:RH}
Given $t \in T$, any stable parabolic connection 
$Q = (E,\nabla,\psi,l) \in \M_t$, upon restricted to 
$\P^1-D_t$, induces a flat connection 
\[
\nabla|_{\P^1-D_t} : E|_{\P^1-D_t} \rightarrow 
E|_{\P^1-D_t} \ot \Om^1_{\P^1-D_t}. 
\]
Let $\rho$ be the Jordan equivalence class of its monodromy 
representation. 
Then the {\it Riemann-Hilbert correspondence} at time $t$ 
is defined by the holomorphic map 
\[
\RH_t : \M_t \to \R_t, \qquad Q \mt \rho.  
\]
By Definition \ref{def:rh0} there exists a commutative 
diagram of holomorophic maps 
\begin{equation} \label{cd:RH1} 
\begin{CD}
\M_t @> \RH_t >> \R_t \\
 @V \pi_t VV @VV \pi_t V \\
\K @>>  \rh > A,  
\end{CD}
\end{equation}
where $\pi_t : \M_t \to \K$ is the map sending each 
parabolic connection to its local exponents and the map 
$\pi_t : \R_t \to A$ is defined by (\ref{eqn:a}). 
As the relative setting of (\ref{cd:RH1}) over $T$, we have 
\begin{equation} \label{cd:RH1a} 
\begin{CD}
\M @> \RH >> \R \\
 @V \pi VV @VV \pi V \\
T \times \K @>>  \id \times \rh > T \times A.   
\end{CD} 
\end{equation}
\end{definition} 
\par
Since $\rh : \K \to A$ is an infinite-to-one map, so is 
the base map of (\ref{cd:RH1a}). 
This fact makes the analysis of (\ref{cd:RH1a}) somewhat 
difficult.  
To avoid this we consider the fiber product $\bR$ defined 
by 
\begin{equation} \label{cd:FP} 
\begin{CD}
\bR @>>> \R \\
 @V \pi VV @VV \pi V \\
T \times \K @>> \id \times \rh > T \times A.  
\end{CD}
\end{equation}
\par
We now set up three versions of Riemann-Hilbert 
correspondence that will be used later. 
\begin{definition}[Three Versions of Riemann-Hilbert 
Correspondence] \label{def:RH2} $\phan{a}$
\begin{enumerate}
\item From (\ref{cd:RH1a}) and (\ref{cd:FP}) we have the 
commutative diagram of holomorphic maps  
\begin{equation} \label{cd:RH3}
\begin{CD} 
\M @> \RH >> \bR \\
@V \pi VV   @VV \pi V \\
T \times \K @= T \times \K, 
\end{CD} 
\end{equation} 
which is called the {\it full-Riemann-Hilbert 
correspondence}.  
\item Fix an exponent $\k \in \K$ and put $a = \rh(\k) \in A$. 
Then (\ref{cd:RH1a}) restricts to the diagram  
\begin{equation} \label{cd:RH2} 
\begin{CD}
\M(\k) @> \RH_{\k} >> \R(a) \\
@V \pi_{\k} VV   @VV \pi_a V \\
T @= T,  
\end{CD}
\end{equation}
which is referred to as the {\it $\k$-Riemann-Hilbert 
correspondence}.  
\item Moreover, upon fixing a time $t \in T$, 
diagram (\ref{cd:RH2}) further restricts to the map 
\begin{equation} \label{eqn:RH4}
\RH_{t,\k} : \M_t(\k) \rightarrow \R_t(a),  
\end{equation} 
which is referred to as the {\it $(t,\k)$-Riemann-Hilbert 
correspondence}. 
\end{enumerate} 
\end{definition} 
\par 
Among the three versions above, the importance of 
(\ref{cd:RH2}) and (\ref{eqn:RH4}) is obvious:  
(\ref{cd:RH2}) will serve as a (quasi-)conjugacy map 
of Painlev\'e flow to isomonodromic flow, while 
(\ref{eqn:RH4}) will give a correspondence between the 
spaces of initial-conditions for these two dynamics. 
On the other hand, although it is not yet clear, 
(\ref{cd:RH3}) will play an important part in constructing 
Painlev\'e flows based on ``codimension-two argument'' 
(see Lemma \ref{lem:codim} and Remark \ref{rem:CTA}). 
\par 
The Riemann-Hilbert problem usually asks the 
surjectivity of Riemann-Hilbert correspondence. 
But the injectivity and properness are also important 
issues in our situation. 
\begin{problem}[Riemann-Hilbert Problem] \label{prob:RH} 
We formulate the problems for $\RH_{\k}$ in (\ref{cd:RH2}). 
\begin{enumerate}
\item Is $\RH_{\k}$ surjective? 
This question is fundamental for the whole development 
of the story. 
\item To what extent $\RH_{\k}$ is injective?  
This question is important for the setup of $\RH_{\k}$  
as a (quasi-)conjugacy map between the 
Painlev\'e flow and the isomonodromic flow. 
\item Is $\RH_{\k}$ a proper map?  
This question is important becuase the properness of 
$\RH_{\k}$ leads to the geometric Painlev\'e property 
of the Painlev\'e flow (see Lemma \ref{lem:proper}).    
\end{enumerate} 
\end{problem} 
In what follows, Riemann-Hilbert problem will often be  
abbreviated to RHP. 
Those for $\RH$, $\RH_{\k}$ and $\RH_{t,\k}$  
will be referred to as full-RHP, $\k$-RHP and 
$(t,\k)$-RHP, respectively. 
\subsection{Affine Weyl Group of 
Type $\mbox{\bm $D_4^{(1)}$}$} 
\label{subsec:weyl}
Before stating our solution to the Riemann-Hilbert problem, 
we introduce an affine Weyl group of type $D_4^{(1)}$ 
acting on the parameter space $\K$ (see 
Definition \ref{def:weyl}) and characterize the 
singularities of $\R_t(a)$ in terms of 
the affine Weyl group structure (see Lemma \ref{lem:sing}).  
In connection with the singularity structure, we introduce 
the concept of Riccati loci (see Defintion \ref{def:riccati}).  
\begin{definition}[Affine Weyl Group] \label{def:weyl}  
The parameter space $\K$ in (\ref{eqn:K}) is an affine space 
modeled on the four-dimensional linear space
\[
K = \{k = (k_0,k_1,k_2,k_3,k_4) \in \C^5 \,:\, 
2 k_0 + k_1 + k_2 + k_3 +k_4 = 0\}, 
\]
endowed with the inner product 
$\la k, \, k' \ra = k_1 k_1'+k_2 k_2'+k_3 k_3'+k_4 k_4'$. 
Let $\si_i$ be the orthogonal affine reflection on $\K$  
having $\{\k \in\K \,:\, \k_i = 0\}$ as its reflecting hyperplane. 
We observe that the group generated by $\si_0$, $\si_1$, 
$\si_2$, $\si_3$, $\si_4$ is an affine Weyl group 
of type $D_4^{(1)}$ (see Figure \ref{fig:dynkin}), 
\[
\W = \langle \si_0, \si_1, \si_2, \si_3, 
\si_4 \rangle.
\]
If $C = (c_{ij})$ is the Cartan matrix of type $D_4^{(1)}$, 
the $i$-th basic reflection $\si_i$ is expressed as 
\begin{equation} \label{eqn:sigma} 
\si_i(\k_j) = \k_j - \k_i c_{ij}. 
\end{equation} 
Let $\Wall \subset \K$ denote the union of the reflecting 
hyperplanes of all reflections in $\W$.  
\end{definition} 
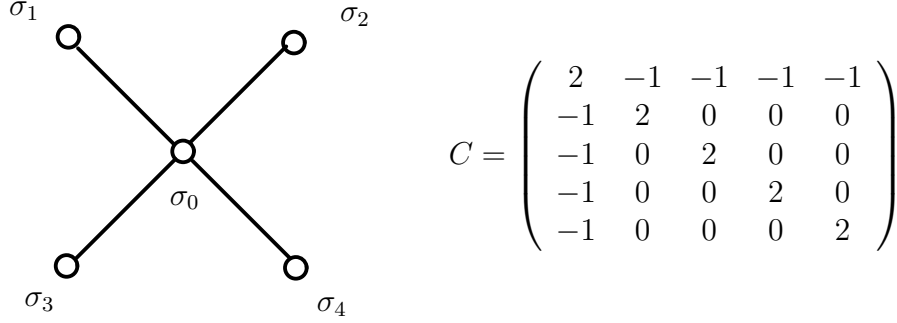
\begin{figure}[t]
\begin{center}
\unitlength 0.1in
\begin{picture}(54.20,15.15)(0.00,-17.25)
%
\special{pn 20}%
\special{ar 1200 995 58 58  6.1412883 6.2831853}%
\special{ar 1200 995 58 58  0.0000000 6.1180366}%
%
\special{pn 20}%
\special{ar 600 395 58 58  6.1412883 6.2831853}%
\special{ar 600 395 58 58  0.0000000 6.1180366}%
%
\special{pn 20}%
\special{ar 1780 425 58 58  6.1412883 6.2831853}%
\special{ar 1780 425 58 58  0.0000000 6.1180366}%
%
\special{pn 20}%
\special{ar 590 1595 58 58  6.1412883 6.2831853}%
\special{ar 590 1595 58 58  0.0000000 6.1180366}%
%
\special{pn 20}%
\special{ar 1790 1605 58 58  6.1412883 6.2831853}%
\special{ar 1790 1605 58 58  0.0000000 6.1180366}%
%
\special{pn 20}%
\special{pa 650 455}%
\special{pa 1150 955}%
\special{fp}%
%
\special{pn 20}%
\special{pa 1250 1045}%
\special{pa 1760 1555}%
\special{fp}%
%
\special{pn 20}%
\special{pa 1740 445}%
\special{pa 1240 945}%
\special{fp}%
%
\special{pn 20}%
\special{pa 1160 1035}%
\special{pa 640 1555}%
\special{fp}%
\put(12.1000,-12.5000){\makebox(0,0){$\sigma_0$}}%
\put(21.8000,-2.4000){\makebox(0,0)[rt]{$\sigma_2$}}%
\put(19.8000,-18.1000){\makebox(0,0){$\sigma_4$}}%
\put(4.5000,-17.8000){\makebox(0,0){$\sigma_3$}}%
\put(2.9000,-2.1000){\makebox(0,0)[lt]{$\sigma_1$}}%
\put(25.9000,-15.2000){\makebox(0,0)[lb]{$C = \left(\begin{array}{ccccc} \-2 & -1 & -1 & -1 & -1 \\ -1 & \-2 & \-0 & \-0 & \-0 \\ -1 & \-0 & \-2 & \-0 & \-0 \\ -1 & \-0 & \-0 & \-2 & \-0 \\ -1 & \-0 & \-0 & \-0 & \-2 \end{array} \right)$}}%
\put(54.2000,-6.1000){\makebox(0,0)[lb]{$\phan{a}$}}%
\end{picture}%
\end{center}  
\par\vspace{0.5cm}\noindent 
\caption{Dynkin diagram and Cartan matrix of type $D_4^{(1)}$}
\label{fig:dynkin} 
\end{figure}   
We remark that a more intrinsic presentation of 
Definition \ref{def:weyl} is possible along the line of Arinkin 
and Lysenko \cite{AL2}, as the Weyl group on the Picard 
lattice of the moduli space $\M_t(\k)$.   
\par
Let $\R_t^{\s}(a)$ be the singular locus and 
$\R_t^{\ci}(a) = \R_t(a) - \R_t^{\s}(a)$ be the smooth locus 
of $\R_t(a)$, respectively.    
The affine Weyl group structure allows us to describe 
the singularities of $\R_t(a)$.   
\begin{lemma}[Singularity] \label{lem:sing} 
Let $\k \in \K$ and put $a = \rh(\k) \in A$.  
\begin{enumerate}
\item The surface $\R_t(a)$ is smooth, that is, 
$\R_t^{\s}(a) = \emptyset$ if and only if 
$\k \not\in\Wall$. 
\item If $\k \in \Wall$, the singular locus $\R_t^{\s}(a)$ 
consists of at most four rational double points. 
\end{enumerate} 
\end{lemma}
The possible types of singularities on $\R_t(a)$ will be 
classified completely in Theorem \ref{thm:classification}.  
In connection with the singular loci of the surfaces  
$\R_t(a)$, we make the following definition. 
\begin{definition}[Riccati/Non-Riccati Loci] 
\label{def:riccati} 
The {\it Riccati loci} are defined by
\[
\bR^{\r} = \coprod_{(t,\k) \in T \times \K} \R^{\s}_t(\rh(\k)), 
\qquad \M^{\r} = \RH^{-1}(\bR^{\r}). 
\]
By Lemma \ref{lem:sing} the disjoint union may be 
taken only over $T \times \Wall$. 
The {\it non-Riccati loci} 
\[
\bR^{\ci} = \bR - \bR^{\r}, \qquad \M^{\ci} = \M - \M^{\r}.  
\]  
are the complements to the Riccati loci. 
These loci are restricted to subspaces $\R(a)$, $\R_t(a)$, 
$\M(\k)$, $\M_t(\k)$ with $a = \rh(\k)$ in an obvious manner: 
The Riccati loci for them are defined by    
\[
\begin{array}{cccccc}
\R^{\r}(a) &=& \ds \coprod_{t \in T} \R_t^{\s}(a), 
\qquad & \M^{\r}(\k) &=& \RH_{\k}^{-1}(\R^{\r}(a)), \\[6mm] 
\R_t^{\r}(a) &=& \R_t^{\s}(a), 
\qquad & \M_t^{\r}(\k) &=& \RH_{t,\k}^{-1}(\R_t^{\r}(a)).     
\end{array} 
\] 
The corresponding non-Riccati loci are the complements to 
them: 
\[ 
\begin{array}{cccccc}
\R^{\ci}(a) &=& \R(a) - \R^{\r}(a), 
\qquad & \M^{\ci}(\k) &=& \M(\k) - \M^{\r}(\k), \\[3mm] 
\R_t^{\ci}(a) &=& \R_t(a) - \R_t^{\r}(a),  
\qquad & \M_t^{\ci}(\k) &=& \M_t(\k) - \M_t^{\r}(\k).    
\end{array} 
\] 
\end{definition} 
\par  
It will turn out that Riccati loci are closely related to 
the so-called Riccati solutions of the Painlev\'e equation. 
This fact motivates the terminology {\it Riccati locus} 
(see \S\ref{subsec:RF}).  
\subsection{Solution to Riemann-Hilbert Problem} 
\label{subsec:RHP}
We are now in a position to state our solution to the 
Riemann-Hilbert problem \cite{IIS2}. 
\begin{theorem}[Solution to Full-RHP] \label{thm:RH} 
$\phan{a}$ 
\begin{enumerate} 
\item $\RH : \M \to \bR$ is a surjective proper holomorphic map,  
and 
\item $\RH : \M^{\ci} \to \bR^{\ci}$ is a biholomorophism. 
\end{enumerate}       
\end{theorem}
\par
Restricting this theorem to each $\k \in \K$, 
we have the following corollary.      
\begin{corollary}[Solution to $\mbox{\bm$\k$}$-RHP] 
\label{cor:RH2} 
Let $\k \in \K$ and put $a = \rh(\k) \in A$. 
\begin{enumerate}
\item $\RH_{\k} : \M(\k) \to \R(a)$ is a surjective 
proper holomorphic map, and
\item $\RH_{\k} : \M^{\ci}(\k) \to \R^{\ci}(a)$ is a 
biholomorphism.    
\end{enumerate}
\end{corollary}
\par 
Furthermore, at each $(t,\k)$-level we have the 
following theorem. 
\begin{theorem}[Solution to $\mbox{\bm$(t,\k)$}$-RHP]
\label{thm:RH3}
Let $(t,\k) \in T \times \K$ and put $a = \rh(\k) \in A$. 
\begin{enumerate} 
\item If $\k \not\in \Wall$, then   
$\RH_{t,\k} : \M_t(\k) \to \R_t(a)$ is a biholomorphic map, and    
\item if $\k \in \Wall$, then 
$\RH_{t,\k} : \M_t(\k) \to \R_t(a)$ is a minimal 
resolution of singularities having the Riccati locus 
$\M_t^{\r}(\k)$ as its exceptional divisor.  
\end{enumerate}
\end{theorem}
\par 
These theorems can be generalized to stable parabolic 
connections of higher rank, with more regular singular points,  
and even on a curve of higher genus. 
We present an essence of the proof, focusing on the surjectivity 
of $\RH$, which remains valid for such generalizations.        
\begin{remark}[How to Prove] \label{rem:howtosolve} 
Given a Jordan equivalence class of representations,  
\begin{enumerate}
\item choose a ``good'' representative from the given 
equivalence class and form the flat connection 
associated to it. 
Since we are working with Jordan equivalence, we can take 
a semisimple representation $\rho_0$ as the good 
representative. 
\item Extend the flat connection to a logarithmic connection 
by Deligne's canonical extension \cite{Deligne} and provide 
it with a parabolic structure. 
If the initial representation $\rho_0$ is irreducible, 
the resulting parabolic connection $Q_0$ is stable and 
so we are done. 
If $\rho_0$ is reducible, we cannot stop here because 
$Q_0$ may be unstable and we should proceed to step (3).  
\item If $\rho_0$ is reducible, take steps (1) 
and (2) relatively, so that we obtain a family of 
parabolic connections $\mathcal{Q} = \{Q_c\}_{c \in C}$ 
parametrized by some curve $C$, with $Q_{c_0} = Q_0$ 
at the reference point $c_0 \in C$, such that the 
monodromy of $Q_c$ is irreducible for every 
$c \in C - \{c_0\}$. 
Then use Langton's technique in Theorem 
\ref{thm:langton} to recast $Q_0$ to a {\it stable} 
parabolic connection. 
\item The family $\mathcal{Q}$ in step (3) is 
constructed as follows. 
Notice that reducible representations occur only on 
a Zariski-closed proper subset $B \subset A$.  
Let $c_0 \in B$ be the local monodromy data of 
$\rho_0$, take a curve $C \subset A$ that  
meets $B$ only at $c_0$, and prolong the representation 
$\rho_0$ along the curve $C$. 
Taking steps (1) and (2) relatively, we obtain the  
desired family $\mathcal{Q}$.            
\end{enumerate}  
\end{remark} 
\par 
Here is the version of Langton's technique \cite{Langton} 
that is needed in the current situation.  
\begin{theorem}[Langton's Technique] \label{thm:langton} 
Let $\mathcal{Q} = \{Q_c\}_{c \in C}$ be a family of parabolic 
connections parametrized by a curve $C$. 
By some applications of elementary transformations, 
$\mathcal{Q}$ can be transformed to a family of 
\ul{stable} parabolic connections, if the monodromy of 
$Q_c$ is \ul{irreducible} for every $c \in C-\{c_0\}$. 
This means that the possible singularity of $\mathcal{Q}$ at 
$c_0$ can be removed by elementary tramsformations, 
provided that all the nearby connections are irreducible.    
\end{theorem} 
\par 
Langton's theorem reminds us of the removable singularity 
theorem of Riemann in complex variable and that of Uhlenbeck
in gauge theory.     
Riemann's classical theorem asserts that an isolated 
singularity of a holomorphic function can be removed, 
if the function is {\it bounded} around the singular point.   
Uhlenbeck's theorem \cite{Uhlenbeck} states that an isolated 
singularity of a Young-Mills connection can be removed by 
applying a gauge transformation, if the curvature of the 
connection is $L^2$-{\it bounded} around the singular point. 
Langton's theorem can be regarded as an 
algebraic-geometry version of such removable singularity 
principles, where the boundedness condition is 
replaced by the {\it irreducibility} of representations. 
\par   
\begin{remark}[Family of \mbox{\bm $(-2)$}-Curves] 
\label{rem:minustwo}  
By Theorem \ref{thm:RH3}, for any 
$(t,\k) \in T \times \Wall$, the $(t,\k)$-Riemann-Hilbert 
correspondence $\RH_{t,\k} : \M_t(\k) \to \R_t(a)$ gives a 
minimal resolution of singularities whose exceptional 
divisor is just the Riccati locus $\M_t^{\r}(\k)$. 
Each irreducible component of $\M_t^{\r}(\k)$    
is a $(-2)$-{\it curve}, that is, a smooth curve 
$C \subset \M_t(\k)$ such that 
\[
C \simeq \P^1, \qquad  C\cdot C = -2. 
\]
Conversely any $(-2)$-curve in $\M_t(\k)$ arises in this 
way, since it must be sent to a singular point by 
$\RH_{t,\k}$.  
Considering this situation relatively for the family 
$\pi_{\k} : \M(\k) \to T$, we see that each irreducible 
component of the Ricatti locus $\M^{\r}(\k) \subset \M(\k)$ 
is a family of $(-2)$-curves over $T$, 
\begin{equation} \label{eqn:minustwo}
\pi_{\k} : \CC \rightarrow T, \qquad 
C_t \subset \M_t(\k) \,\,:\,\, \mbox{$(-2)$-curve}. 
\end{equation}
\end{remark} 
\par
To apply Hartog's extension theorem later, we state the 
following simple lemma. 
\begin{lemma}[Codimension Two] \label{lem:codim} 
The Riccati locus $\M^{\r}$ is of codimension two in 
$\M$. 
\end{lemma}  
This is intuitively clear: By Lemma \ref{lem:sing} the 
Riccati locus $\M^{\r}$ can lie only over the 
codimension-one subset $T \times \Wall \subset T \times \K$ 
with respect to the fibration (\ref{eqn:family1}). 
On the other hand Remark \ref{rem:minustwo} implies that 
for each $(t,\k) \in T \times \Wall$, the Riccati locus 
$\M_t^{\r}(\k)$ is of codimension one in $\M_t(\k)$. 
In total $\M^{\r}$ is of codimension two in $\M$. 
Lemma \ref{lem:codim} will be used in Remark \ref{rem:CTA}. 
\section{Isomonodromic Flow and Painlev\'e Flow} 
\label{sec:painleve}
\begin{figure}[t]
\begin{center}
\unitlength 0.1in
\begin{picture}(45.40,29.60)(1.70,-33.20)
%
\special{pn 20}%
\special{pa 2620 600}%
\special{pa 2620 2600}%
\special{fp}%
%
\special{pn 20}%
\special{pa 1450 600}%
\special{pa 4440 600}%
\special{pa 4440 2610}%
\special{pa 1450 2610}%
\special{pa 1450 600}%
\special{fp}%
%
\special{pn 20}%
\special{pa 1450 3196}%
\special{pa 4440 3196}%
\special{fp}%
%
\special{pn 20}%
\special{pa 2630 3146}%
\special{pa 2630 3256}%
\special{fp}%
%
\special{pn 20}%
\special{pa 2250 2200}%
\special{pa 3490 2210}%
\special{fp}%
\special{pa 3490 2210}%
\special{pa 3490 2210}%
\special{fp}%
%
\special{pn 20}%
\special{pa 3430 1000}%
\special{pa 2720 1010}%
\special{fp}%
\special{sh 1}%
\special{pa 2720 1010}%
\special{pa 2787 1029}%
\special{pa 2773 1009}%
\special{pa 2786 989}%
\special{pa 2720 1010}%
\special{fp}%
%
\special{pn 20}%
\special{ar 3430 1610 600 600  4.7287810 6.2831853}%
\special{ar 3430 1610 600 600  0.0000000 1.5199926}%
%
\special{pn 13}%
\special{pa 2240 2090}%
\special{pa 3230 2100}%
\special{fp}%
\special{sh 1}%
\special{pa 3230 2100}%
\special{pa 3164 2079}%
\special{pa 3177 2099}%
\special{pa 3163 2119}%
\special{pa 3230 2100}%
\special{fp}%
%
\special{pn 13}%
\special{pa 2240 1970}%
\special{pa 3230 1980}%
\special{fp}%
\special{sh 1}%
\special{pa 3230 1980}%
\special{pa 3164 1959}%
\special{pa 3177 1979}%
\special{pa 3163 1999}%
\special{pa 3230 1980}%
\special{fp}%
%
\special{pn 20}%
\special{sh 0.600}%
\special{ar 2630 2200 10 10  0.0000000 6.2831853}%
%
\special{pn 20}%
\special{sh 0.600}%
\special{ar 2620 2200 48 48  0.0000000 6.2831853}%
%
\special{pn 20}%
\special{sh 0.600}%
\special{ar 2620 1000 48 48  0.0000000 6.2831853}%
%
\special{pn 13}%
\special{pa 2630 3136}%
\special{pa 4000 3146}%
\special{fp}%
%
\special{pn 13}%
\special{pa 4010 3076}%
\special{pa 2650 3076}%
\special{fp}%
\special{sh 1}%
\special{pa 2650 3076}%
\special{pa 2717 3096}%
\special{pa 2703 3076}%
\special{pa 2717 3056}%
\special{pa 2650 3076}%
\special{fp}%
%
\special{pn 13}%
\special{ar 4000 3106 36 36  4.1243864 6.2831853}%
\special{ar 4000 3106 36 36  0.0000000 2.0344439}%
%
\special{pn 13}%
\special{pa 2630 2700}%
\special{pa 2630 3040}%
\special{fp}%
\special{sh 1}%
\special{pa 2630 3040}%
\special{pa 2650 2973}%
\special{pa 2630 2987}%
\special{pa 2610 2973}%
\special{pa 2630 3040}%
\special{fp}%
\put(47.0000,-16.9000){\makebox(0,0)[lb]{$\R(a)$}}%
\put(22.9000,-24.3000){\makebox(0,0)[lb]{$\rho$}}%
\put(22.6000,-10.8000){\makebox(0,0)[lb]{$\rho'$}}%
\put(47.1000,-32.6000){\makebox(0,0)[lb]{$T$}}%
\put(25.9000,-34.9000){\makebox(0,0)[lb]{$t$}}%
\put(22.9000,-29.4000){\makebox(0,0)[lb]{$\pi_a$}}%
\put(30.9000,-30.5000){\makebox(0,0)[lb]{$\beta$}}%
%
\special{pn 13}%
\special{pa 2550 2140}%
\special{pa 2525 2112}%
\special{pa 2501 2085}%
\special{pa 2477 2057}%
\special{pa 2453 2029}%
\special{pa 2430 2002}%
\special{pa 2407 1974}%
\special{pa 2386 1947}%
\special{pa 2365 1919}%
\special{pa 2345 1891}%
\special{pa 2327 1864}%
\special{pa 2310 1836}%
\special{pa 2295 1808}%
\special{pa 2281 1781}%
\special{pa 2269 1753}%
\special{pa 2259 1726}%
\special{pa 2250 1698}%
\special{pa 2245 1670}%
\special{pa 2241 1643}%
\special{pa 2240 1615}%
\special{pa 2241 1588}%
\special{pa 2245 1560}%
\special{pa 2252 1532}%
\special{pa 2260 1505}%
\special{pa 2271 1477}%
\special{pa 2284 1450}%
\special{pa 2298 1422}%
\special{pa 2314 1395}%
\special{pa 2331 1367}%
\special{pa 2350 1340}%
\special{pa 2370 1312}%
\special{pa 2391 1285}%
\special{pa 2413 1257}%
\special{pa 2436 1230}%
\special{pa 2459 1202}%
\special{pa 2482 1175}%
\special{pa 2506 1147}%
\special{pa 2530 1120}%
\special{sp 0.070}%
%
\special{pn 13}%
\special{pa 2530 1110}%
\special{pa 2560 1070}%
\special{dt 0.045}%
\special{sh 1}%
\special{pa 2560 1070}%
\special{pa 2504 1111}%
\special{pa 2528 1113}%
\special{pa 2536 1135}%
\special{pa 2560 1070}%
\special{fp}%
\put(24.4000,-5.3000){\makebox(0,0)[lb]{$\R_t(a)$}}%
\put(1.7000,-15.5000){\makebox(0,0)[lb]{nonlinear}}%
\put(1.7000,-18.2000){\makebox(0,0)[lb]{monodromy}}%
\put(38.0000,-29.3000){\makebox(0,0)[lb]{braid}}%
%
\special{pn 20}%
\special{sh 0.600}%
\special{ar 2620 3200 48 48  0.0000000 6.2831853}%
%
\special{pn 8}%
\special{pa 3700 2880}%
\special{pa 3350 2910}%
\special{fp}%
\special{sh 1}%
\special{pa 3350 2910}%
\special{pa 3418 2924}%
\special{pa 3403 2905}%
\special{pa 3415 2884}%
\special{pa 3350 2910}%
\special{fp}%
%
\special{pn 8}%
\special{pa 1230 1610}%
\special{pa 2070 1610}%
\special{fp}%
\special{sh 1}%
\special{pa 2070 1610}%
\special{pa 2003 1590}%
\special{pa 2017 1610}%
\special{pa 2003 1630}%
\special{pa 2070 1610}%
\special{fp}%
\end{picture}%
$\phan{aaaaaaa}$
\end{center}
\caption{Isomonodromic flow} \label{fig:IMF}
\end{figure}
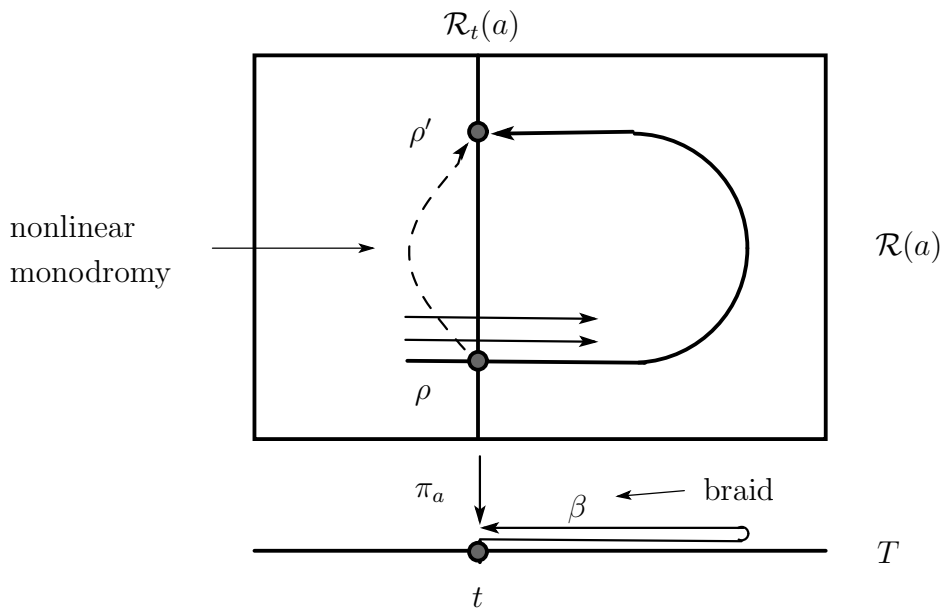
From our dynamical point of view, we should consiously 
distingush the Painlev\'e flow on the moduli space of stable 
parabolic connections from the isomonodromic flow on the 
moduli space of monodromy representations and throw a  
bridge between these two dynamics via the Riemann-Hilbert 
correspondence. 
We begin with the isomonodromic flow.   
\subsection{Isomonodromic Flow} \label{subsec:isomonodromy}
Fix a base point $t \in T$ and take the loops 
$\ga_i \in \pi_1(\P^1-D_t,*)$ as in Figure \ref{fig:loops}. 
Let $U$ be a sufficiently small simply-connected neighborhood 
of $t$ in $T$. 
Then, having $\{\ga_i\}$ as common generators, all the 
fundamental groups $\pi_1(\P^1-D_s,*)$ with $s \in U$ are 
identified with the reference group $\pi_1(\P^1-D_t,*)$. 
Passing to moduli spaces of representations, we have 
isomorphisms
\begin{equation} \label{eqn:psits}
\psi_t^s : \R_t(a) \rightarrow \R_s(a) \qquad (s \in U). 
\end{equation}
This means that the fibration $\pi_a : \R(a) \to T$ is 
locally trivial, where a local trivialization over $U$ is 
given by $\psi_t : \R_t(a) \times U \to \R(a)|_{U}, \, 
(\rho, s) \mt \psi_t^s(\rho)$. 
Then there exists the trivial foliation on $\R(a)|_{U}$ 
whose leaves are the slices $\psi(\{\rho\} \times U)$ 
parametrized by $\rho \in \R_t(a)$. 
These local foliations for various simply-connected open 
subsets $U \subset T$ are patched together to form a global 
foliation on $\R(a)$. 
Moreover, patching together various local isomorphisms of 
the form (\ref{eqn:psits}), we can associate to each 
path $\ell$ in $T$ an isomorphism 
\begin{equation} \label{eqn:ell}
\ell_* : \R_t(a) \rightarrow \R_s(a),  
\end{equation}
where $t$ and $s$ are the initial and terminal points 
of $\ell$, respectively.  
Note that the isomorphism $\ell_*$ depends only on the homotopy 
class of the path $\ell$. 
\begin{definition}[Isomonodromic Flow] \label{def:IMF}
The foliation on $\R(a)$ induced from 
the local triviality of the fibration $\pi : \R(a) \to T$ is 
called the {\it $a$-isomonodromic flow} and is denoted by 
$\F_{\IMF(a)}$.   
It is a time-dependent Hamiltonian dynamics in the 
sense of Definition \ref{def:Hamilton}. 
Namely each fiber $\R_t(a)$ is a symplectic manifold, whose 
symplectic structure $\Om_{\R_t(a)}$ will be described 
in \S\ref{subsec:symplectic}, and the isomorphism 
(\ref{eqn:ell}) is a symplectic isomorphism. 
The dynamical system $(\R(a), \F_{\IMF(a)})$ is denoted 
by $\IMF(a)$, whose fundamental $2$-form $\Om_{\R(a)}$ is 
defined by the following conditions:  
\begin{enumerate} 
\item $\Om_{\R(a)}$ is restricted to the symplectic structure 
$\Om_{\R_t(a)}$ on $\R_t(a)$ for every $t \in T$.  
\item $\iota_v \Om_{\R(a)} = 0$ for any 
$\F_{\IMF(a)}$-horizontal vector filed $v$. 
\end{enumerate}
\end{definition}
\begin{definition}[Family of Isomonodromic Flows] 
\label{def:family} 
Relative versions of Definition \ref{def:IMF}:  
\begin{enumerate}
\item There exists a (unique) family 
$\IMF = (\R, \F_{\IMF})$ of isomonodromic flows over $A$, 
where $\F_{\IMF}$ is a relative foliation on the fibration 
$\R \to A$ that restricts to the foliation $\F_{\IMF(a)}$ 
on each fiber $\R(a)$.    
Moreover there exists a relative $2$-form $\Om_{\R}$ 
on $\R$ that restricts to the fundamental $2$-form 
$\Om_{\R(a)}$ on ${\R(a)}$. 
\item By the fiber-product morphism (\ref{cd:FP}), 
$\IMF$ is pulled back to a relative foliation 
$\bIMF = (\bR, \F_{\sc{\bIMF}})$ on $\bR$, with the 
corresponding relative $2$-form $\Om_{\bR}$.   
\end{enumerate} 
\end{definition} 
\par 
Although it is almost trivial from the purely topological 
nature of the isomonodromic flow, the following lemma is 
worth stating explicitly.  
\begin{lemma}[Geometric Painlev\'e Property] \label{lem:pp1}
For each $a \in A$, the isomonodromic flow $\IMF(a)$ has 
geometric Painlev\'e property. 
\end{lemma}
\par 
It is clear from the construction that the Riccati locus 
$\R^{\r}(a)$ and the non-Riccati locus $\R^{\ci}(a)$ are 
stable under the isomonodromic flow $\IMF(a)$. 
\begin{definition}[Riccati/Non-Riccati Flow] \label{def:RF1} 
For each $a \in A$ (actually for $a \in \rh(\Wall)$), 
\begin{enumerate} 
\item the isomonodromic flow $\IMF(a)$ restricted to 
the Riccati locus $\R^{\r}(a)$ is referred to as the 
{\it Riccati flow} and is denoted by $\IMF^{\r}(a)$, 
and  
\item the isomonodromic flow $\IMF(a)$ restricted to the 
non-Riccati locus $\R^{\ci}(a)$ is referred to as the 
{\it non-Riccati flow} and is denoted by $\IMF^{\ci}(a)$. 
\end{enumerate} 
\end{definition}
\subsection{Symplectic Structure on $\mbox{\bm $\R_t(a)$}$}
\label{subsec:symplectic}
The symplectic nature of moduli spaces of monodromy 
representations was first discussed by Goldman \cite{Goldman}. 
It has been used to study Painlev\'e-type equations by 
Iwasaki \cite{Iwasaki1,Iwasaki2}, Hitchin \cite{Hitchin}, 
Kawai \cite{Kawai1,Kawai2}, Boalch \cite{Boalch1} and others. 
Now we recall the topological description of the symplectic 
structure $\Om_{\R_t(a)}$ on the smooth locus of 
$\R_t(a)$ (under some generic condition on $a$).  
A more comprehensive sheaf-cohomological description, 
which allows for every values of $a \in A$, can be 
found in Inaba, Iwasaki and Saito \cite{IIS2}. 
\par
\begin{figure}[t]
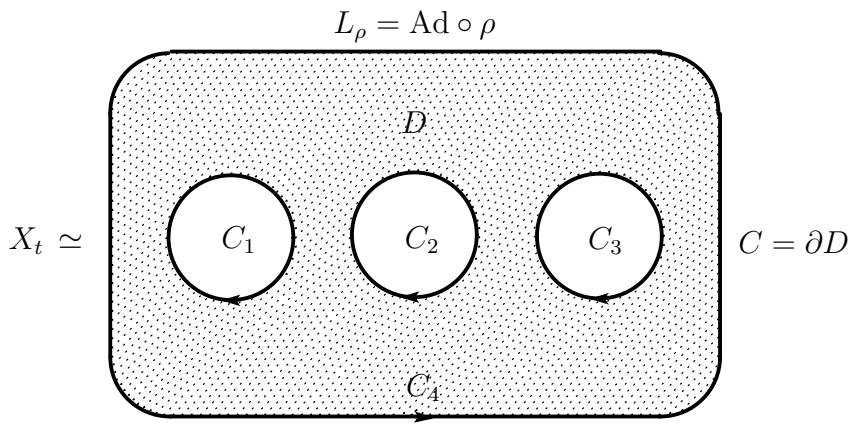

\begin{center}
\unitlength 0.1in
%
\end{center}
\caption{The Riemann sphere with four disks deleted} 
\label{fig:D} 
\end{figure}

 
In general, given a topological space $X$, let $\R(X)$ denote 
the set of all Jordan equivalence classes of 
$SL_2(\C)$-representations of $\pi_1(X)$. 
In stead of using the $4$-punctured Riemann sphere $\P^1-D_t$, 
we employ a homotopically equivalent 
domain\footnote{somewhat confusing notation: $D$ should not 
be confused with $D_t$.} $D$ obtained from $\P^1$ by deleting 
four disjoint open disks centered at $t_1$, $t_2$, $t_3$, $t_4$.   
The boundary $C$ of $D$ consists of four disjoint circles 
$C_1$, $C_2$, $C_2$, $C_4$, (see Figure \ref{fig:D}). 
Then $\R(D)$ is identified with $\R_t = \R(\P^1-D_t)$, while 
$\R(C)$ is identified with $A = \C_a^4$ by the isomorphism 
\[
\R(C) \rightarrow A, \qquad \rho \mt a = 
(\Tr\,\rho(C_1), \, \Tr\,\rho(C_2), \, \Tr\,\rho(C_3), \, 
\Tr\,\rho(C_4)).
\]
Restricting representations of $\pi_1(D)$ to $\pi_1(C)$,  
we have the restriction mapping
\[ 
r : \R(D) \to \R(C), \qquad \rho \mt \rho|_{C}. 
\]
Then $\R_t(a)$ is identified with the fiber $r^{-1}(a)$ over 
$a \in A = \R(C)$ and the Zariski tangent space of 
$\R_t(a)$ at a point $\rho \in \R_t(a)$ is given by 
\[
T_{\rho}\R_t(a) = \Ker\, 
[\,(dr)_{\rho} : T_{\rho} \R(D) \rightarrow T_{r(\rho)} \R(C)\,], 
\]
\par
Let $L_{\rho}$ be the locally constant system on $D$ 
associated to the representation $\Ad\circ\rho$, where 
$\Ad : SL_2(\C) \to GL(\fraksl_2(\C))$ is the adjoint 
representation of $SL_2(\C)$. 
Then the standard infinitesimal deformation theory tells us that 
the Zariski tangent spaces $T_{\rho}\R(D)$ and $T_{r(\rho)} \R(C)$ 
are identified with the first cohomology groups $H^1(D,L_{\rho})$ and 
$H^1(C;L_{\rho})$, and that the tangent map $(dr)_{\rho}$ is 
identified with the homomorphism $j^*$ in the cohomology long 
exact sequence 
\[ 
\begin{CD}
@>>> H^0(C;L_{\rho}) @> \delta^* >> H^1(D,C;L_{\rho}) @> i^* >> 
H^1(D;L_{\rho}) @> j^* >> H^1(C;L_{\rho}) @>>>   
\end{CD}
\]
for the space pair $(D,C)$ with coefficients in $L_{\rho}$. 
Thus we have an isomorphism
\begin{equation} \label{eqn:tangent1}
T_{\rho}\R_t(a) \cong \Ker\, 
[\, j^* : H^1(D;L_{\rho}) \rightarrow H^1(C;L_{\rho}) \,], 
\end{equation}
Moreover the cohomology long exact sequence yields another 
isomorphism induced by $i^*$, 
\begin{equation} \label{eqn:tangent2}
T_{\rho}\R_t(a) \cong 
\dfrac{H^1(D,C;L_{\rho})}{\delta^* H^0(C;L_{\rho})}.   
\end{equation}
\par 
By the Poincar\'e-Lefschetz duality, there exists a  
nondegenerate bilinear form
\[
\begin{CD}
H^1(D;L_{\rho}) \ot H^1(D,C;L_{\rho}) 
@> \scriptstyle{\mathrm{cup \,\, product}} >> 
H^2(D,C;L_{\rho} \ot L_{\rho}) \\
@> \scriptstyle{\mathrm{Killing \,\, form}} >> 
H^2(D,C;\C_{D}) \cong \C, 
\end{CD} 
\]
which induces a nondegenerate pairing between the righthand 
sides of (\ref{eqn:tangent1}) and (\ref{eqn:tangent2}), and hence 
a nondegenerate skew-symmetric bilinear form on the tangent 
space $T_{\rho} \R_t(a)$,  
\[
\Om_{\R_t(a), \, \rho} : 
T_{\rho}\R_t(a) \times T_{\rho}\R_t(a) \to \C. 
\]
In this manner we have obtained an alomost symplectic structure 
$\Om_{\R_t(a)}$ on $\R_t(a)$, which in fact is 
a symplectic structure. 
This fact, namely, the closedness of $\Om_{\R_t(a)}$ is 
trivial in our $4$-point case where $\R_t(a)$ is a surface.  
It can be proved in the general $n$-point situation. 
\subsection{Nonlinear Monodromy of Isomonodromic Flow}
\label{subsec:NMIF}
Given a base point $t \in T$, we consider isomorphisms 
(\ref{eqn:ell}) when the $\ell$'s are loops in $T$ with 
base point at $t$.  
Then they become automorphisms of $\R_t(a)$ and yield a 
group homomorphism
\begin{equation} \label{eqn:NMIF}
\pi_1(T,t) \rightarrow \Aut\,\R_t(a), \qquad 
\ell \mt \ell_*, 
\end{equation}
which is nothing but the nonlinear monodromy of the 
isomonodromic flow $\IMF(a)$ 
(see Definition \ref{def:PRM}). 
This homomorphism can be described in terms of braid 
groups on three strings (see Dubrovin and Mazzocco \cite{DM}, 
Iwasaki \cite{Iwasaki4} and Boalch \cite{Boalch2,Boalch3}).  
To recall this description, we put $t_4$ at infinity and 
redefine the time-variable space as the configuration space 
of distinct ordered three points in $\C$, that is,   
\begin{equation}  \label{eqn:T2} 
T = \{\, t = (t_1,t_2,t_3) \in \C^3\,:\, 
t_i \neq t_j \,\,\, \mbox{for} \,\,\, i \neq j \, \}. 
\end{equation} 
Then the fundamental group $\pi_1(T,t)$ is isomorphic to 
the pure braid group $P_3$ on three strings. 
If $T$ is replaced by the configuration space of distinct 
{\it unordered} three points in $\C$, then $\pi_1(T,t)$ 
is isomorphic to the ordinary braid group $B_3$ on 
three strings (see e.g. Birman \cite{Birman}). 
Recall that there exists the natural exact sequence 
$1 \to P_3 \to B_3 \to S_3 \to 1$, where $S_3$ represents  
the permutations of $t_1$, $t_2$, $t_3$. 
For later convenience, we employ the following terminology.  
\begin{definition}[Full-Monodromy and Half-Monodromy] 
\label{def:half} 
Monodromy in terms of pure braids are referred 
to as {\it full-monodromy}, while monodromy in terms of 
ordinary braids are referred to as {\it half-monodromy},   
respectively. 
\end{definition}
Using half-monodromy is convenient for shorter presentation   
and the full-monodromy is just obtained by restricting 
the ordinary braid group to its pure subgroup. 
The full-monodromy in (\ref{eqn:NMIF}) makes 
sense for each individual $\IMF(a)$, while the 
half-monodromy only makes sense for $\IMF$.  
To describe the half-monodromy, we introduce the following 
natural action of $B_3$ on $\R_t$.  
\begin{definition}[Action of Braids on Representations] 
\label{def:braidaction}
The action of the braid group $B_3$ on the moduli space 
$\R_t$ of monodromy representations,     
\[ 
B_3 \times \R_t \rightarrow \R_t, \qquad 
(\b, \rho) \mt \rho^{\b}, 
\]
is defined by the following condition, which we call 
the {\it global isomonodromy condition}, 
\begin{equation} \label{eqn:gimdc}
\phan{a} \rho^{\b}(\ga^{\b}) = \rho(\ga). 
\end{equation}
Here $\ga \mt \ga^{\beta}$ is the natural 
action of $\b \in B_3$ on $\pi_1(X_t,*)$ defined 
in Definition \ref{def:braidaction2}, where  
\[
X_t = \C - \{t_1,t_2,t_3\}.
\]
\end{definition}
\begin{figure}[p] 
\begin{center}
\unitlength 0.1in
%
\end{center}
\caption{Basic braid $\b_i$, where $(i,j,k)$ is a cyclic 
of $(1,2,3)$} 
\label{fig:braid} 
\vspace{1.6cm}  
\begin{center}  
\unitlength 0.1in
%
%
\end{center} 
\caption {The braid action $\beta_i : (\ga_i,\ga_j,\ga_k) 
\mt (\ga_i',\ga_j',\ga_k')$}
\label{fig:action}
\end{figure}
\begin{definition}[Action of Braids on Fundamental Group]
\label{def:braidaction2}
Let $\b_i$ be the braid as indicated in Figure \ref{fig:braid}, 
where $(i,j,k)$ is any cyclic permutation of $(1,2,3)$. 
Then the braid group $B_3$ is generated by the basic braids 
$\b_1$, $\b_2$, $\b_3$. 
On the other hand the fundamental group $\pi_1(X_t,*)$ is 
the free group generated by the loops $\ga_1$, 
$\ga_2$, $\ga_3$ in Figure \ref{fig:loops}. 
Thus we have 
\[
B_3 = \la \b_1,\b_2,\b_3\ra, \qquad 
\pi_1(X_t,*) = \la \ga_1,\ga_2,\ga_3 \ra. 
\]
In terms of these generators, the action of $B_3$ on 
$\pi_1(X_t,*)$ is given as in Figure \ref{fig:action}. 
Namely the action of the $i$-th basic baraid, 
$\b_i : (\ga_i,\ga_j,\ga_k) \mt (\ga_i',\ga_j',\ga_k')$, 
is expressed as  
\[
\ga_i' = \ga_i^{-1}\ga_j \ga_i, \qquad 
\ga_j' = \ga_i', \qquad \ga_k' = \ga_k.  
\]
where the composition of loops is taken from right to left. 
\end{definition} 
\par 
The following theorem is clear from the manner in which 
the action is defined as in (\ref{eqn:gimdc}).  
\begin{theorem}[Nonlinear Monodromy] \label{thm:NMofIMF} 
The half-monodromy of $\IMF$ is given by the $B_3$-action 
on $\R_t$ defined as in Definition $\ref{def:braidaction}$ 
and the full-monodromy of $\IMF(a)$ is the $P_3$-action on 
$\R_t(a)$ that is the restriction of the $B_3$-action above 
to $P_3$ $(see $Table $\ref{tab:NMofIMF}$$)$.  
\end{theorem} 
\begin{table}[h] 
\begin{center} 
\begin{tabular}{|c|c|}
\hline 
 half-monodromy  & full-monodromy     \\ 
\hline 
 $B_3 \car \R_t$ & $P_3 \car \R_t(a)$ \\
\hline 
\end{tabular}
\end{center}
\caption{Half-monodromy of $\IMF$ and full-monodromy of $\IMF(a)$} 
\label{tab:NMofIMF}    
\end{table}
\par 
As in Dubrovin and Mazzocco \cite{DM} and 
Iwasaki \cite{Iwasaki4}, we make the following remark.   
\begin{remark}[Reduction to Modular Group] \label{rem:modular}
It is well known that the center $Z(B_3)$ of $B_3$ is the 
infinite cyclic group $\la (\b_i\b_j)^3 \ra$ generated by 
$(\b_i\b_j)^3$ and the quotient group $B_3/Z(B_3)$ is 
isomorphic to the full modular group $\G \simeq PSL(2,\Z)$.  
An inspection shows that our braid group action is trivial 
on the center $Z(B_3)$. 
Hence it is reduced to an action of the full modular group 
$\G$ on $\R_t$. 
In view of Remark \ref{rem:reduction}, this reduction is quite 
possible since the fundamental group of $\P^1-\{0,1,\infty\}$ 
is isomorphic to the level-two principal congruence subgroup 
$\G(2)$ of $\G$. 
The resulting modular group action will be described  
explicitly in Definition \ref{def:automorphism}. 
\end{remark} 
\subsection{Painlev\'e Flow} \label{subsec:PF} 
From our point of view, the Painlev\'e flow should be defined 
as the pull-back of the isomonodromic flow by the 
Riemann-Hilbert correspondence. 
This standpoint was first adopted by Iwasaki 
\cite{Iwasaki1,Iwasaki2}, though things were still 
looked at locally.    
Currently a completely global formulation is feasible, 
now that we have such a neat result as in Theorem \ref{thm:RH}.       
\begin{theorem}[Painlev\'e Flow] \label{thm:PF} 
For any $\k \in \K$, put $a = \rh(\k) \in A$.   
\begin{enumerate}
\item There exists a unique holomorphic foliation 
$\F_{\PVI(\k)}$ on $\M(\k)$ such that the 
$\k$-Riemann-Hilbert correspondence 
$\RH_{\k} : \M(\k) \to \R(a)$ gives a quasi-conjugacy map  
\begin{equation} \label{eqn:qRHconj} 
\RH_{\k} : (\M(\k), \, \F_{\PVI(\k)}) \rightarrow 
(\R(a), \, \F_{\IMF(a)}).    
\end{equation} 
The dynamical system $\PVI(\k) = (\M(\k), \F_{\PVI(\k)})$ is 
called the $\k$-Painlev\'e flow.  
\item The quasi-conjugacy map $(\ref{eqn:qRHconj})$ induces 
a conjugacy map  
\[
\RH_{\k} : (\M^{\ci}(\k), \, \F_{\PVI(\k)}) \rightarrow 
(\R^{\ci}(a), \, \F_{\IMF(a)}), 
\]    
when restricted to the Riccati locus.  
\item The fundamental $2$-form $\Om_{\M(\k)}$ for the 
$\k$-Painlev\'e flow $\PVI(\k)$ is the unique holomorophic 
$2$-form on $\M(\k)$ that satisfies the condition    
\[
\Om_{\M(\k)} = \RH_{\k}^* \, \Om_{\R(a)} \qquad \mbox{on} \quad 
\M^{\ci}(\k). 
\] 
\end{enumerate} 
\end{theorem}
\par 
The point of Theorem \ref{thm:PF} is explained in the 
following manner.  
\begin{remark}[Codimension-Two Argument] \label{rem:CTA} 
If $\k \in \K - \Wall$, this theorem immediately follows 
from Corollary \ref{cor:RH2}, since in this case 
there is no Riccati locus and $\RH_{\k} : \M(\k) \to \R(a)$ 
is biholomorphic.  
However, if $\k \in \Wall$, things are not so simple because 
$\RH_{\k}$ fails to be injective on the Riccati locus 
$\M^{\r}(\k)$, which is of codimension one in $\M(\k)$. 
In this case it is not immediately clear as to whether the 
Painlev\'e flow extends to the Riccati locus.     
To avoid this difficulty, we should consider the 
full-Riemann-Hilbert correspondence $\RH : \M \to \bR$ in 
(\ref{cd:RH3}). 
By Definition \ref{def:family} we have relative foliation 
$\F_{\bIMF}$ and relative $2$-form $\Om_{\bR}$ on $\bR$. 
Since $\RH : \M^{\ci} \to \bR^{\ci}$ is biholomorphic 
by Theorem \ref{thm:RH}, $\F_{\bIMF}$ and $\Om_{\bR}$ 
can be pulled back to a holomorphic relative foliation 
$\F_{\PVI}$ and a holomorphic relative $2$-form  
$\Om_{\M}$ on $\M^{\ci}$.    
Since the complement $\M^{\r} = \M - \M^{\ci}$ 
is of codimension two in $\M$ (see Lemma \ref{lem:codim}), 
Hartog's extension theorem implies that $\F_{\PVI}$ and 
$\Om_{\M}$ can be extended to the whole space $\M$ 
holomorophically. 
Restricting these extensions to each $\M(\k)$ yields  
a holomorphic flow $\F_{\PVI(\k)}$ and a holomorphic 
$2$-form $\Om_{\M(\k)}$ on $\M(\k)$. 
These are just what we have been seeking.    
\end{remark}   
\begin{theorem}[Geometric Painlev\'e Property] \label{thm:pp2} 
For any $\k \in \K$ the Painlev\'e flow $\PVI(\k)$ 
enjoys geometric Painlev\'e property. 
\end{theorem}
This theorem readily follows from the geometric Painlev\'e 
property for the isomonodromic flow $\IMF(a)$ with 
$a = \rh(\k)$ (see Lemma \ref{lem:pp1}) and from the fact that 
$\RH_{\k}$ is a quasi-conjugacy map between 
$\PVI(\k)$ and $\IMF(a)$, especially from the properness 
of $\RH_{\k}$ (see Lemma \ref{lem:proper}).  
\par
It is clear that the Riccati locus $\M^{\r}(\k)$ and 
the non-Riccati locus $\M^{\ci}(\k)$ are stable under 
the Painlev\'e flow $\PVI(\k)$. 
As a counterpart of Definition \ref{def:RF1} we make the 
following definition. 
\begin{definition}[Riccati/Non-Riccati Flow] 
\label{def:RF2} For each $\k \in \K$ 
(actually for $\k \in \Wall$), 
\begin{enumerate} 
\item the Painlev\'e flow $\PVI(\k)$ restricted to 
the Riccati locus $\M^{\r}(\k)$ is referred to as the 
{\it Riccati flow} and is denoted by $\PVI^{\r}(\k)$, 
and  
\item the Painlev\'e flow $\PVI(\k)$ restricted to the 
non-Riccati locus $\M^{\ci}(\k)$ is referred to as the 
{\it non-Riccati flow} and is denoted by $\PVI^{\ci}(\k)$. 
\end{enumerate}
\end{definition} 
\par
The assertion (2) of Theorem \ref{thm:PF} is now restated 
as follows. 
\begin{theorem}[Conjugacy for Non-Riccati Flows] 
\label{thm:RHconj} 
For any $\k \in \K$ put $a = \rh(\k) \in A$.  
The Riemann-Hilbert correspondence $\RH_k$ yields a 
conjugacy between the non-Riccati Painlev\'e flow 
$\PVI^{\ci}(\k)$ and the non-Riccati isomonodromic 
flow $\IMF^{\ci}(a)$.  
Inparicular the nonlinear monodromy of $\PVI^{\ci}(\k)$ 
is faithfully represented by that of $\IMF^{\ci}(a)$, 
where the latter is described by Theorem 
$\ref{thm:NMofIMF}$ restricted to the non-Riccati loci.   
\end{theorem} 
\par 
The discussions of this subsection are summarized as follows.  
The Riemann-Hilbert correspondence gives an analytic 
quasi-conjugacy between the Painlev\'e flow and the 
isomonodromic flow. 
It gives an analytic conjugacy in the strict sense 
{\it outside} the Riccati locus, while it collapses the 
Riccati locus to a family of singularities. 
Thus we have almost arrived at the situation described in 
the Guiding Diagram in Figure \ref{fig:guiding}, though 
subtle details on the Riccati flow are not depicted there.  
One point yet to be discussed in Figure \ref{fig:guiding} 
is the isomorphism $\R_t(\k) \simeq \mathcal{S}(\th)$,  
which will be established in Theorem \ref{thm:identify}.      
\subsection{Riccati Flows and Hypergeometric Equations} 
\label{subsec:RF} 
This subsection is devoted to the linearization of 
Riccati-Painlev\'e flows. 
This procedure will clearly explain why Riccati flows 
are called so.    
Throughout this subsection we fix $\k \in \Wall$. 
\par 
The Riccati-Painlev\'e flow $\PVI^{\r}(\k)$ is confined 
in the Riccati locus $\M^{\r}(\k)$. 
By Remark \ref{rem:minustwo} each irreducible component 
$\CC \subset \M^{\r}(\k)$, which is stable under the flow, 
is a family of $(-2)$-curves over $T$ as in 
(\ref{eqn:minustwo}). 
Thus $\PVI^{\r}(\k)$ restricts to a dynamical system on 
the $\P^1$-bundle $\pi_{\k} : \CC \to T$.   
For this, we have the following theorem.  
\begin{theorem}[Hypergeometric Equation] \label{thm:HGE}  
On each irreducible component of $\M^{\r}(\k)$ the 
Riccati-Painlev\'e flow $\PVI^{\r}(\k)$ is linearizable 
in terms of a Gauss hypergeometric equation. 
\end{theorem}
\par 
To understand what this means, we should recall 
the following famous theorem.
\begin{theorem}[Fuchs-Poincar\'e] \label{thm:FP}
Let $F(x,y,z)$ be a polynomial of $(y,z)$ whose 
coefficients are meromorphic functions of $x$ in a domain 
$U \subset \C$. 
Let $g$ be the genus of the algebraic curve 
\[
C_x = \{\, (y,z) \in \C^2 \,:\, F(x,y,z) = 0 \,\} 
\]
at a generic point $x \in U$.   
If the first-order nonlinear differential equation 
\begin{equation} \label{eqn:ADE} 
F(x,y,y') = 0, \qquad y' = dy/dx,   
\end{equation}
has analytic Painlev\'e property, then there exists 
the following trichotomy.   
\begin{enumerate}
\item if $g = 0$ then $(\ref{eqn:ADE})$ can be reduced 
to a Riccati equation 
\begin{equation} \label{eqn:RicE}
y' = a(x) \, y^2 + b(x) \, y + c(x),   
\end{equation} 
\item if $g = 1$ then $(\ref{eqn:ADE})$ can be reduced 
to the differential equation of an elliptic curve 
\[
(y')^2 = 4 \, y^3 - g_2 \, y - g_3, 
\]
\item if $g \ge 2$ then $(\ref{eqn:ADE})$ can be 
solved by algebraic quadratures.    
\end{enumerate} 
\end{theorem}
\par 
This theorem means that first-order dynamical systems with 
Painlev\'e property are classified by the genera of spaces 
of initial-conditions.   
In the case of genus zero, it asserts that the dynamics is 
governed by a Riccati equation. 
The Riccati equation (\ref{eqn:RicE}) is linearized as  
\begin{equation} \label{eqn:LinE} 
y = - \dfrac{1}{a(x)} \dfrac{Y'}{Y}, \qquad 
a(x) \, Y''- \{a'(x) + a(x) \, b(x) \}\, Y' 
+ a^2(x) \, c(x) \, Y = 0. 
\end{equation}    
\par 
Let us return to the situation in Theorem \ref{thm:HGE}, 
where we were discussing the Riccati-Painlev\'e flow 
$\PVI^{\r}(\k)$ restricted to an irreducible component 
$\CC \subset \M^{\r}(\k)$. 
Then obviously we are in the genus-zero case of Theorem 
\ref{thm:FP}\footnote{To apply Theorem \ref{thm:FP}, we 
should recast the bundle $\pi_{\k} : \CC \to T$ to a 
$\P^1$-bundle over $U = \P^1-\{0,1,\infty\}$ by using the 
symplectic reduction in Remark \ref{rem:reduction}.}. 
If we use the coordinate expression of $\PVI(\k)$, we can 
see that in our case the linear equation (\ref{eqn:LinE}) 
is (essentially) a Gauss hypergeometric equation. 
Here the coordinate expression of $\PVI(\k)$ will be given 
in Theorem \ref{thm:HVI}. 
Since using coordinate expressions is not beautiful, 
we may pose the following problem.    
\begin{problem}[Linearization] \label{prob:riccati}   
For each irreducible component $\CC \subset \M^{\r}(\k)$, 
show that there exist a rank-two vector bundle $E$ on 
$(\P^1)^4$ and an integrable connection $\nabla$ 
on it, having regular singularities along the diagonal 
$(\P^1)^4-T$, such that the 
Riccati-Painlev\'e flow $\PVI^{\r}(\k)$ restricted to 
$\CC$ is the flat projective connection induced from the 
flat linear connection $(E, \nabla)|_{T}$.  
\end{problem}   
\par 
Of course we have to solve it conceptually without using 
coordinate expressions. 
In any case, it is now clear that 
classifying Riccati solutions amounts to classifying 
irreducible components of Riccati loci. 
We may consider this problem at a fixed $t \in T$. 
So the problem is to classify $(-2)$-curves on  
moduli spaces $\M_t(\k)$, $\k \in \Wall$. 
Originally, the relation between Riccati solutions to  
Painlev\'e equations and $(-2)$-curves on spaces of 
initial-conditions was clarified by Saito and 
Terajima \cite{STe} and Sakai \cite{Sakai}. 
In particular Saito and Terajima gave a complete 
classification of $(-2)$-curves. 
We can now amplify their viewpoint by the picture 
of resolution of singularities by Riemann-Hilbert 
correspondence. 
To do so we should pose the following problem. 
\begin{problem}[Classification of 
$\mbox{\bm $(-2)$}$-Curves] \label{prob:minustwo2} 
Given any $(t,\k) \in T \times \Wall$, classify all 
$(-2)$-curves on $\M_t(\k)$ in connection with the 
resolution of singularities by the Riemann-Hilbert 
correspondence $\RH_{t,\k} : \M_t(\k) \to \R_t(a)$. 
\end{problem}  
\par 
This problem will be settled in 
Theorem \ref{thm:classification}.  
This subsection is closed with the following:  
\begin{remark}[History] \label{rem:history}
Attempts at generalizing Theorem \ref{thm:FP} to second-order 
equations led Painlev\'e to discover his famous equations.
\end{remark}  
\section{Family of Affine Cubic Surfaces}  \label{sec:cubic} 
All the constructions described so far can be made more 
explicit if we consider a family of affine cubic surfaces 
defined as a certain categorical quotient. 
We present the construction of the family, following 
Iwasaki \cite{Iwasaki4}. 
Throughout this section we fix a time $t \in T$.  
\subsection{Categorical Quotient} 
Let $\Hom_t = \Hom(\pi_1(\P^1-D_t,*), SL_2(\C))$ be the set 
of all representations of $\pi_1(\P^1-D_t,*)$ into $SL_2(\C)$. 
Then $\Hom_t$ is naturally an affine algebraic variety and 
admits the adjoint action 
\[
\Ad : SL_2(\C) \times \Hom_t \to \Hom_t, \qquad (P, \rho) 
\mt \Ad(P) \rho, 
\]
defined by $(\Ad(P) \rho)(\ga) = P \rho(\ga) P^{-1}$ for 
$\ga \in \pi_1(\P^1-D_t,*)$. 
It is known that the moduli space $\R_t$ of Jordan equivalence 
classes of representations is isomorphic to the categorical 
quotient 
\[
\Hom_t/\!/\Ad = \Spec\,\C[\Hom_t]^{\Ad}. 
\]
If the generators $\ga_i$ of $\pi_1(X_t,*)$ are chosen as 
in Figure \ref{fig:loops}, then $\Hom_t$ can be identified with 
\[ 
\RR = \{\, M = (M_1,M_2,M_3,M_4) \in SL_2(\C)^4  
\,:\, M_4M_3M_2M_1 = I \,\}, 
\]
through the map $\Hom_t \to \RR, \,\, \rho \mt M$ defined 
by $M_i = \rho(\ga_i)$. 
With this identification, the moduli space of representations 
$\R_t$ is isomorphic to the categorical quotient 
\begin{equation} \label{eqn:catquot}
\RR/\!/\Ad = \mathrm{Spec}\, \C[\RR]^{\Ad},  
\end{equation}
where $\Ad$ represents the diagonal adjoint action of 
$SL_2(\C)$ on $\RR$. 
\par
The invariant ring $\C[\RR]^{\Ad}$ has generators 
$(x,a) = (x_1,x_2,x_3, a_1,a_2,a_3,a_4)$ given by
\[
\left\{ 
\begin{array}{rcll}
x_i &=& \Tr(M_jM_k) \,\, & (\{i,j,k\}=\{1,2,3\}), \\[2mm] 
a_i &=& \Tr\, M_i   \,\, & (i = 1,2,3,4). 
\end{array}
\right.
\]
Note that $a = (a_1,a_2,a_3,a_4) \in A$ is just the local 
monodromy data defined in (\ref{eqn:a}). 
We may refer to $x = (x_1,x_2,x_3) \in \C^3_x$ as the 
{\it global monodromy data}, since $x_i$ comes from the  
monodromy matrix $M_jM_k$ along the global loop $\ga_i\ga_j$ 
surrounding the two points $t_i$ and $t_j$ simultaneously. 
The generators $(x,a)$ have only one algebraic relation 
$f(x,\th(a)) = 0$, where $f(x,\th)$ is the cubic polynomial 
of $x$ with coefficients $\th = (\th_1,\th_2,\th_3,\th_4)$ 
defined by 
\begin{equation} \label{eqn:f}
f(x,\th) = x_1x_2x_3+x_1^2+x_2^2+x_3^2 
-\th_1 x_1 -\th_2 x_2 -\th_3 x_3 + \th_4, 
\end{equation}
which can be found in the book \cite{FK} 
of Fricke and Klein. 
In terms of local monodromy data $a \in A$, the coefficients 
$\th = \th(a)$ are expressed as 
\begin{equation}  \label{eqn:theta}
\th_i = 
\left\{
\begin{array}{ll}
a_i a_4 + a_j a_k \qquad & (i = 1,2,3), \\[2mm]
a_1 a_2 a_3 a_4 + a_1^2 + a_2^2 + a_3^2 + a_4^2 - 4  
\qquad & (i=4). 
\end{array}
\right. 
\end{equation}
Let $\Th := \C^4_{\th}$ denote the complex $4$-space 
parametrizing the coefficients $\th$ of $f(x,\th)$. 
\subsection{Correspondences of Parameters}
So far we have encountered three kinds of parameters, 
that is, the parameters $\k \in \K$ of $\PVI$, which is 
nothing but the exponents of parabolic connections; 
the local monodromy data $a \in A$; and the coefficients 
$\th \in \Th$ of the cubic polynomial $f(x,\th)$. 
Relations among them are depicted in 
Table \ref{tab:parameters}, where $\k \mt a$ is 
given by (\ref{eqn:rh2}) and $a \mt \th$ is given 
by (\ref{eqn:theta}), respectively. 
\begin{table}[h] 
\begin{center}  
\begin{tabular}{ccccc}
{parameters of} && {local mono-} && {parameters} \\[1mm]
{Painlev\'e VI} && {dromy data}  && {of cubics} \\[4mm] 
$\K$ & $\lra$ & $A$  & $\lra$ & $\Th$  \\
\rot{90}{$\in$} & & \rot{90}{$\in$} & & \rot{90}{$\in$} \\
$\k$ & $\lmt$ & $a$  & $\lmt$ & $\th$ 
\end{tabular} 
\end{center}
\caption{Correspondences of parameters}
\label{tab:parameters}
\end{table}
\par 
In many respects the parameters of cubics $\th \in \Th$ 
are more essential than the local monodromy data $a \in A$. 
One reason lies in the following observation due 
to Terajima \cite{Terajima}.  
\begin{lemma}[Basis of $\mbox{\bm $\W$}$-Invariants] 
\label{lem:basis} 
As a function of exponents $\k \in \K$, the coefficients 
$\th = (\th_1,\th_2,\th_3,\th_4)$ of the cubic polynomial 
$f(x,\th)$ form a basis of $\W$-invariants. 
\end{lemma}
Here, by $\th$ being a basis of $\W$-invariants, we mean that 
any $\W$-invariant entire functions on $\K$ is an entire 
function of $\th$.  
So far the map $\rh : \K \to A$ in (\ref{eqn:rh3}) has been 
called the Riemann-Hilbert correspondence in the parameter 
level (see Definition \ref{def:rh0}).   
However, in view of Lemma \ref{lem:basis} the following 
revised definition would be more appropriate. 
\begin{definition}[Riemann-Hilbert Correspondence in 
Parameter Level] \label{def:rh4} 
From now on the composite $\K \to \Th$ of two maps 
$\K \to A$ in (\ref{eqn:rh3}) and $A \to \Th$ in 
(\ref{eqn:theta}) is referred to as the 
{\it Riemann-Hilbert correspondence in the parameter level}. 
Hereafter we write $\rh : \K \to \Th$. 
\end{definition}  
\subsection{Family of Affine Cubic Surfaces} 
\label{subsec:cubicsurface}
The cubic equation $f(x,\th) = 0$ defines a family of 
affine cubic surfaces, that is, the variety 
\[
\mathcal{S} = 
\{\,(x,\th)\in \C_x^3 \times \Th \,:\, f(x,\th) = 0 \,\},   
\]
together with the projection $\pi : \mathcal{S} \to \Th$, 
$(x,\th) \mt \th$. 
The previous discussions imply that the categorical quotient 
$\RR/\!/\Ad$ in (\ref{eqn:catquot}) is realized as the fiber 
product of $\mathcal{S}$ and $A$ over $\Th$ relative to 
the natural projections $\pi : \mathcal{S} \to \Th$ and 
$A \to \Th$. 
Namely we have isomorphisms 
\begin{equation}  \label{eqn:RtS}
\R_t \simeq \RR/\!/\Ad \simeq \mathcal{S} \times_{\Th} A.   
\end{equation}
We write $f_{\th}(x) = f(x,\th)$ regarding it as a polynomial 
of $x$ depending on parameters $\th$. 
For each $\th \in \Th$ the fiber of 
$\pi : \mathcal{S} \to \Th$ over $\th$ is an affine cubic 
surface 
\[
\mathcal{S}(\th) =
\{\, x\in \C^3\,:\, f_{\th}(x) = 0 \,\}. 
\] 
Now (\ref{eqn:RtS}) means that $\R_t(a)$ is isomorphic to 
the cubic surface $\mathcal{S}(\th)$ provided that $\th$ 
is given by (\ref{eqn:theta}) in terms of $a$. 
Thus we have the following definition.  
\begin{definition}[Reformulation of $\mbox{\bm $\RH$}$] 
\label{redef:RH} 
The isomorphism (\ref{eqn:RtS}) and Definition \ref{def:rh4} 
allow us to reformulate the $t$-Riemann-Hilbert correspondence 
(\ref{cd:RH1}) as the commutative diagram 
\[
\begin{CD}
\M_t @> \RH_t >> \mathcal{S} \\
      @V \pi_t VV    @VV \pi V  \\
    \K  @>> \rh > \Th 
\end{CD}
\] 
In a similar manner the $(t,\k)$-Riemann-Hilbert correspondence 
(\ref{eqn:RH4}) is reformulated as 
\begin{equation}  \label{eqn:RHSth}  
\RH_{t,k} : \M_t(\k) \rightarrow \mathcal{S}(\th),  
\qquad \th = \rh(\k). 
\end{equation}
\end{definition} 
\begin{definition}[Poincar\'e Residue] \label{PR}
A natural symplectic structure or an area form on the cubic 
surface $\mathcal{S}(\th)$ is the {\it Poincar\'e residue} 
defined by  
\[
\omega_{\th} = 
\dfrac{dx_i \wedge dx_j}{(\partial f_{\th}/\partial x_k)}, 
\] 
where $(i,j,k)$ is any cyclic permutation of 
$(1,2,3)$; it does not depends on the cyclic permutation 
chosen. 
The smooth and singular loci of $\mathcal{S}(\th)$ 
are denoted by $\mathcal{S}^{\ci}(\th)$ and 
$\mathcal{S}^{\s}(\th)$, respectively.   
Then the Poincar\'e residue $\omega_{\th}$ is holomorphic 
on $\mathcal{S}^{\ci}(\th)$, having singularities along  
$\mathcal{S}^{\s}(\th)$.  
\end{definition}
\par
A complete characterization of the singular locus 
$\mathcal{S}^{\s}(\th)$ will be 
presented in Section \ref{sec:riccati}. 
The following theorem is found in \cite{Iwasaki4}.  
\begin{theorem}[Moduli of Representations and Cubic Surface] 
\label{thm:identify} 
For any $(t,a) \in T \times A$, let 
$\th = \th(a)$ be defined by $(\ref{eqn:theta})$. 
Then there exists an identification of symplectic manifolds 
\begin{equation}  \label{eqn:identify} 
i \,:\,\, (\R_t(a), \, \Om_{\R_t(a)}) \simeq 
(\mathcal{S}(\th), \, \omega_{\th}). 
\end{equation} 
\end{theorem}
\par
The main ingredient of the proof is the de Rham theorem. 
At the end of this section, the following remark would be 
of some interests. 
\begin{remark}[Moduli of Cubic Surfaces] \label{rem:cayley} 
It is well known in classical algebraic geometry that the 
isomorphism classes of cubic surfaces in $\P^3$ have a 
$4$-dimensional moduli space and that there exists a 
$4$-parameter family of general cubic surfaces, known as 
Cayley's normal form \cite{Cayley}. 
Some computations imply that our family $\mathcal{S}$ and 
Cayley's normal form, as modified by Naruki and Sekiguchi 
\cite{Naruki, NS}, have a common algebraic cover and hence 
our family captures general moduli (see Iwasaki \cite{Iwasaki3} 
and Terajima \cite{Terajima}). 
Thus the family $\mathcal{S}$ can be taken as another normal 
form than Cayley's.  
It is remarkable that general sixth Painlev\'e equations are 
connected with general cubic surfaces through the 
Riemann-Hilbert correspondence.   
\end{remark} 
\section{B\"acklund Transformations} \label{sec:backlund} 
Symmetries of the Painlev\'e equation are called 
B\"acklund transformations. 
Na\"{\i}vely, a B\"acklund transformation is a birational 
transformation that converts one Painlev\'e equation $\PVI(\k)$ 
to another Painlev\'e equation $\PVI(\k')$, where $\k$ and 
$\k'$ may differ. 
More precisely, it is a birational map from one phase space 
$\M(\k)$ to another phase space $\M(\k')$ that commutes with 
the Painlev\'e flows.  
There are at least two approaches to understand B\"acklund 
transformations. 
\begin{remark}[Two Approaches to B\"acklund Transformations] 
\label{rem:backlund1} 
$\phan{a}$ 
\begin{enumerate}
\item birational canonical transformations, 
\vspace{-0.2cm}
\item covering transformations of the Riemann-Hilbert 
correspondence. 
\end{enumerate}
\end{remark}
\par 
The first approach (1) has been employed by such authors  
as Lukashevich and Yablonski \cite{LY} Fokas and 
Ablowitz \cite{FA} and Okamoto \cite{Okamoto4} in the 
style of explicit calculations.    
In particular, Okamoto discovered that $\PVI$ 
admits affine Weyl group symmetries of type $D_4^{(1)}$. 
He expressed them as birational canonical transformations 
of Hamiltonian systems;  
Noumi and Yamada \cite{NY} systematized them by a symmetric 
description in terms of a new Lax pair;  
Arinkin and Lysenko \cite{AL2} geometrized them as isomorphisms 
between moduli spaces of $SL_2(\C)$-connections;  
Sakai \cite{Sakai} also geometrized them in his Cremona 
framework;  
Saito and Umemura \cite{SU} characterized them as flops. 
In fact, these various viewpoints are too diverse to be 
tagged with the same label.   
\par
Nonetheless, they still have a common feature to the effect 
that they look only on a one-side of the Riemann-Hilbert 
correspondence, that is, the moduli space of parabolic 
connections or its relatives (spaces where $\PVI$ is defined), 
with no attentions to the moduli space of representations.      
On the other hand, the second approach (2) is interested 
in the interaction between the source space and the target 
space of Riemann-Hilbert correspondence, asking what the 
B\"acklund transformations look like through the 
telescope of Riemann-Hilbert correspondence.  
In this section we take approach (2), following the 
exposition of Inaba, Iwasaki and Saito \cite{IIS1}. 
\par 
Take any $\k \in \K$ and put $\th = \rh(\k) \in \Th$. 
Given any element $\si \in \W$, 
we consider the Riemann-Hilbert correspondences 
(\ref{eqn:RHSth}) for the parameter $\k$ and for its 
$\si$-translate $\si(\k)$, 
\[
\RH_{t,\k} : \M_t(\k) \to \mathcal{S}(\th), \qquad 
\RH_{t,\si(\k)} : \M_t(\si(\k)) \to 
\mathcal{S}(\si(\th)). 
\]
By our solution to the Riemann-Hilbert problem 
(see Theorem \ref{thm:RH3}), both of them are biholomorphic 
maps  if $\k \not\in \Wall$, and are minimal resolutions of 
singularities if $\k \in \Wall$, respectively; in any case 
they are bimeromorhpic morphisms. 
On the other hand, by the $\W$-invariance of $\th$ 
(see Lemma \ref{lem:basis}), we have $\th = \si(\th)$ and 
hence the cubic surfaces $\mathcal{S}(\th)$ and 
$\mathcal{S}(\si(\th))$ are identical. 
Therefore there exists a unique bimeromorphic map 
\begin{equation} \label{eqn:backlund}
s_{\si} : \M_t(\k) \rightarrow \M_t(\si(\k)) 
\end{equation}
that makes the diagram in Figure \ref{fig:backlund} 
commute, that is, the unique lift that covers the 
{\it identity} on 
$\mathcal{S}(\th) = \mathcal{S}(\si(\th))$ through 
the Riemann-Hilbert correspondence. 
\begin{figure}[t]
\[
\begin{CD}
\M_t(\k) @> \phan{aa} s_{\si} \phan{aa} >> \M_t(\si(\k)) \\
  @V \RH_{t,\k} VV               @VV \RH_{t,\si(\k)} V \\
\mathcal{S}(\th) @>> \mathrm{\sc identity} > 
\mathcal{S}(\si(\th))
\end{CD}
\]
\caption{B\"acklund transformations} \label{fig:backlund}
\end{figure}
In the sprit of approach (2) it is natural to define the concept 
of B\"acklund transformations in the following manner. 
\begin{definition}[B\"acklund Transformations] 
\label{def:backlund} 
By a {\it B\"acklund transformation}, we mean the lift 
$s_{\si}$ of an element $\si \in \W$ as 
in (\ref{eqn:backlund}). 
The group of B\"acklund transformations is, by definition, the 
group consisting of all those lifts $s_{\si}$ with 
$\si \in \W$, 
\[
G = \la \, s_{\si} \,|\, \si \in \W \, \ra = 
\la s_0, s_1, s_2, s_3, s_4 \ra \simeq \W. 
\]
where $s_i$ is the lift of the basic reflection 
$\si_i$ for $i = 0,1,2,3,4$ (see (\ref{eqn:sigma}) and 
Figure \ref{fig:dynkin}); we refer to $s_i$ as the $i$-th 
{\it basic} B\"acklund transformation. 
For each $t \in T$ the group $G$ acts on the moduli space 
$\M_t$ over $\K$ in such a manner that there exists the 
commutative diagram 
\[
\begin{CD}
\M_t @> G >> \M_t \\
     @V \pi_t VV @VV \pi_t V \\
\K @>> \W > \K 
\end{CD}
\]
\end{definition}
\begin{remark}[Advantage and Disadvantage] \label{rem:AD}
The advantage of Definition \ref{def:backlund} is clear.
\begin{enumerate} 
\item The character of B\"acklund transformations is  
transparent, as the covering transformations of the 
Riemann-Hilbert correspondence, in the sense 
of Figure \ref{fig:backlund}. 
\item The origin of the affine Weyl group structure 
of B\"acklund transformations is clear: 
It just comes from the fact that the Riemann-Hilbert 
correspondence in the parameter level $\rh : \K \to \Th$ 
is a branched $\W$-covering.  
\end{enumerate} 
There is also a disadvantage of this definition. 
\begin{enumerate}
\item[(3)] The birational character of B\"acklund 
transformations is far from trivial, because our 
definition makes use of the Riemann-Hilbert 
correspondence which is highly transcendental. 
From Definition \ref{def:backlund} we only know that 
B\"acklund transformations are bimeromorphic, while   
their birationality is {\sl a priori} clear from the 
viewpoints of approach (1).     
\end{enumerate} 
\end{remark}  
\par 
So we are obliged to discuss the relation between these two 
approaches and to unify them.   
In this respect, Inaba, Iwasaki and Saito \cite{IIS1} proved 
the following result.       
\begin{theorem}[Coincidence of Two Approaches] 
\label{thm:backlund}
The two approaches in Remark $\ref{rem:backlund1}$ coincide. 
Namely the B\"acklund transformations in the sense of 
Definition $\ref{def:backlund}$ are exactly those which 
have been known as the birational canonical transformations.  
\end{theorem}  
\par 
We remark that a different proof of this theorem was given 
later by Boalch \cite{Boalch2}. 
There exists an explicit formula for the basic B\"acklund 
transformations $s_i$ in terms of certain canonical 
coordinates on $\M(\k)$ (see Theorem \ref{thm:BVI}).  
We can calculate the lift $s_i$ of $\si_i$, overcoming 
the transcendental nature of the Riemann-Hilbert 
correspondence (see \cite{IIS1}). 
As a matter of fact, $s_1$, $s_2$, $s_3$, $s_4$ are easy 
to handle and the true difficulty lies in the treatment 
of $s_0$. 
As for this the following remark might be helpful.  
\begin{remark}[Gauge Transformations] \label{rem:gauge} 
Let $W'$ be the subgroup of $\W$ 
stabilizing the local monodromy data $a = a(\k)$ as 
a function of $\k$ (see (\ref{eqn:rh3}) and (\ref{eqn:rh2})).  
The subgroup $G' \subset G$ corresponding to $W'$ is called 
the group of {\it gauge transformations}.  
Note that a B\"acklund transformation is a gauge 
transformation if and only if it does not change monodromy.    
\begin{enumerate}
\item $s_1$, $s_2$, $s_3$, $s_4$ are very simple gauge 
transformations; see e.g. \cite{IIS1}.  
\vspace{-0.2cm} 
\item $s_0$ is {\it not} a gauge transformation. 
Arinkin and Lysenko \cite{AL2} described it in 
the level of {\it abstract} isomorphism between 
moduli spaces of $SL(2)$-connections. 
Boalch \cite{Boalch2} was able to characterize it 
as a {\it concrete} transformation of $2 \times 2$ 
Fuchsian systems (\ref{eqn:schlesinger}), passing 
through $3 \times 3$ irregular singular systems via 
Fourier transformation. 
It is desirable to realize it as a natural 
transformation of moduli functors of stable 
parabolic connections. 
\end{enumerate} 
\end{remark}
\par   
As is mentioned in \cite{Okamoto4, AL2, NY, Boalch3}, 
the following remark should be made at this stage. 
\begin{remark}[Extended Affine Weyl Group] 
\label{rem:backlund2}
The affine Weyl group symmetry can be enlarged to an 
{\it extended} affine Weyl group symmetry, if we 
allow some permutations of time variables 
$t = (t_1,t_2,t_3,t_4)$. 
Let $\Kl \simeq \Z_2 \times \Z_2$ be Klein's 
$4$-group of permutations 
\[
\Kl = \{\, 1, (12)(34), \, (13)(24), \, (14)(23) \,\} 
\subset S_4, 
\]
acting on $T \times \K$ by permuting their components.  
Note that the semi-direct product 
\[ 
\Kl \ltimes \W = \widetilde{W}(D_4^{(1)}). 
\] 
is the extended affine Weyl group of type $D_4^{(1)}$ acting 
on $T \times \K$. 
This action is lifted to the moduli space $\M$ over 
$T \times \K$ through the Riemann-Hilbert correspondence, as
\[ 
\begin{CD}
\M @> \widetilde{G} >> \M \\
     @V \pi VV @VV \pi V \\
T \times \K @>> \widetilde{W}(D_4^{(1)}) > T \times \K. 
\end{CD}
\]
\end{remark}
The original independent variable of $\PVI$ in (\ref{eqn:PVI}), 
namely, the cross ratio $x$ in (\ref{eqn:x}) is $\Kl$-invariant 
and hence remains invariant under the 
$\widetilde{W}(D_4^{(1)})$-symmetry. 
Moreover, if one allows the full $S_4$-permutations of 
$t = (t_1,t_2,t_3,t_4)$, one gets the {\it full}  
$W(F_4^{(1)})$-{\it symmetry} of $\PVI$. 
Its connection with the mapping class group of the $4$-holed 
sphere is discussed in Remark 12 of Boalch \cite{Boalch2}. 
\section{Nonlinear Monodromy} \label{sec:nlmonod}
The nonlinear monodromy of the Painlev\'e flow,   
or more precisely, that of the non-Riccati Painlev\'e flow,   
can be represented explicitly in terms of a certain 
modular group action on cubic surfaces. 
In this seciton we are concerned with this description.     
\subsection{Modular Group Action} \label{subsec:modular}
The action is first defined on the ambient space 
$\C^7 = \C^3_x \times \Th$ and then restricted to 
$\mathcal{S}$.  
In what follows $(i,j,k)$ stands for any cyclic permutation 
of $(1,2,3)$. 
We start with the symmetric group $S_3$ of degree 
$3$ acting on $\Th$ by permuting the first three 
components of $\th = (\th_1,\th_2,\th_3,\th_4)$.  
If $\tau_i = (ij) \in S_3$ denotes the 
transposition\footnote{Note that $\tau_3$ is not $(34)$  
but $(31)$. } of $\th_i$ and $\th_j$, then $S_3$ is 
generated by $\tau_1$, $\tau_2$, $\tau_3$, 
\[
S_3 = \la \tau_1, \, \tau_2, \, \tau_3 \ra. 
\]
Next we introduce a lift of $\tau_i$ to $\C^7$ relative 
to the second projection $\C^7 = \C^3_x \times \Th \to \Th$.      
\begin{definition}[Group of Polynomial Automorphisms] 
\label{def:automorphism}
For each $i = 1,2,3$, let 
\[
g_i : \C^7 \to \C^7, \quad (x,\th) \mt (x',\th')
\]
be the polynomial automorphism defined by the formula  
\[
(x_i',\, x_j',\, x_k',\, 
\th_i',\, \th_j',\, \th_k',\, \th_4') = 
(\th_j - x_j - x_k x_i,\, x_i,\, x_k,\, 
\th_j,\, \th_i,\, \th_k,\, \th_4). 
\]
Moreover let $G$ denote the transformation 
group\footnote{The group of B\"acklund transformations is 
also denoted by $G$ in \S\ref{sec:backlund}, but 
no confusion might occur.} generated by $g_1$, $g_2$, $g_3$, 
that is, 
\[  
G = \la g_1, \, g_2, \, g_3 \ra.
\]
\end{definition} 
\begin{figure}[t]
\begin{center}
\unitlength 0.1in
\begin{picture}(48.20,30.90)(5.10,-32.30)
%
\special{pn 20}%
\special{pa 1930 420}%
\special{pa 5130 420}%
\special{pa 5130 2220}%
\special{pa 1930 2220}%
\special{pa 1930 420}%
\special{fp}%
%
\special{pn 20}%
\special{pa 1930 3010}%
\special{pa 5130 3010}%
\special{fp}%
%
\special{pn 20}%
\special{pa 3120 420}%
\special{pa 3120 2220}%
\special{fp}%
%
\special{pn 20}%
\special{pa 4320 410}%
\special{pa 4320 2220}%
\special{fp}%
%
\special{pn 13}%
\special{pa 4320 2210}%
\special{pa 4320 3020}%
\special{da 0.070}%
%
\special{pn 13}%
\special{pa 3120 2220}%
\special{pa 3120 3020}%
\special{da 0.070}%
%
\special{pn 20}%
\special{pa 3720 2380}%
\special{pa 3720 2820}%
\special{fp}%
\special{sh 1}%
\special{pa 3720 2820}%
\special{pa 3740 2753}%
\special{pa 3720 2767}%
\special{pa 3700 2753}%
\special{pa 3720 2820}%
\special{fp}%
\put(29.2000,-3.2000){\makebox(0,0)[lb]{$\mathcal{S}(\th)$}}%
\put(40.3000,-3.1000){\makebox(0,0)[lb]{$\mathcal{S}(\tau(\th))$}}%
\put(38.0000,-26.5000){\makebox(0,0)[lb]{$\pi$}}%
\put(53.1000,-30.8000){\makebox(0,0)[lb]{$\Th$}}%
\put(53.3000,-14.3000){\makebox(0,0)[lb]{$\mathcal{S}$}}%
\put(30.5000,-33.9000){\makebox(0,0)[lb]{$\th$}}%
\put(42.1000,-34.0000){\makebox(0,0)[lb]{$\tau(\th) = \th'$}}%
%
\special{pn 20}%
\special{sh 0.600}%
\special{ar 3120 3010 81 81  0.0000000 6.2831853}%
%
\special{pn 20}%
\special{sh 0.600}%
\special{ar 4310 3000 81 81  0.0000000 6.2831853}%
%
\special{pn 20}%
\special{sh 0.600}%
\special{ar 3120 1900 81 81  0.0000000 6.2831853}%
%
\special{pn 20}%
\special{sh 0.600}%
\special{ar 4320 880 81 81  0.0000000 6.2831853}%
\put(28.2000,-7.3000){\makebox(0,0)[lb]{$\omega_{\th}$}}%
\put(44.0000,-20.7000){\makebox(0,0)[lb]{$\omega_{\tau(\th)}$}}%
%
\special{pn 13}%
\special{ar 4750 2630 1771 1771  3.6801001 4.3391741}%
%
\special{pn 13}%
\special{pa 4090 980}%
\special{pa 4160 950}%
\special{fp}%
\special{sh 1}%
\special{pa 4160 950}%
\special{pa 4091 958}%
\special{pa 4111 971}%
\special{pa 4107 995}%
\special{pa 4160 950}%
\special{fp}%
\special{pa 4160 950}%
\special{pa 4160 950}%
\special{fp}%
\put(33.9000,-12.7000){\makebox(0,0)[lb]{$g$}}%
\put(5.1000,-14.8000){\makebox(0,0)[lb]{$G(2) \,\,\, \subset \,\,\, G \,\,\, \car$}}%
\put(5.1000,-8.1000){\makebox(0,0)[lb]{$\G(2) \,\,\, \subset \,\,\, \G$}}%
\put(5.4000,-30.7000){\makebox(0,0)[lb]{$\{1\} \,\,\, \subset \,\,\, S_3 \,\,\, \car$}}%
%
\special{pn 13}%
\special{pa 720 860}%
\special{pa 720 1170}%
\special{fp}%
\special{sh 1}%
\special{pa 720 1170}%
\special{pa 740 1103}%
\special{pa 720 1117}%
\special{pa 700 1103}%
\special{pa 720 1170}%
\special{fp}%
%
\special{pn 13}%
\special{pa 720 1630}%
\special{pa 720 2710}%
\special{fp}%
\special{sh 1}%
\special{pa 720 2710}%
\special{pa 740 2643}%
\special{pa 720 2657}%
\special{pa 700 2643}%
\special{pa 720 2710}%
\special{fp}%
%
\special{pn 13}%
\special{pa 1305 850}%
\special{pa 1305 1160}%
\special{fp}%
\special{sh 1}%
\special{pa 1305 1160}%
\special{pa 1325 1093}%
\special{pa 1305 1107}%
\special{pa 1285 1093}%
\special{pa 1305 1160}%
\special{fp}%
%
\special{pn 13}%
\special{pa 1305 1650}%
\special{pa 1305 2730}%
\special{fp}%
\special{sh 1}%
\special{pa 1305 2730}%
\special{pa 1325 2663}%
\special{pa 1305 2677}%
\special{pa 1285 2663}%
\special{pa 1305 2730}%
\special{fp}%
\put(21.9000,-12.3000){\makebox(0,0)[lb]{area form}}%
\put(25.4000,-19.8000){\makebox(0,0)[lb]{$(x,\th)$}}%
\put(37.2000,-8.4000){\makebox(0,0)[lb]{$(x',\th')$}}%
%
\special{pn 8}%
\special{pa 2720 1000}%
\special{pa 2850 760}%
\special{fp}%
\special{sh 1}%
\special{pa 2850 760}%
\special{pa 2801 809}%
\special{pa 2825 807}%
\special{pa 2836 828}%
\special{pa 2850 760}%
\special{fp}%
\end{picture}%
\end{center}
\caption{A total picture of the modular group action} 
\label{fig:total}  
\end{figure}
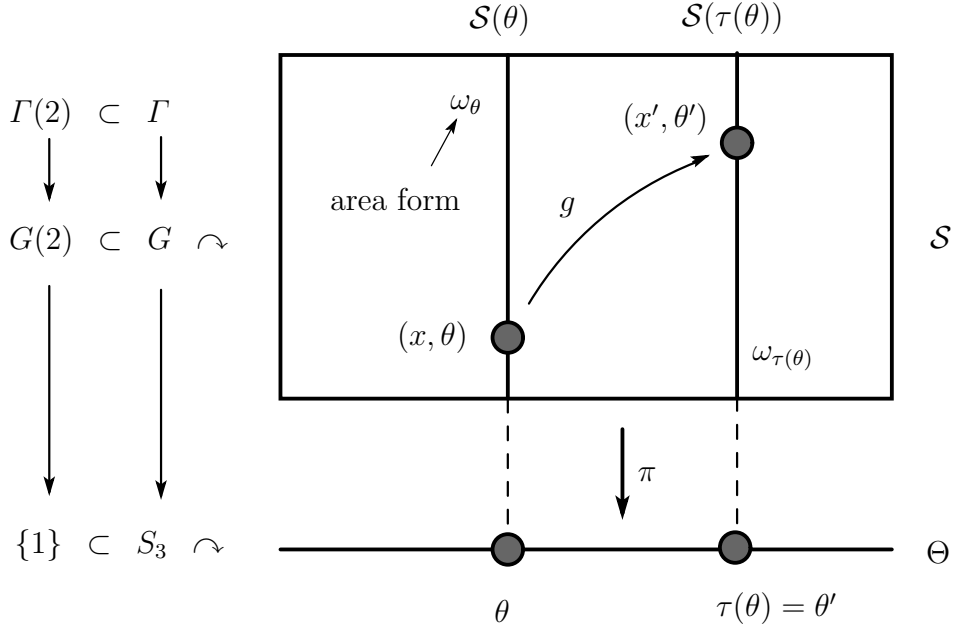
\par
A direct check shows that the generators satisfy three 
relations
\[
g_i g_j g_i = g_j g_i g_j, \qquad (g_i g_j)^3 = 1, \qquad 
g_k = g_i g_j g_i^{-1},   
\] 
which are exactly the defining relations of the full 
modular group
\[ 
\G = PSL_2(\Z)  
= \left\{\, z \mt \dfrac{az+b}{cz+d}\,:\, 
a, \ b, \ c, \ d \in \Z, \,\, ad-bc = 1 \, \right\}. 
\]
Hence there exists a group homomorphism $\G \to G$, through 
which the modular group $\G$ acts on $\C^7$. 
This action is restricted to the principal congruence 
subgroup of level $2$,
\[
\G(2) 
= \left\{\, z \mt \dfrac{az+b}{cz+d} \in \G 
\,:\, a \equiv d \equiv 1, \, b \equiv c \equiv 0 \,\,  
(\mathrm{mod} \,\, 2) \, \right\}. 
\]
The subgroup of $G$ corresponding to $\G(2) \subset \G$ 
is given by 
\[
G(2) = \la g_1^2, \, g_2^2, \, g_3^2 \ra, 
\]  
which is referred to as the transformation group of 
level $2$. 
Note that there exists the natural homomorphism 
$G \rightarrow S_3$ defined by $g_i \mt \tau_i = (ij)$, 
whose kernel is just $G(2)$.   
Thus we have a $4$-parameter family of $\G(2)$-actions 
on $\C_x^3$ parametrized by $\th \in \Th$.  
\par
We observe that the polynomial $f(x,\th)$ in (\ref{eqn:f}) 
is $g_i$-invariant and hence the family of cubic surfaces 
$\mathcal{S}$ is stable under the action of $G$. 
Moreover, for each $\th \in \Th$, the cubic surface 
$\mathcal{S}(\th)$ is stable under the action of $\G(2)$. 
So the above action can be restricted to cubic surfaces. 
\begin{definition}[Modular Group Actions]  
\label{def:modular} $\phan{a}$
\begin{enumerate} 
\item The symmetric group $S_3$ of degree $3$ acts on the 
base space $\Th$ by permuting the first three components  
$(\th_1,\th_2,\th_3)$ of $\th$, while 
keeping the fourth component $\th_4$ always fixed. 
\item The full modular group $\G$ acts on the family 
$\pi : \mathcal{S} \to \Th$ of affine cubic surfaces  
through the transformation group $G$, covering the action 
of $S_3$ on $\Th$. 
\item The congruence subgroup $\G(2)$ of level $2$ acts on 
each cubic surface $\mathcal{S}(\th)$ through 
the transformation group $G(2)$, {\it area-preservingly} 
with respect to the Poincar\'e residue $\omega_{\th}$.     
\end{enumerate} 
\end{definition}  
\par
A total picture of these actions is presented in 
Figure \ref{fig:total}. 
The identification of these actions with those in 
\S\ref{subsec:NMIF} (see Table \ref{tab:NMofIMF}) is 
stated in the following manner. 
\begin{lemma}[Braid Versus Modular Group Actions]
\label{thm:identify2}   
With the isomorphisms $(\ref{eqn:RtS})$ and  
$(\ref{eqn:identify})$, the braid group actions and 
the modular group actions are identified, including 
their symplectic structures, as in 
Table $\ref{tab:identify}$. 
\end{lemma}
\begin{table}[h] 
\begin{center} 
\begin{tabular}{|c|c|c|c|}
\hline 
               & half-monodromy & 
full-monodromy &  area form      \\ 
\hline 
braid group action   & $B_3 \car \R_t$ & 
$P_3 \car \R_t(a)$ &  $\Om_{\R_t(a)}$ \\ 
\hline 
modular group action &  $\G \car \mathcal{S}$ & 
$\G(2) \car \mathcal{S}(\th)$ & $\omega_{\th}$ \\ 
\hline 
\end{tabular}
\end{center}
\caption{Identification of the braid and modular 
group actions}  
\label{tab:identify}  
\end{table}
\subsection{Nonlinear Monodromy of Painlev\'e Flow} 
\label{subsec:nmpvi} 
Having Theorems \ref{thm:RHconj} and 
\ref{thm:identify2} in hands, we can easily describe the 
nonlinear monodromy of the non-Riccati part of $\PVI$ in 
terms of the modular group action 
in Definition \ref{def:modular}. 
\begin{theorem}[Nonlinear Monodromy] \label{thm:nlm} 
For any $\k \in \K$, put $\th = \rh(\k) \in \Th$. 
Then the nonlinear monodromy of the non-Riccati Painlev\'e 
flow $\PVI^{\ci}(\k)$ is faithfully represented by the 
$\G(2)$-action on the smooth locus $\mathcal{S}^{\ci}(\th)$ 
of the cubic surface $\mathcal{S}(\th)$ through the 
Riemann-Hilbert correspondence $\RH_{t,\k}$. 
Namely we have the intertwining isomorphism 
\begin{equation} \label{eqn:intertwine}
\begin{CD}
[\, \mbox{$\mathrm{NM}$ $\mathrm{of}$} \,\ \PVI^{\ci}(\k) 
\car \M_t^{\ci}(\k) \,] 
@> \RH_{t,\k} >> [\, \G(2) \car \mathcal{S}^{\ci}(\th) \,],   
\end{CD} 
\end{equation}
where $\mathrm{NM}$ means nonlinear monodromy. 
An image picture of $(\ref{eqn:intertwine})$  
is given in Figure $\ref{fig:nm}$.    
\end{theorem} 
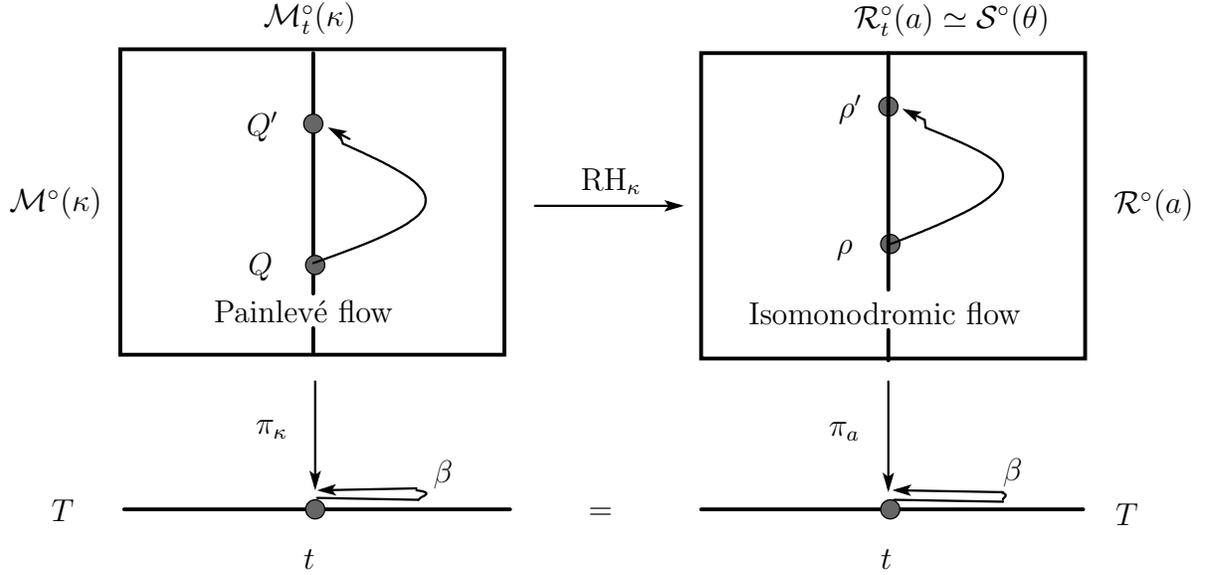
\begin{figure}[t]
\begin{center}
\unitlength 0.1in
\begin{picture}(57.90,28.00)(4.20,-33.40)
%
\special{pn 20}%
\special{pa 1020 3200}%
\special{pa 3040 3200}%
\special{fp}%
%
\special{pn 13}%
\special{pa 3170 1610}%
\special{pa 3930 1610}%
\special{fp}%
\special{sh 1}%
\special{pa 3930 1610}%
\special{pa 3863 1590}%
\special{pa 3877 1610}%
\special{pa 3863 1630}%
\special{pa 3930 1610}%
\special{fp}%
%
\special{pn 20}%
\special{pa 4040 810}%
\special{pa 6050 810}%
\special{pa 6050 2410}%
\special{pa 4040 2410}%
\special{pa 4040 810}%
\special{fp}%
%
\special{pn 20}%
\special{pa 4040 3200}%
\special{pa 6040 3200}%
\special{fp}%
%
\special{pn 20}%
\special{pa 2010 810}%
\special{pa 2010 2070}%
\special{fp}%
%
\special{pn 13}%
\special{pa 2020 2530}%
\special{pa 2020 3080}%
\special{fp}%
\special{sh 1}%
\special{pa 2020 3080}%
\special{pa 2040 3013}%
\special{pa 2020 3027}%
\special{pa 2000 3013}%
\special{pa 2020 3080}%
\special{fp}%
%
\special{pn 13}%
\special{pa 5020 2530}%
\special{pa 5020 3090}%
\special{fp}%
\special{sh 1}%
\special{pa 5020 3090}%
\special{pa 5040 3023}%
\special{pa 5020 3037}%
\special{pa 5000 3023}%
\special{pa 5020 3090}%
\special{fp}%
%
\special{pn 8}%
\special{sh 0.600}%
\special{ar 2020 1920 50 50  0.0000000 6.2831853}%
%
\special{pn 8}%
\special{sh 0.600}%
\special{ar 2020 3200 50 50  0.0000000 6.2831853}%
%
\special{pn 8}%
\special{sh 0.600}%
\special{ar 5030 1810 50 50  0.0000000 6.2831853}%
%
\special{pn 8}%
\special{sh 0.600}%
\special{ar 5030 3200 50 50  0.0000000 6.2831853}%
\put(17.5000,-7.1000){\makebox(0,0)[lb]{$\M_t^{\ci}(\k)$}}%
\put(19.6000,-35.1000){\makebox(0,0)[lb]{$t$}}%
\put(34.1000,-15.5000){\makebox(0,0)[lb]{$\RH_{\k}$}}%
\put(4.2000,-16.8000){\makebox(0,0)[lb]{$\M^{\ci}(\k)$}}%
\put(62.0000,-16.9000){\makebox(0,0)[lb]{$\R^{\ci}(a)$}}%
\put(48.4000,-7.2000){\makebox(0,0)[lb]{$\R_t^{\ci}(a) \simeq \mathcal{S}^{\ci}(\th)$}}%
\put(6.4000,-32.6000){\makebox(0,0)[lb]{$T$}}%
\put(62.1000,-32.8000){\makebox(0,0)[lb]{$T$}}%
\put(17.1000,-28.3000){\makebox(0,0)[lb]{$\pi_{\k}$}}%
\put(47.1000,-28.4000){\makebox(0,0)[lb]{$\pi_a$}}%
\put(16.7000,-20.0000){\makebox(0,0)[lb]{$Q$}}%
\put(47.5000,-18.9000){\makebox(0,0)[lb]{$\rho$}}%
\put(49.8000,-35.1000){\makebox(0,0)[lb]{$t$}}%
\put(14.9000,-22.2000){\makebox(0,0)[lb]{Painlev\'e flow}}%
\put(42.9000,-22.3000){\makebox(0,0)[lb]{Isomonodromic flow}}%
\put(34.5000,-32.5000){\makebox(0,0)[lb]{$=$}}%
%
\special{pn 8}%
\special{sh 0.600}%
\special{ar 2010 1180 50 50  0.0000000 6.2831853}%
%
\special{pn 8}%
\special{sh 0.600}%
\special{ar 5020 1090 50 50  0.0000000 6.2831853}%
\put(16.6000,-12.7000){\makebox(0,0)[lb]{$Q'$}}%
\put(47.6000,-11.8000){\makebox(0,0)[lb]{$\rho'$}}%
%
\special{pn 13}%
\special{pa 2030 3140}%
\special{pa 2560 3150}%
\special{fp}%
%
\special{pn 13}%
\special{pa 2550 3090}%
\special{pa 2050 3100}%
\special{fp}%
\special{sh 1}%
\special{pa 2050 3100}%
\special{pa 2117 3119}%
\special{pa 2103 3099}%
\special{pa 2116 3079}%
\special{pa 2050 3100}%
\special{fp}%
\special{pa 2050 3100}%
\special{pa 2050 3100}%
\special{fp}%
%
\special{pn 13}%
\special{pa 2550 3100}%
\special{pa 2592 3109}%
\special{pa 2610 3120}%
\special{pa 2591 3136}%
\special{pa 2560 3150}%
\special{sp}%
%
\special{pn 13}%
\special{pa 5050 3150}%
\special{pa 5620 3160}%
\special{fp}%
%
\special{pn 13}%
\special{pa 5620 3110}%
\special{pa 5060 3100}%
\special{fp}%
\special{sh 1}%
\special{pa 5060 3100}%
\special{pa 5126 3121}%
\special{pa 5113 3101}%
\special{pa 5127 3081}%
\special{pa 5060 3100}%
\special{fp}%
%
\special{pn 13}%
\special{pa 5610 3110}%
\special{pa 5640 3127}%
\special{pa 5620 3150}%
\special{pa 5620 3150}%
\special{sp}%
%
\special{pn 13}%
\special{pa 2010 1910}%
\special{pa 2053 1894}%
\special{pa 2095 1877}%
\special{pa 2137 1861}%
\special{pa 2178 1844}%
\special{pa 2219 1828}%
\special{pa 2259 1811}%
\special{pa 2298 1795}%
\special{pa 2335 1778}%
\special{pa 2371 1762}%
\special{pa 2405 1745}%
\special{pa 2438 1729}%
\special{pa 2468 1712}%
\special{pa 2497 1696}%
\special{pa 2522 1679}%
\special{pa 2546 1662}%
\special{pa 2566 1646}%
\special{pa 2584 1629}%
\special{pa 2598 1612}%
\special{pa 2609 1596}%
\special{pa 2617 1579}%
\special{pa 2621 1562}%
\special{pa 2621 1545}%
\special{pa 2617 1528}%
\special{pa 2609 1511}%
\special{pa 2597 1494}%
\special{pa 2582 1477}%
\special{pa 2564 1460}%
\special{pa 2543 1443}%
\special{pa 2519 1425}%
\special{pa 2492 1408}%
\special{pa 2463 1391}%
\special{pa 2432 1374}%
\special{pa 2399 1356}%
\special{pa 2365 1339}%
\special{pa 2329 1322}%
\special{pa 2293 1305}%
\special{pa 2255 1287}%
\special{pa 2217 1270}%
\special{pa 2178 1252}%
\special{pa 2150 1240}%
\special{sp}%
%
\special{pn 13}%
\special{pa 2200 1260}%
\special{pa 2104 1209}%
\special{fp}%
\special{sh 1}%
\special{pa 2104 1209}%
\special{pa 2153 1258}%
\special{pa 2151 1234}%
\special{pa 2172 1223}%
\special{pa 2104 1209}%
\special{fp}%
%
\special{pn 13}%
\special{pa 5030 1810}%
\special{pa 5073 1793}%
\special{pa 5116 1775}%
\special{pa 5158 1758}%
\special{pa 5200 1740}%
\special{pa 5241 1723}%
\special{pa 5281 1705}%
\special{pa 5320 1688}%
\special{pa 5358 1671}%
\special{pa 5394 1653}%
\special{pa 5429 1636}%
\special{pa 5461 1619}%
\special{pa 5492 1601}%
\special{pa 5520 1584}%
\special{pa 5546 1567}%
\special{pa 5569 1550}%
\special{pa 5590 1533}%
\special{pa 5607 1515}%
\special{pa 5621 1498}%
\special{pa 5632 1481}%
\special{pa 5639 1464}%
\special{pa 5643 1447}%
\special{pa 5642 1431}%
\special{pa 5638 1414}%
\special{pa 5629 1397}%
\special{pa 5616 1380}%
\special{pa 5600 1363}%
\special{pa 5581 1347}%
\special{pa 5558 1330}%
\special{pa 5533 1314}%
\special{pa 5505 1297}%
\special{pa 5474 1281}%
\special{pa 5442 1264}%
\special{pa 5408 1248}%
\special{pa 5373 1231}%
\special{pa 5337 1215}%
\special{pa 5299 1198}%
\special{pa 5261 1182}%
\special{pa 5223 1165}%
\special{pa 5210 1160}%
\special{sp}%
%
\special{pn 13}%
\special{pa 5210 1160}%
\special{pa 5110 1110}%
\special{fp}%
\special{sh 1}%
\special{pa 5110 1110}%
\special{pa 5161 1158}%
\special{pa 5158 1134}%
\special{pa 5179 1122}%
\special{pa 5110 1110}%
\special{fp}%
%
\special{pn 20}%
\special{pa 1000 790}%
\special{pa 3010 790}%
\special{pa 3010 2390}%
\special{pa 1000 2390}%
\special{pa 1000 790}%
\special{fp}%
\put(26.4000,-30.9000){\makebox(0,0)[lb]{$\b$}}%
\put(56.2000,-30.8000){\makebox(0,0)[lb]{$\b$}}%
%
\special{pn 20}%
\special{pa 5020 2260}%
\special{pa 5020 2420}%
\special{fp}%
%
\special{pn 20}%
\special{pa 5020 810}%
\special{pa 5020 2050}%
\special{fp}%
%
\special{pn 20}%
\special{pa 2010 2250}%
\special{pa 2010 2380}%
\special{fp}%
\end{picture}%
\end{center}
\caption{Nonlinear monodromy of $\PVI^{\ci}(\k)$} 
\label{fig:nm}
\end{figure}
\begin{remark}[Dichotomy] \label{rem:dichotomy} 
Theorem \ref{thm:nlm} means that the global nature 
of $\PVI$ is well understood according to the {\it dichotomy} 
into the Riccati and  non-Riccati components.  
\begin{enumerate} 
\item On the Riccati component $\PVI^{\r}$, the 
flows are linearizable in terms 
of Gauss hypergeometric equations (see Theorem \ref{thm:HGE}), 
whose global nature is well understood classically. 
\item On the non-Riccati component $\PVI^{\ci}$, the 
nonlinear monodromy is faithfully represented by an explicit 
modular group action on cubic surfaces, from which we can 
extract the global nature of $\PVI^{\ci}$.      
\end{enumerate}
After giving a total picture of Riccati solutions 
in \S\ref{sec:riccati}, we shall give a more complete description 
of the nonlinear monodromy of $\PVI(\k)$, when it contains 
the Riccati component, in Theorem \ref{thm:nlm2}. 
\end{remark}
\par
The monodromy problem for $\PVI$ was discussed in     
Dubrovin and Mazzocco \cite{DM} for a special 
one-parameter family and in Iwasaki \cite{Iwasaki4} 
for the full-family. 
Then the solution in Iwasaki \cite{Iwasaki4} has been 
completed in Inaba, Iwasaki and Saito \cite{IIS2} by solving 
the Riemann-Hilbert problem precisely and is now 
presented in this article.      
The global nature of $\PVI$ can also be investigated 
from more analytical point of view, as the 
{\it connection problem}.  
For the latter subject we refer to the important papers 
by Jimbo \cite{Jimbo} and Guzzetti \cite{Guzzetti}. 
We also remark that Jimbo's asymptotic formula was 
corrected by Boalch \cite{Boalch2}. 
\section{Singularities and Riccati Solutions} 
\label{sec:riccati}
We shall classify $(-2)$-curves on moduli spaces 
in terms of resolutions of singularities of cubic surfaces 
by Riemann-Hilbert correspondence. 
Together with the modular group action on 
cubic surfaces, this makes it possible to get a 
total picture of Riccati solutions to $\PVI$. 
\subsection{Singularities of Cubic Surfaces}  
\label{subsec:singularity} 
For our family of cubic surfaces, the discriminant locus 
was calculated by Iwasaki \cite{Iwasaki3}.  
\begin{definition}[Discriminant] \label{def:discriminant} 
Let $\varDelta(\th)$ be the discriminant of the cubic 
surface $\mathcal{S}(\th)$, which is an irreducible polynomial 
of $\th \in \Th$. 
When lifted by the map $A \to \Th$ in $(\ref{eqn:theta})$, 
$\varDelta(\th)$ factors as      
\[
\begin{array}{rcl}
\varDelta(\th) &=& w(a)^2 \ds \prod_{i=1}^4(a_i^2-4), \\[5mm]
w(a) &=& \ds \prod_{\ve_1\ve_2\ve_3=1} 
(\ve_1 a_1+\ve_2 a_2+\ve_3 a_3+a_4) 
-\ds \prod_{i=1}^3 (a_i a_4 - a_j a_k),  
\end{array} 
\]
where $\ve_i = \pm1$ and $\{i,j,k\} = \{1,2,3\}$D
\end{definition}
\begin{lemma}[Discriminant Locus] \label{lem:singularity}  
Tthe Riemann-Hilbert correspondence in the parameter 
level, $\rh : \K \to \Th$, maps $\Wall$ onto the discriminant 
locus $V = \{\, \th \in \Th \,:\, \varDelta(\th) = 0 \,\}$  
$($see Figure $\ref{fig:rh}$$)$.    
For any $\th \in \Th$, the cubic surface $\mathcal{S}(\th)$ 
is singular if and only if $\th \in V$.   
\end{lemma} 
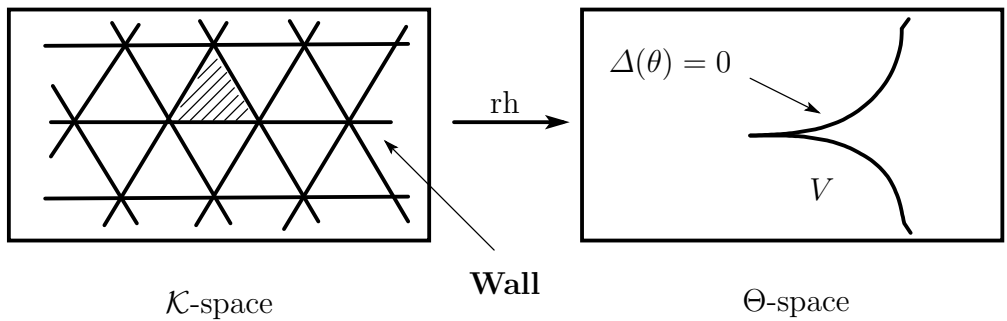
\begin{figure}[t]
\begin{center}
\unitlength 0.1in
\begin{picture}(52.00,14.50)(2.00,-16.60)
%
\special{pn 20}%
\special{pa 200 210}%
\special{pa 2400 210}%
\special{pa 2400 1420}%
\special{pa 200 1420}%
\special{pa 200 210}%
\special{fp}%
%
\special{pn 20}%
\special{pa 3200 210}%
\special{pa 5400 210}%
\special{pa 5400 1420}%
\special{pa 3200 1420}%
\special{pa 3200 210}%
\special{fp}%
%
\special{pn 20}%
\special{pa 420 800}%
\special{pa 2200 800}%
\special{fp}%
%
\special{pn 20}%
\special{pa 420 980}%
\special{pa 880 290}%
\special{fp}%
%
\special{pn 20}%
\special{pa 740 280}%
\special{pa 1360 1330}%
\special{fp}%
%
\special{pn 20}%
\special{pa 1340 280}%
\special{pa 700 1340}%
\special{fp}%
%
\special{pn 20}%
\special{pa 1200 270}%
\special{pa 1820 1320}%
\special{fp}%
%
\special{pn 20}%
\special{pa 1820 280}%
\special{pa 1180 1340}%
\special{fp}%
%
\special{pn 20}%
\special{pa 1670 270}%
\special{pa 2290 1320}%
\special{fp}%
%
\special{pn 20}%
\special{pa 2280 300}%
\special{pa 1640 1360}%
\special{fp}%
%
\special{pn 20}%
\special{pa 870 1330}%
\special{pa 430 610}%
\special{fp}%
%
\special{pn 20}%
\special{pa 390 400}%
\special{pa 2290 390}%
\special{fp}%
%
\special{pn 20}%
\special{pa 390 1200}%
\special{pa 2290 1190}%
\special{fp}%
%
\special{pn 8}%
\special{pa 1380 620}%
\special{pa 1210 790}%
\special{fp}%
\special{pa 1360 580}%
\special{pa 1150 790}%
\special{fp}%
\special{pa 1340 540}%
\special{pa 1090 790}%
\special{fp}%
\special{pa 1310 510}%
\special{pa 1100 720}%
\special{fp}%
\special{pa 1290 470}%
\special{pa 1190 570}%
\special{fp}%
\special{pa 1400 660}%
\special{pa 1270 790}%
\special{fp}%
\special{pa 1430 690}%
\special{pa 1330 790}%
\special{fp}%
\special{pa 1450 730}%
\special{pa 1390 790}%
\special{fp}%
\put(27.2000,-7.5000){\makebox(0,0)[lb]{$\rh$}}%
\put(33.4000,-5.8000){\makebox(0,0)[lb]{$\varDelta(\theta) = 0$}}%
\put(10.2000,-18.3000){\makebox(0,0)[lb]{$\K$-space}}%
\put(40.4000,-18.2000){\makebox(0,0)[lb]{$\Theta$-space}}%
\put(26.1000,-17.0000){\makebox(0,0)[lb]{$\Wall$}}%
%
\special{pn 8}%
\special{pa 4080 570}%
\special{pa 4440 750}%
\special{fp}%
\special{sh 1}%
\special{pa 4440 750}%
\special{pa 4389 702}%
\special{pa 4392 726}%
\special{pa 4371 738}%
\special{pa 4440 750}%
\special{fp}%
%
\special{pn 20}%
\special{pa 4080 870}%
\special{pa 4112 868}%
\special{pa 4145 865}%
\special{pa 4177 863}%
\special{pa 4209 860}%
\special{pa 4241 857}%
\special{pa 4273 853}%
\special{pa 4305 849}%
\special{pa 4337 843}%
\special{pa 4368 837}%
\special{pa 4399 831}%
\special{pa 4430 823}%
\special{pa 4461 814}%
\special{pa 4491 803}%
\special{pa 4521 792}%
\special{pa 4550 779}%
\special{pa 4579 765}%
\special{pa 4607 749}%
\special{pa 4635 732}%
\special{pa 4662 714}%
\special{pa 4688 694}%
\special{pa 4713 674}%
\special{pa 4737 652}%
\special{pa 4760 628}%
\special{pa 4781 604}%
\special{pa 4801 578}%
\special{pa 4820 551}%
\special{pa 4837 523}%
\special{pa 4853 495}%
\special{pa 4867 465}%
\special{pa 4879 435}%
\special{pa 4889 404}%
\special{pa 4898 372}%
\special{pa 4904 341}%
\special{pa 4908 309}%
\special{pa 4910 277}%
\special{pa 4910 260}%
\special{sp}%
%
\special{pn 20}%
\special{pa 4080 870}%
\special{pa 4113 870}%
\special{pa 4145 871}%
\special{pa 4178 871}%
\special{pa 4210 872}%
\special{pa 4242 873}%
\special{pa 4275 875}%
\special{pa 4307 877}%
\special{pa 4339 880}%
\special{pa 4370 883}%
\special{pa 4402 888}%
\special{pa 4434 894}%
\special{pa 4465 900}%
\special{pa 4496 908}%
\special{pa 4527 917}%
\special{pa 4557 927}%
\special{pa 4587 939}%
\special{pa 4617 952}%
\special{pa 4647 967}%
\special{pa 4675 983}%
\special{pa 4703 1001}%
\special{pa 4730 1020}%
\special{pa 4755 1041}%
\special{pa 4779 1063}%
\special{pa 4801 1087}%
\special{pa 4821 1111}%
\special{pa 4838 1138}%
\special{pa 4854 1165}%
\special{pa 4867 1194}%
\special{pa 4879 1223}%
\special{pa 4889 1254}%
\special{pa 4898 1285}%
\special{pa 4906 1316}%
\special{pa 4913 1348}%
\special{pa 4920 1380}%
\special{sp}%
%
\special{pn 8}%
\special{pa 2740 1470}%
\special{pa 2170 900}%
\special{fp}%
\special{sh 1}%
\special{pa 2170 900}%
\special{pa 2203 961}%
\special{pa 2208 938}%
\special{pa 2231 933}%
\special{pa 2170 900}%
\special{fp}%
%
\special{pn 20}%
\special{pa 2530 800}%
\special{pa 3100 800}%
\special{fp}%
\special{sh 1}%
\special{pa 3100 800}%
\special{pa 3033 780}%
\special{pa 3047 800}%
\special{pa 3033 820}%
\special{pa 3100 800}%
\special{fp}%
\put(43.9000,-12.1000){\makebox(0,0)[lb]{$V$}}%
\end{picture}%
\end{center}
\caption{Riemann-Hilbert correspondence in parameter level}
\label{fig:rh} 
\end{figure} 
\begin{table}[b]
\begin{center}
\begin{tabular}{|c||c|c|c|c|c|}
\hline 
number of nodes & 4 & 3 & 2 & 1 & 0 \\ 
\hline 
Dynkin diagram & $D_4$ & $A_3$ & $A_2$ & $A_1$ & $\emptyset$\\ 
\cline{2-6}  
& $A_1^{\op 4}$ & $A_1^{\op 3}$ & $A_1^{\op 2}$ & $-$ & $-$\\ 
\hline  
\end{tabular}
\end{center}
\caption{Dynkin types of singularities} 
\label{tab:type}
\end{table}
\par 
In order to classify the singularities of $\mathcal{S}(\th)$, 
we introduce a stratification of $\K$. 
\begin{definition}[Stratification] 
\label{def:stratification} 
Let $\I$ be the set of all {\sl proper} subset of 
$\{0,1,2,3,4\}$ including the empty set $\emptyset$. 
For each subset $I \in \I$ we put 
\[ 
\begin{array}{rcl}
\K_I &=& \W\mbox{-translates of the subset} \,\, 
\{\, \k \in \K \,:\, \k_i = 0 \,\, (i \in I), \,\, 
\k_i \neq 0 \,\, (i \not\in I) \, \}. 
\\[2mm]
D_I &=& \mbox{Dynkin subdiagram of $D_4^{(1)}$ that has nodes 
precisely in $I$}. 
\end{array}
\]
Then the parameter space $\K$ admits a stratification   
\[
\K = \coprod_{I \in \I} \K_I, \qquad 
\K_{\emptyset} = \K - \Wall \quad (\mbox{the nonsingular locus}). 
\] 
Those Dynkin diagrams which are realized as $D_I$ 
for some $I \in \I$ are precisely the proper subdiagrams 
of $D_4^{(1)}$, tabulated in Table \ref{tab:type}. 
We are interested not only in the subdiagram $D_I$ but also 
in the inclusion pattern $D_I \hookrightarrow D_4^{(1)}$.    
Some typical patterns are illustrated in 
Table \ref{tab:strata}.   
\end{definition} 
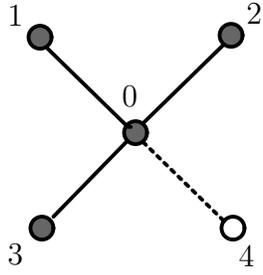
\begin{table}[p]
\begin{center}
\begin{tabular}{ccccc}
& 
\unitlength 0.1in
\begin{picture}(12.50,13.00)(3.30,-15.60)
%
\special{pn 20}%
\special{sh 0.600}%
\special{ar 500 500 60 60  0.0000000 6.2831853}%
%
\special{pn 20}%
\special{sh 0.600}%
\special{ar 1500 490 60 60  0.0000000 6.2831853}%
%
\special{pn 20}%
\special{sh 0.600}%
\special{ar 510 1500 60 60  0.0000000 6.2831853}%
%
\special{pn 20}%
\special{ar 1510 1500 60 60  0.0000000 6.2831853}%
%
\special{pn 20}%
\special{sh 0.600}%
\special{ar 1000 1000 60 60  0.0000000 6.2831853}%
%
\special{pn 20}%
\special{pa 540 560}%
\special{pa 960 970}%
\special{fp}%
\special{pa 970 970}%
\special{pa 970 970}%
\special{fp}%
%
\special{pn 20}%
\special{pa 1040 1050}%
\special{pa 1460 1460}%
\special{dt 0.054}%
\special{pa 1460 1460}%
\special{pa 1460 1460}%
\special{dt 0.054}%
\special{pa 1460 1460}%
\special{pa 1460 1460}%
\special{dt 0.054}%
%
\special{pn 20}%
\special{pa 1460 540}%
\special{pa 1040 950}%
\special{fp}%
%
\special{pn 20}%
\special{pa 950 1050}%
\special{pa 570 1440}%
\special{fp}%
\put(3.3000,-4.4000){\makebox(0,0)[lb]{$1$}}%
\put(15.8000,-4.3000){\makebox(0,0)[lb]{$2$}}%
\put(3.3000,-16.9000){\makebox(0,0)[lb]{$3$}}%
\put(15.4000,-17.0000){\makebox(0,0)[lb]{$4$}}%
\put(9.3000,-8.6000){\makebox(0,0)[lb]{$0$}}%
\end{picture}%
& \aa & 
\unitlength 0.1in
\begin{picture}(12.50,13.00)(3.30,-15.60)
%
\special{pn 20}%
\special{sh 0.600}%
\special{ar 500 500 60 60  0.0000000 6.2831853}%
%
\special{pn 20}%
\special{sh 0.600}%
\special{ar 1500 490 60 60  0.0000000 6.2831853}%
%
\special{pn 20}%
\special{sh 0.600}%
\special{ar 510 1500 60 60  0.0000000 6.2831853}%
%
\special{pn 20}%
\special{sh 0.600}%
\special{ar 1510 1500 60 60  0.0000000 6.2831853}%
%
\special{pn 20}%
\special{ar 1000 1000 60 60  0.0000000 6.2831853}%
%
\special{pn 20}%
\special{pa 540 560}%
\special{pa 960 970}%
\special{dt 0.054}%
\special{pa 960 970}%
\special{pa 960 970}%
\special{dt 0.054}%
\special{pa 970 970}%
\special{pa 970 970}%
\special{dt 0.054}%
%
\special{pn 20}%
\special{pa 1040 1050}%
\special{pa 1460 1460}%
\special{dt 0.054}%
\special{pa 1460 1460}%
\special{pa 1460 1460}%
\special{dt 0.054}%
\special{pa 1460 1460}%
\special{pa 1460 1460}%
\special{dt 0.054}%
%
\special{pn 20}%
\special{pa 1460 540}%
\special{pa 1040 950}%
\special{dt 0.054}%
\special{pa 1040 950}%
\special{pa 1040 950}%
\special{dt 0.054}%
%
\special{pn 20}%
\special{pa 950 1050}%
\special{pa 570 1440}%
\special{dt 0.054}%
\special{pa 570 1440}%
\special{pa 570 1440}%
\special{dt 0.054}%
\put(3.3000,-4.4000){\makebox(0,0)[lb]{$1$}}%
\put(15.8000,-4.3000){\makebox(0,0)[lb]{$2$}}%
\put(3.3000,-16.9000){\makebox(0,0)[lb]{$3$}}%
\put(15.4000,-17.0000){\makebox(0,0)[lb]{$4$}}%
\put(9.3000,-8.6000){\makebox(0,0)[lb]{$0$}}%
\end{picture}%
& \\[6mm]
&  $D_4$ : $I = \{0,1,2,3\}$ & \aa &
  $A_1^{\oplus 4}$ : $I = \{1,2,3,4\}$ & \\[9mm] 
&
\unitlength 0.1in
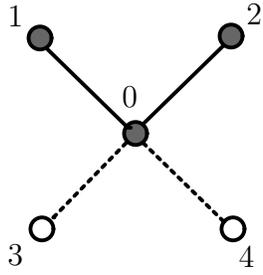
\begin{picture}(12.50,13.00)(3.30,-15.60)
%
\special{pn 20}%
\special{sh 0.600}%
\special{ar 500 500 60 60  0.0000000 6.2831853}%
%
\special{pn 20}%
\special{sh 0.600}%
\special{ar 1500 490 60 60  0.0000000 6.2831853}%
%
\special{pn 20}%
\special{ar 510 1500 60 60  0.0000000 6.2831853}%
%
\special{pn 20}%
\special{ar 1510 1500 60 60  0.0000000 6.2831853}%
%
\special{pn 20}%
\special{sh 0.600}%
\special{ar 1000 1000 60 60  0.0000000 6.2831853}%
%
\special{pn 20}%
\special{pa 540 560}%
\special{pa 960 970}%
\special{fp}%
\special{pa 970 970}%
\special{pa 970 970}%
\special{fp}%
%
\special{pn 20}%
\special{pa 1040 1050}%
\special{pa 1460 1460}%
\special{dt 0.054}%
\special{pa 1460 1460}%
\special{pa 1460 1460}%
\special{dt 0.054}%
\special{pa 1460 1460}%
\special{pa 1460 1460}%
\special{dt 0.054}%
%
\special{pn 20}%
\special{pa 1460 540}%
\special{pa 1040 950}%
\special{fp}%
%
\special{pn 20}%
\special{pa 950 1050}%
\special{pa 570 1440}%
\special{dt 0.054}%
\special{pa 570 1440}%
\special{pa 570 1440}%
\special{dt 0.054}%
\put(3.3000,-4.4000){\makebox(0,0)[lb]{$1$}}%
\put(15.8000,-4.3000){\makebox(0,0)[lb]{$2$}}%
\put(3.3000,-16.9000){\makebox(0,0)[lb]{$3$}}%
\put(15.4000,-17.0000){\makebox(0,0)[lb]{$4$}}%
\put(9.3000,-8.6000){\makebox(0,0)[lb]{$0$}}%
\end{picture}%
& \aa  & 
\unitlength 0.1in
\begin{picture}(12.50,13.00)(3.30,-15.60)
%
\special{pn 20}%
\special{sh 0.600}%
\special{ar 500 500 60 60  0.0000000 6.2831853}%
%
\special{pn 20}%
\special{sh 0.600}%
\special{ar 1500 490 60 60  0.0000000 6.2831853}%
%
\special{pn 20}%
\special{sh 0.600}%
\special{ar 510 1500 60 60  0.0000000 6.2831853}%
%
\special{pn 20}%
\special{ar 1510 1500 60 60  0.0000000 6.2831853}%
%
\special{pn 20}%
\special{ar 1000 1000 60 60  0.0000000 6.2831853}%
%
\special{pn 20}%
\special{pa 540 560}%
\special{pa 960 970}%
\special{dt 0.054}%
\special{pa 960 970}%
\special{pa 960 970}%
\special{dt 0.054}%
\special{pa 970 970}%
\special{pa 970 970}%
\special{dt 0.054}%
%
\special{pn 20}%
\special{pa 1040 1050}%
\special{pa 1460 1460}%
\special{dt 0.054}%
\special{pa 1460 1460}%
\special{pa 1460 1460}%
\special{dt 0.054}%
\special{pa 1460 1460}%
\special{pa 1460 1460}%
\special{dt 0.054}%
%
\special{pn 20}%
\special{pa 1460 540}%
\special{pa 1040 950}%
\special{dt 0.054}%
\special{pa 1040 950}%
\special{pa 1040 950}%
\special{dt 0.054}%
%
\special{pn 20}%
\special{pa 950 1050}%
\special{pa 570 1440}%
\special{dt 0.054}%
\special{pa 570 1440}%
\special{pa 570 1440}%
\special{dt 0.054}%
\put(3.3000,-4.4000){\makebox(0,0)[lb]{$1$}}%
\put(15.8000,-4.3000){\makebox(0,0)[lb]{$2$}}%
\put(3.3000,-16.9000){\makebox(0,0)[lb]{$3$}}%
\put(15.4000,-17.0000){\makebox(0,0)[lb]{$4$}}%
\put(9.3000,-8.6000){\makebox(0,0)[lb]{$0$}}%
\end{picture}%
 & \\[6mm]
& $A_3$ : $I = \{0,1,2\}$  & \aa  & 
$A_1^{\oplus 3}$ : $I = \{1,2,3\}$     & \\[9mm] 
&
\unitlength 0.1in
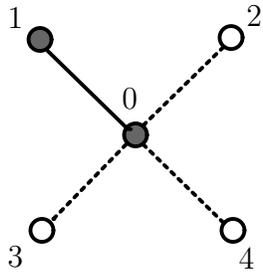
\begin{picture}(12.50,13.00)(3.30,-15.60)
%
\special{pn 20}%
\special{sh 0.600}%
\special{ar 500 500 60 60  0.0000000 6.2831853}%
%
\special{pn 20}%
\special{ar 1500 490 60 60  0.0000000 6.2831853}%
%
\special{pn 20}%
\special{ar 510 1500 60 60  0.0000000 6.2831853}%
%
\special{pn 20}%
\special{ar 1510 1500 60 60  0.0000000 6.2831853}%
%
\special{pn 20}%
\special{sh 0.600}%
\special{ar 1000 1000 60 60  0.0000000 6.2831853}%
%
\special{pn 20}%
\special{pa 540 560}%
\special{pa 960 970}%
\special{fp}%
\special{pa 970 970}%
\special{pa 970 970}%
\special{fp}%
%
\special{pn 20}%
\special{pa 1040 1050}%
\special{pa 1460 1460}%
\special{dt 0.054}%
\special{pa 1460 1460}%
\special{pa 1460 1460}%
\special{dt 0.054}%
\special{pa 1460 1460}%
\special{pa 1460 1460}%
\special{dt 0.054}%
%
\special{pn 20}%
\special{pa 1460 540}%
\special{pa 1040 950}%
\special{dt 0.054}%
\special{pa 1040 950}%
\special{pa 1040 950}%
\special{dt 0.054}%
%
\special{pn 20}%
\special{pa 950 1050}%
\special{pa 570 1440}%
\special{dt 0.054}%
\special{pa 570 1440}%
\special{pa 570 1440}%
\special{dt 0.054}%
\put(3.3000,-4.4000){\makebox(0,0)[lb]{$1$}}%
\put(15.8000,-4.3000){\makebox(0,0)[lb]{$2$}}%
\put(3.3000,-16.9000){\makebox(0,0)[lb]{$3$}}%
\put(15.4000,-17.0000){\makebox(0,0)[lb]{$4$}}%
\put(9.3000,-8.6000){\makebox(0,0)[lb]{$0$}}%
\end{picture}%
& \aa  & 
\unitlength 0.1in
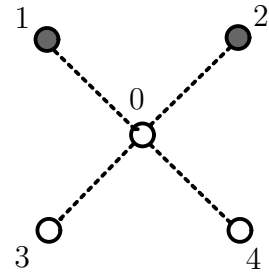
\begin{picture}(12.50,13.00)(3.30,-15.60)
%
\special{pn 20}%
\special{sh 0.600}%
\special{ar 500 500 60 60  0.0000000 6.2831853}%
%
\special{pn 20}%
\special{sh 0.600}%
\special{ar 1500 490 60 60  0.0000000 6.2831853}%
%
\special{pn 20}%
\special{ar 510 1500 60 60  0.0000000 6.2831853}%
%
\special{pn 20}%
\special{ar 1510 1500 60 60  0.0000000 6.2831853}%
%
\special{pn 20}%
\special{ar 1000 1000 60 60  0.0000000 6.2831853}%
%
\special{pn 20}%
\special{pa 540 560}%
\special{pa 960 970}%
\special{dt 0.054}%
\special{pa 960 970}%
\special{pa 960 970}%
\special{dt 0.054}%
\special{pa 970 970}%
\special{pa 970 970}%
\special{dt 0.054}%
%
\special{pn 20}%
\special{pa 1040 1050}%
\special{pa 1460 1460}%
\special{dt 0.054}%
\special{pa 1460 1460}%
\special{pa 1460 1460}%
\special{dt 0.054}%
\special{pa 1460 1460}%
\special{pa 1460 1460}%
\special{dt 0.054}%
%
\special{pn 20}%
\special{pa 1460 540}%
\special{pa 1040 950}%
\special{dt 0.054}%
\special{pa 1040 950}%
\special{pa 1040 950}%
\special{dt 0.054}%
%
\special{pn 20}%
\special{pa 950 1050}%
\special{pa 570 1440}%
\special{dt 0.054}%
\special{pa 570 1440}%
\special{pa 570 1440}%
\special{dt 0.054}%
\put(3.3000,-4.4000){\makebox(0,0)[lb]{$1$}}%
\put(15.8000,-4.3000){\makebox(0,0)[lb]{$2$}}%
\put(3.3000,-16.9000){\makebox(0,0)[lb]{$3$}}%
\put(15.4000,-17.0000){\makebox(0,0)[lb]{$4$}}%
\put(9.3000,-8.6000){\makebox(0,0)[lb]{$0$}}%
\end{picture}%
& \\[6mm]
& $A_2$ : $I = \{0,1\}$    & \aa  & 
$A_1^{\oplus 2}$ : $I = \{1,2\}$       & \\[9mm] 
&
\unitlength 0.1in
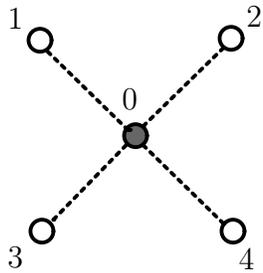
\begin{picture}(12.50,13.00)(3.30,-15.60)
%
\special{pn 20}%
\special{ar 500 500 60 60  0.0000000 6.2831853}%
%
\special{pn 20}%
\special{ar 1500 490 60 60  0.0000000 6.2831853}%
%
\special{pn 20}%
\special{ar 510 1500 60 60  0.0000000 6.2831853}%
%
\special{pn 20}%
\special{ar 1510 1500 60 60  0.0000000 6.2831853}%
%
\special{pn 20}%
\special{sh 0.600}%
\special{ar 1000 1000 60 60  0.0000000 6.2831853}%
%
\special{pn 20}%
\special{pa 540 560}%
\special{pa 960 970}%
\special{dt 0.054}%
\special{pa 960 970}%
\special{pa 960 970}%
\special{dt 0.054}%
\special{pa 970 970}%
\special{pa 970 970}%
\special{dt 0.054}%
%
\special{pn 20}%
\special{pa 1040 1050}%
\special{pa 1460 1460}%
\special{dt 0.054}%
\special{pa 1460 1460}%
\special{pa 1460 1460}%
\special{dt 0.054}%
\special{pa 1460 1460}%
\special{pa 1460 1460}%
\special{dt 0.054}%
%
\special{pn 20}%
\special{pa 1460 540}%
\special{pa 1040 950}%
\special{dt 0.054}%
\special{pa 1040 950}%
\special{pa 1040 950}%
\special{dt 0.054}%
%
\special{pn 20}%
\special{pa 950 1050}%
\special{pa 570 1440}%
\special{dt 0.054}%
\special{pa 570 1440}%
\special{pa 570 1440}%
\special{dt 0.054}%
\put(3.3000,-4.4000){\makebox(0,0)[lb]{$1$}}%
\put(15.8000,-4.3000){\makebox(0,0)[lb]{$2$}}%
\put(3.3000,-16.9000){\makebox(0,0)[lb]{$3$}}%
\put(15.4000,-17.0000){\makebox(0,0)[lb]{$4$}}%
\put(9.3000,-8.6000){\makebox(0,0)[lb]{$0$}}%
\end{picture}%
& \aa  & 
\unitlength 0.1in
\begin{picture}(12.50,13.00)(3.30,-15.60)
%
\special{pn 20}%
\special{sh 0.600}%
\special{ar 500 500 60 60  0.0000000 6.2831853}%
%
\special{pn 20}%
\special{ar 1500 490 60 60  0.0000000 6.2831853}%
%
\special{pn 20}%
\special{ar 510 1500 60 60  0.0000000 6.2831853}%
%
\special{pn 20}%
\special{ar 1510 1500 60 60  0.0000000 6.2831853}%
%
\special{pn 20}%
\special{ar 1000 1000 60 60  0.0000000 6.2831853}%
%
\special{pn 20}%
\special{pa 540 560}%
\special{pa 960 970}%
\special{dt 0.054}%
\special{pa 960 970}%
\special{pa 960 970}%
\special{dt 0.054}%
\special{pa 970 970}%
\special{pa 970 970}%
\special{dt 0.054}%
%
\special{pn 20}%
\special{pa 1040 1050}%
\special{pa 1460 1460}%
\special{dt 0.054}%
\special{pa 1460 1460}%
\special{pa 1460 1460}%
\special{dt 0.054}%
\special{pa 1460 1460}%
\special{pa 1460 1460}%
\special{dt 0.054}%
%
\special{pn 20}%
\special{pa 1460 540}%
\special{pa 1040 950}%
\special{dt 0.054}%
\special{pa 1040 950}%
\special{pa 1040 950}%
\special{dt 0.054}%
%
\special{pn 20}%
\special{pa 950 1050}%
\special{pa 570 1440}%
\special{dt 0.054}%
\special{pa 570 1440}%
\special{pa 570 1440}%
\special{dt 0.054}%
\put(3.3000,-4.4000){\makebox(0,0)[lb]{$1$}}%
\put(15.8000,-4.3000){\makebox(0,0)[lb]{$2$}}%
\put(3.3000,-16.9000){\makebox(0,0)[lb]{$3$}}%
\put(15.4000,-17.0000){\makebox(0,0)[lb]{$4$}}%
\put(9.3000,-8.6000){\makebox(0,0)[lb]{$0$}}%
\end{picture}%
& \\[6mm]
& $A_1$ : $I = \{0\}$      & \aa  & 
$A_1$ : $I = \{1\}$                    & \\[4mm] 
\end{tabular}
\end{center}
\caption{Stratification}
\label{tab:strata}
\end{table}
\par 
Using the stratification in 
Definition \ref{def:stratification}, we can clearly 
describe all the possible singularities. 
\begin{theorem}[Classfication of Singularities] 
\label{thm:classification} 
For any $I \in \I - \{ \emptyset \}$ and any $\k \in \K_I$, 
the surface $\mathcal{S}(\th)$ with 
$\th = \rh(\k)$ has simple singularities of type $D_I$.    
\end{theorem} 
In this situation the Riemann-Hilbert correspondence 
$\RH_{t,\k} : \M_t(\k) \to \mathcal{S}(\th)$ gives a 
minimal resolution of singularities (Theorem \ref{thm:RH3}). 
Thus the exceptional divisor of $\RH_{t,k}$, namely, 
the Riccati locus $\M_t^{\r}(\k) \subset \M_t(\k)$  
has the dual graph of Dynkin type $D_I$. 
We give an example. 
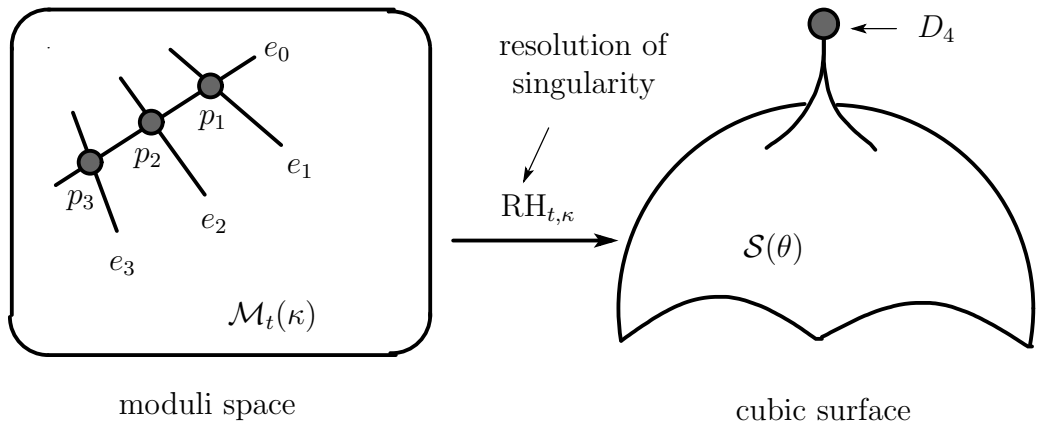
\begin{figure}[t]
\begin{center}
\unitlength 0.1in
\begin{picture}(53.57,19.80)(2.00,-21.80)
%
\special{pn 20}%
\special{pa 400 200}%
\special{pa 2210 200}%
\special{fp}%
%
\special{pn 20}%
\special{pa 400 2010}%
\special{pa 2220 2010}%
\special{fp}%
%
\special{pn 20}%
\special{pa 210 400}%
\special{pa 210 1800}%
\special{fp}%
%
\special{pn 20}%
\special{pa 2400 400}%
\special{pa 2400 1810}%
\special{fp}%
%
\special{pn 8}%
\special{ar 400 400 0 0  3.1415927 3.1941757}%
%
\special{pn 20}%
\special{ar 400 390 190 190  3.0419240 4.7649720}%
%
\special{pn 20}%
\special{ar 400 1810 200 200  1.5707963 3.1415927}%
%
\special{pn 20}%
\special{ar 2200 1800 200 200  6.2332269 6.2831853}%
\special{ar 2200 1800 200 200  0.0000000 1.6162196}%
%
\special{pn 20}%
\special{ar 2200 410 201 201  4.6217291 6.1835167}%
%
\special{pn 20}%
\special{pa 1480 450}%
\special{pa 440 1120}%
\special{fp}%
\put(15.3000,-4.7000){\makebox(0,0)[lb]{$e_0$}}%
\put(16.5000,-10.8000){\makebox(0,0)[lb]{$e_1$}}%
\put(12.0000,-13.7000){\makebox(0,0)[lb]{$e_2$}}%
\put(7.2000,-15.9000){\makebox(0,0)[lb]{$e_3$}}%
\put(7.7000,-23.5000){\makebox(0,0)[lb]{moduli space}}%
%
\special{pn 20}%
\special{pa 4460 276}%
\special{pa 4460 308}%
\special{pa 4459 341}%
\special{pa 4458 373}%
\special{pa 4457 405}%
\special{pa 4455 437}%
\special{pa 4451 469}%
\special{pa 4447 501}%
\special{pa 4441 532}%
\special{pa 4434 563}%
\special{pa 4424 594}%
\special{pa 4413 624}%
\special{pa 4401 654}%
\special{pa 4386 683}%
\special{pa 4370 712}%
\special{pa 4352 739}%
\special{pa 4334 766}%
\special{pa 4313 792}%
\special{pa 4292 817}%
\special{pa 4269 840}%
\special{pa 4245 863}%
\special{pa 4220 883}%
\special{pa 4195 903}%
\special{pa 4168 921}%
\special{pa 4160 926}%
\special{sp}%
%
\special{pn 20}%
\special{pa 4460 276}%
\special{pa 4460 310}%
\special{pa 4460 345}%
\special{pa 4460 379}%
\special{pa 4461 413}%
\special{pa 4462 446}%
\special{pa 4464 480}%
\special{pa 4467 512}%
\special{pa 4472 545}%
\special{pa 4477 576}%
\special{pa 4483 607}%
\special{pa 4491 637}%
\special{pa 4501 667}%
\special{pa 4513 695}%
\special{pa 4526 722}%
\special{pa 4542 749}%
\special{pa 4559 774}%
\special{pa 4579 798}%
\special{pa 4600 821}%
\special{pa 4623 844}%
\special{pa 4647 866}%
\special{pa 4671 887}%
\special{pa 4697 909}%
\special{pa 4722 930}%
\special{pa 4730 936}%
\special{sp}%
%
\special{pn 20}%
\special{ar 4460 1766 1075 1075  2.9825695 4.6192017}%
%
\special{pn 20}%
\special{ar 4490 1766 1067 1067  4.7504658 6.2831853}%
\special{ar 4490 1766 1067 1067  0.0000000 0.1790153}%
%
\special{pn 20}%
\special{pa 3400 1936}%
\special{pa 3429 1916}%
\special{pa 3458 1895}%
\special{pa 3487 1875}%
\special{pa 3516 1855}%
\special{pa 3545 1836}%
\special{pa 3573 1817}%
\special{pa 3602 1800}%
\special{pa 3631 1782}%
\special{pa 3660 1766}%
\special{pa 3689 1751}%
\special{pa 3718 1738}%
\special{pa 3747 1725}%
\special{pa 3776 1714}%
\special{pa 3805 1705}%
\special{pa 3834 1697}%
\special{pa 3863 1692}%
\special{pa 3892 1688}%
\special{pa 3921 1686}%
\special{pa 3950 1687}%
\special{pa 3979 1689}%
\special{pa 4008 1694}%
\special{pa 4037 1700}%
\special{pa 4066 1709}%
\special{pa 4095 1719}%
\special{pa 4124 1730}%
\special{pa 4153 1743}%
\special{pa 4182 1757}%
\special{pa 4211 1773}%
\special{pa 4240 1789}%
\special{pa 4269 1806}%
\special{pa 4298 1824}%
\special{pa 4328 1843}%
\special{pa 4357 1862}%
\special{pa 4386 1882}%
\special{pa 4415 1902}%
\special{pa 4444 1922}%
\special{pa 4450 1926}%
\special{sp}%
%
\special{pn 20}%
\special{pa 4450 1916}%
\special{pa 4480 1899}%
\special{pa 4510 1882}%
\special{pa 4541 1866}%
\special{pa 4571 1849}%
\special{pa 4601 1833}%
\special{pa 4631 1818}%
\special{pa 4661 1803}%
\special{pa 4691 1790}%
\special{pa 4721 1776}%
\special{pa 4751 1764}%
\special{pa 4781 1753}%
\special{pa 4811 1744}%
\special{pa 4841 1735}%
\special{pa 4871 1728}%
\special{pa 4901 1723}%
\special{pa 4931 1719}%
\special{pa 4960 1716}%
\special{pa 4990 1716}%
\special{pa 5020 1718}%
\special{pa 5049 1721}%
\special{pa 5079 1726}%
\special{pa 5108 1733}%
\special{pa 5138 1741}%
\special{pa 5167 1751}%
\special{pa 5196 1762}%
\special{pa 5225 1774}%
\special{pa 5255 1788}%
\special{pa 5284 1802}%
\special{pa 5313 1818}%
\special{pa 5342 1834}%
\special{pa 5371 1851}%
\special{pa 5400 1869}%
\special{pa 5429 1887}%
\special{pa 5458 1906}%
\special{pa 5487 1925}%
\special{pa 5516 1944}%
\special{pa 5545 1963}%
\special{pa 5550 1966}%
\special{sp}%
\put(40.0000,-23.5000){\makebox(0,0)[lb]{cubic surface }}%
\put(27.7000,-13.3000){\makebox(0,0)[lb]{$\RH_{t,\k}$}}%
\put(27.6000,-4.4000){\makebox(0,0)[lb]{resolution of}}%
\put(28.3000,-6.7000){\makebox(0,0)[lb]{singularity}}%
%
\special{pn 20}%
\special{ar 4460 276 0 0  1.4056476 6.2831853}%
\special{ar 4460 276 0 0  0.0000000 1.2490458}%
%
\special{pn 20}%
\special{sh 0.600}%
\special{ar 4460 276 76 76  0.0000000 6.2831853}%
\put(40.4000,-15.4000){\makebox(0,0)[lb]{$\mathcal{S}(\theta)$}}%
%
\special{pn 20}%
\special{pa 2520 1410}%
\special{pa 3350 1410}%
\special{fp}%
\special{sh 1}%
\special{pa 3350 1410}%
\special{pa 3283 1390}%
\special{pa 3297 1410}%
\special{pa 3283 1430}%
\special{pa 3350 1410}%
\special{fp}%
\put(13.4000,-18.9000){\makebox(0,0)[lb]{$\M_t(\k)$}}%
%
\special{pn 20}%
\special{pa 1040 410}%
\special{pa 1620 910}%
\special{fp}%
%
\special{pn 20}%
\special{pa 780 560}%
\special{pa 1220 1170}%
\special{fp}%
%
\special{pn 20}%
\special{pa 530 740}%
\special{pa 760 1360}%
\special{fp}%
\put(49.5000,-3.7000){\makebox(0,0)[lb]{$D_4$}}%
%
\special{pn 8}%
\special{pa 4840 300}%
\special{pa 4610 300}%
\special{fp}%
\special{sh 1}%
\special{pa 4610 300}%
\special{pa 4677 320}%
\special{pa 4663 300}%
\special{pa 4677 280}%
\special{pa 4610 300}%
\special{fp}%
%
\special{pn 8}%
\special{pa 3020 820}%
\special{pa 2890 1100}%
\special{fp}%
\special{sh 1}%
\special{pa 2890 1100}%
\special{pa 2936 1048}%
\special{pa 2912 1052}%
\special{pa 2900 1031}%
\special{pa 2890 1100}%
\special{fp}%
%
\special{pn 20}%
\special{sh 0.600}%
\special{ar 1250 600 60 60  0.0000000 6.2831853}%
%
\special{pn 20}%
\special{sh 0.600}%
\special{ar 940 790 60 60  0.0000000 6.2831853}%
%
\special{pn 20}%
\special{sh 0.600}%
\special{ar 620 1000 60 60  0.0000000 6.2831853}%
\put(11.9000,-8.3000){\makebox(0,0)[lb]{$p_1$}}%
\put(8.5000,-10.3000){\makebox(0,0)[lb]{$p_2$}}%
\put(5.0000,-12.5000){\makebox(0,0)[lb]{$p_3$}}%
\end{picture}%
\end{center}
\caption{A singularity of type $D_4$ for $\k = (0,0,0,0,1)$ 
and $\th = (8,8,8,28)$} 
\label{fig:resolution} 
\end{figure}    
\begin{example}[Singularity of Type $\mbox{\bm $D_4$}$] 
\label{ex:D4} 
For $\th = (8,8,8,28)$ the cubic surface 
$\mathcal{S}(\th)$ has a simple singularity of type $D_4$. 
If $\k$ is a $\W$-translate of $(0,0,0,0,1)$, 
then one has $\rh(\k) = \th$ and the 
$(t,\k)$-Riemann-Hilbert correspondence 
$\RH_{t,\k} : \M_t(\k) \to \mathcal{S}(\th)$ gives a minimal 
resolution of singularity as in Figure \ref{fig:resolution}. 
In this case we have the Riccati locus of type $D_4$,
\begin{equation} \label{eqn:excep} 
\M_t^{\r}(\k) = e_0 \cup e_1 \cup e_2 \cup e_3.
\end{equation}
\end{example}
\subsection{A Total Picture of Riccati Solutions} 
\label{subsec:yonmi} 
As is mentioned in \S\ref{subsec:RF}, 
Saito and Terajima \cite{STe} established the relation 
between $(-2)$-curves and Riccati solutions; this is 
an event on the Painlev\'e equation side.  
On the other hand, on the cubic surface side, 
Iwasaki \cite{Iwasaki3,Iwasaki4} pointed out that the 
singular points on $\mathcal{S}(\th)$ are precisely the 
fixed points of the ${\G(2)}$-action. 
Then Inaba, Iwasaki and Saito \cite{IIS2} added one 
more aspect, namely, the aspect of resolution of singularities 
by Riemann-Hilbert correspondence; the last one finds 
itself in the middle of the previous two aspects.    
Combining all these three, we are now able to get a 
total picture of Riccati solutions. 
\begin{theorem}[A Total Picture] \label{thm:total} 
Let $(t,\k) \in T \times \Wall$ and put 
$\th = \rh(\k) \in \Th$. 
\begin{enumerate}  
\item The germs at $t$ of Riccati solutions to $\PVI(\k)$ 
are in one-to-one correspondence with the points on 
the $(-2)$-curves on the moduli space $\M_t(\k)$. 
\item Each $(-2)$-curve on the moduli space $\M_t(\k)$ 
is sent to a singular point on the cubic surface 
$\mathcal{S}(\th)$ by the Riemann-Hilbert 
correspondence 
$\RH_{t,\k} : \M_t(\k) \to \mathcal{S}(\th)$.  
\item Conversely any $(-2)$-curve on the moduli space 
$\M_t(\k)$ arises as an irreducible component of the 
exceptional divisor of the resolution of singularities 
$\RH_{t,\k} : \M_t(\k) \to \mathcal{S}(\th)$.   
\item The singular points on $\mathcal{S}(\th)$ are exactly 
the fixed points of the $\G(2)$-action on $\mathcal{S}(\th)$.
\item The singular points on $\mathcal{S}(\th)$, as well 
as the $(-2)$-curves on $\M_t(\k)$, are completely 
classified as in Theorem $\ref{thm:classification}$.   
\end{enumerate}   
\end{theorem} 
\par 
Theorem \ref{thm:total} can be visualized as  
in Figure \ref{fig:yonmi}. 
The following is a simple application.   
\begin{corollary}[Single-Valued Solutions] \label{cor:riccati} 
Any single-valued solution to $\PVI$ is a Riccati solution 
and moreover it is a rational solution.  
\end{corollary}
{\it Proof}. 
The proof is very easy by now. 
Not a Riccati solution implies not a fixed point, 
implies not single-valued, since the Riemann-Hilbert 
correspondence is one-to-one outside the Riccati locus. 
Hence any single-valued solution must be Riccati. 
Now recall that a Riccati solution is a logarithmic 
derivative of a hypergeometric function 
(see Theorem \ref{thm:HGE}). 
Such a function of regular singular type can be 
single-valued only if it is a rational function. 
The proof ends. \hfill $\Box$
\par\medskip 
The rational solutions to $\PVI$ were classified by 
Mazzocco \cite{Mazzocco}. 
\begin{figure}[t]
\begin{center}
\unitlength 0.1in
\begin{picture}(61.30,29.00)(3.80,-29.00)
%
\special{pn 20}%
\special{pa 390 500}%
\special{pa 2870 500}%
\special{pa 2870 1290}%
\special{pa 390 1290}%
\special{pa 390 500}%
\special{fp}%
%
\special{pn 20}%
\special{pa 2930 880}%
\special{pa 3980 880}%
\special{fp}%
\special{pa 2922 932}%
\special{pa 3980 932}%
\special{fp}%
\put(7.5000,-9.9000){\makebox(0,0)[lb]{Riccati solutions of $\PVI$ }}%
\put(46.5000,-8.7000){\makebox(0,0)[lb]{Singular points of }}%
\put(48.4000,-11.2000){\makebox(0,0)[lb]{cubic surface}}%
%
\special{pn 20}%
\special{pa 2900 2500}%
\special{pa 3950 2500}%
\special{fp}%
\special{pa 2892 2552}%
\special{pa 3950 2552}%
\special{fp}%
\put(6.1000,-25.8000){\makebox(0,0)[lb]{$(-2)$-curves on moduli space}}%
\put(42.9000,-25.9000){\makebox(0,0)[lb]{Fixed points of $\G(2)$-action}}%
\put(26.3000,-17.8000){\makebox(0,0)[lb]{resolution of singularity}}%
%
\special{pn 20}%
\special{pa 380 2110}%
\special{pa 2860 2110}%
\special{pa 2860 2900}%
\special{pa 380 2900}%
\special{pa 380 2110}%
\special{fp}%
%
\special{pn 20}%
\special{pa 4030 2110}%
\special{pa 6510 2110}%
\special{pa 6510 2900}%
\special{pa 4030 2900}%
\special{pa 4030 2110}%
\special{fp}%
%
\special{pn 20}%
\special{pa 4030 500}%
\special{pa 6510 500}%
\special{pa 6510 1290}%
\special{pa 4030 1290}%
\special{pa 4030 500}%
\special{fp}%
%
\special{pn 20}%
\special{pa 1470 1360}%
\special{pa 1470 2060}%
\special{fp}%
\special{pa 1530 1350}%
\special{pa 1530 2070}%
\special{fp}%
%
\special{pn 20}%
\special{pa 5270 1360}%
\special{pa 5270 2060}%
\special{fp}%
\special{pa 5330 1350}%
\special{pa 5330 2070}%
\special{fp}%
\put(9.0000,-1.7000){\makebox(0,0)[lb]{{\bf Painlev\'e side}}}%
\put(44.8000,-1.8000){\makebox(0,0)[lb]{{\bf Cubic surface side}}}%
%
\special{pn 20}%
\special{pa 3570 1620}%
\special{pa 3970 1340}%
\special{fp}%
\special{sh 1}%
\special{pa 3970 1340}%
\special{pa 3904 1362}%
\special{pa 3926 1371}%
\special{pa 3927 1395}%
\special{pa 3970 1340}%
\special{fp}%
%
\special{pn 20}%
\special{pa 2880 2070}%
\special{pa 3320 1790}%
\special{fp}%
\put(33.0000,-1.7000){\makebox(0,0)[lb]{$\mbox{\bm $\RH$}$}}%
%
\special{pn 20}%
\special{pa 2870 120}%
\special{pa 3140 120}%
\special{fp}%
%
\special{pn 20}%
\special{pa 3720 130}%
\special{pa 4050 130}%
\special{fp}%
\special{sh 1}%
\special{pa 4050 130}%
\special{pa 3983 110}%
\special{pa 3997 130}%
\special{pa 3983 150}%
\special{pa 4050 130}%
\special{fp}%
\end{picture}%
\end{center}
\caption{A total picture of Riccati solutions} 
\label{fig:yonmi} 
\end{figure}
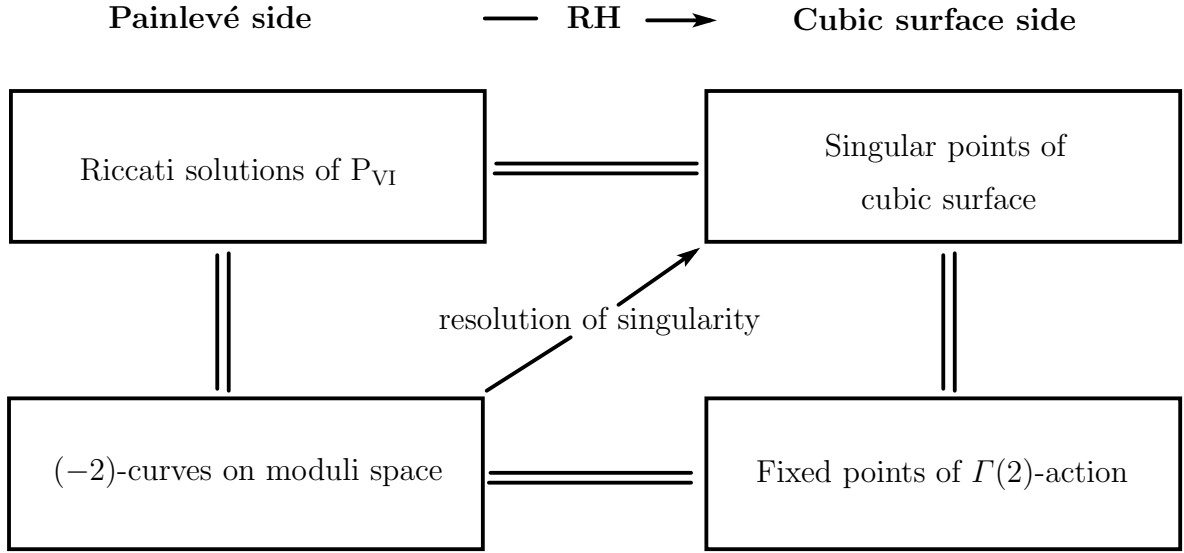
\begin{example}[Some Rational Solutions] \label{ex:rational} 
In Example \ref{ex:D4} any solution on the $(-2)$-curve  
$e_0$ is rational. 
Indeed any half-monodromy $\a$ preserves the Riccati 
configuration (\ref{eqn:excep}) and hence induces an 
automorphism $\b$ of $e_0 \simeq \P^1$ that permutes the 
three points $p_1$, $p_2$, $p_3$ in 
Figure \ref{fig:resolution}. 
If $\a$ is a full-monodromy, then $\b$ fixes each of 
$p_1$, $p_2$, $p_3$ and hence is identity on $e_0$. 
This means that any solution on $e_0$ is single-valued. 
By Corollary \ref{cor:riccati} it is a rational solution.    
\end{example}
\par 
We refer to Lukashevich and Yablonski \cite{LY}, 
Fokas and Ablowitz \cite{FA}, Okamoto \cite{Okamoto4}, 
Watanabe \cite{Watanabe}, Gromak, Laine and 
Shimomura \cite{GLS} and the references therein for explicit 
calculations of Riccati solutions. 
\par
Having established the total picture of Riccati solutions, 
especially, of their connection with Klein sigularities, 
we shall revisit the nonlinear monodromy of $\PVI$ discussed 
in \S\ref{sec:nlmonod}. 
For any $\k \in \K-\Wall$, Theorem \ref{thm:nlm} completely 
describes the nonlinear monodromy of $\PVI(\k)$, since 
$\PVI(\k)$ has no Riccati locus. 
We now deal with the case $\k \in \Wall$. 
While the Riemann-Hilbert correspondence $\RH_{t,\k} : \M_t(\k) 
\to \mathcal{S}(\th)$ with $ \th = \rh(\k)$ gives an 
{\it analytic} minimal resolution of singularities, 
there exists an {\it algebraic} minimal resolution of 
singularities 
\begin{equation} \label{eqn:amrs}
\widetilde{\mathcal{S}}(\th) \rightarrow \mathcal{S}(\th),
\end{equation}
constructed in a standard manner for Klein singularities 
as in Brieskorn \cite{Brieskorn}. 
Then the action of $\G(2)$ on $\mathcal{S}(\th)$ can be 
lifted to $\widetilde{\mathcal{S}}(\th)$ uniquely. 
Combining this fact with Theorem \ref{thm:nlm}, we have:
\begin{theorem}[Nonlinear Monodromy Revisited] \label{thm:nlm2} 
For any $\k \in \Wall$, put $\th = \rh(\k) \in \Th$. 
Then the nonlinear monodromy of $\PVI(\k)$ is faithfully 
represented by the $\G(2)$-action on 
$\widetilde{\mathcal{S}}(\th)$. 
\end{theorem} 
\section{Canonical Coordinates} \label{sec:coordinates}  
The moduli space $\M(\k)$ admits a natural canonical 
coordinate system whose local charts are laveled by 
the affine Weyl group $\W$. 
In this section we shall construct such coordinates 
and write down the Painlev\'e dynamics explicitly in 
terms of them. 
The principle of producing canonical coordinates is 
the Wronskian construction that converts a stable parabolic 
connection to a second-order single Fuchsian differntial 
equation. 
In this section we mean by $T$ the configuration space of 
distinct ordered three points in $\C$ as in (\ref{eqn:T2}) 
upon putting $t_4 = \infty$; hence
\[  
t = (t_1,t_2,t_3) = (t_1,t_2,t_3,\infty) \in T.
\]   
\subsection{Space of Fuchsian Equations} 
\label{subsec:fuchsian}
We start with spaces of Fuchsian equations from which local 
coordinates are to be extracted.   
\begin{definition}[Fuchsian Equations] \label{def:fuchsian} 
For any $\k \in \K$, let $\E(\k)$ be the set of all  
second-order Fuchsian differential equations of the form 
\begin{equation} \label{eqn:fuchsian}
\dfrac{d^2f}{dz^2}-v_1(z)\dfrac{df}{dz}+v_2(z)f = 0, 
\end{equation}
with four regular singular points 
$t = (t_1,t_2,t_3,t_4) \in T$ and an apparent 
singular point $q$, having Riemann scheme as 
in Table \ref{tab:riemann2}, where $\k$ is fixed while 
$t$ and $q$ may vary in such a manner that $q$ does not 
meet any of $t_1$, $t_2$, $t_3$, $t_4$.   
\end{definition} 
\begin{table}[t]
\begin{center}
\begin{tabular}{|c||c|c|c|c||c|} 
\hline 
singularity & $t_1$ & $t_2$ & $t_3$ & $t_4 = \infty$ & $q$ \\ 
\hline 
first exponent & $0$ & $0$ & $0$ & $\k_0$ & $0$ \\ 
\hline 
second exponent & $\k_1$ & $\k_2$ & $\k_3$ & $\k_4+\k_0$ & $2$ \\ 
\hline
\end{tabular}
\end{center}  
\caption{Riemann Scheme}
\label{tab:riemann2}
\end{table}
\par 
The affine linear relation $2 \k_0+\k_1+\k_2+\k_3+\k_4=1$ 
in (\ref{eqn:K}) is exactly Fuchs' relation for 
Fuchsian differential equation (\ref{eqn:fuchsian}). 
The classical Fuchs-Frobenius method in the theory of 
Fuchsian differential equations allows us to determine   
the coefficients $v_1(z)$ and $v_2(z)$ as  
\begin{equation} \label{eqn:coeff} 
v_1(z) = \dfrac{1}{z-q} + \sum_{i=1}^3 \dfrac{\k_i-1}{z-t_i}, 
\qquad  
v_2(z) = \dfrac{p}{z-q} + \sum_{i=1}^3 \dfrac{H_i(\k)}{z-t_i}.  
\end{equation}  
The condition that $q$ is apparent with exponents $0$ and 
$2$ implies that 
$H_i(\k) = H_i(q,p,t;\k)$ is a functions of $(q,p,t,\k)$. 
This function, called the $i$-th {\it Hamiltonian}, is  
explicitly determined as follows.   
\begin{lemma}[Hamiltonians] \label{lem:Hamitonians} 
The $i$-th Hamiltonian $H_i(\k) = H_i(q,p,t;\k)$ is given by 
\begin{equation} \label{eqn:Hamiltonian} 
(t_{ij}t_{ik}) \, H_i(\k) 
= (q_iq_jq_k) p^2  
- \{(\k_i-1) q_jq_k + \k_j q_kq_i + \k_k q_iq_j\} p 
+ \k_0(\k_0+\k_4) q_i, 
\end{equation} 
with $\{i,j,k\}=\{1,2,3\}$, where $q_i := q - t_i$ and 
$t_{ij} := t_i - t_j$. 
\end{lemma}  
\begin{remark}[Polynomial Hamiltonians] \label{rem:polynh}
Note that Hamiltonians (\ref{eqn:Hamiltonian}) 
are polynomials of $(q,p)$.  
This is because one exponent at each finite singular 
point, $t_1$, $t_2$, $t_3$, $q$, is {\it zero} in the 
Riemann scheme of Table \ref{tab:riemann2}.  
We call (\ref{eqn:Hamiltonian}) {\it polynomial 
Hamiltonians} (see Okamoto \cite{Okamoto2}).     
\end{remark} 
\par
Formulas (\ref{eqn:coeff}) and (\ref{eqn:Hamiltonian}) 
tell us that for a fixed $\k$, the Fuchsian equation 
(\ref{eqn:fuchsian}) is determined uniquely by 
the data $(q,p,t)$. 
Thus the following definition would be natural.  
\begin{definition}[Canonical Coordinates] \label{def:CC} 
The set $\E(\k)$ is identified with the affine variety  
\begin{equation} \label{eqn:U}   
U = \{\,(q,p,t) \in \C_q \times \C_p \times \C_t^3\,:\, 
q \neq t_i, \,\, 
t_i \neq t_j \,\,\, \mbox{for} \,\,\, i \neq j \,\}, 
\end{equation} 
having coordinates $(q,p,t)$, on which the fundamental 
$2$-form is defined by 
\begin{equation} \label{eqn:OE}
\Om_{\E(\k)} = 
dq \wedge dp - \ds \sum_{i=1}^3 dH_i(q,p;\k) \wedge dt_i  
\end{equation}  
\end{definition}
\subsection{Wronskian Construction}  
\label{subsec:wronskian}
First we shall define the concept of apparent singular 
point of a stable parabolic connection.   
Let $Q = (E,\nabla,\psi,l) \in \M(\k)$ be any stable 
parabolic connection.   
By the stability of $Q$, we can show that there exists 
a unique line subbundle $F \subset E$ of maximal degree.  
The line bundle $F$ is called the maximal subbundle  
of $E$ and the quotient bundle $L = E/F$ is called the 
minimal quotient bundle of $E$.  
Note that $F$ and $L$ are of degrees $0$ and $-1$, 
respectively. 
Let 
\begin{equation} \label{eqn:proj}
\pi : E \to L = E/F
\end{equation}
be the canonical projection.   
We see that the composite 
$u : F \to L \ot \Om^1_{\P^1}(D_t)$ of the sequence  
\[
\begin{CD}
F @> \sc{\mathrm{inclusion}} >> E @> \nabla >> 
E \ot \Om^1_{\P^1}(D_t) @> \pi \ot 1 >> 
L \ot \Om^1_{\P^1}(D_t)    
\end{CD}
\] 
is an $\O_{\P^1}$-homomorphism and gives a holomorphic 
section of the line bundle 
$\Hom(F, L) \ot \Om^1_{\P^1}(D_t)$, where $t \in T$ is 
the regular singular points of $Q$. 
Then the stability of $Q$ implies that $u$ is a 
nontrivial section. 
Since the line bundle 
$\Hom(F, L) \ot \Om^1_{\P^1}(D_t)$ is of degree 
one, the section $u$ has a unique simple 
zero $q$. 
Since the construction so far is canonical, the 
point $q = q(Q) \in \P^1$ is uniquely determined by 
$Q \in \M(\k)$. 
Hence we have a well-defined morphism
\begin{equation} \label{eqn:apparent} 
q : \M(\k) \rightarrow \P^1, \qquad Q \mt q = q(Q).  
\end{equation}      
\begin{definition}[Apparent Singular Point] 
\label{def:apparent} The point $q = q(Q) \in \P^1$ in 
(\ref{eqn:apparent}) is called the 
{\it apparent singular point} of the stable parabolic 
connection $Q \in \M(\k)$.   
\end{definition} 
\par
Using the morphism (\ref{eqn:apparent}), we can consider  
the locus $\M^{\id}(\k) \subset \M(\k)$ where the apparent 
singular point $q$ does not meet any regular singular 
point $t_i$, $i = 1,2,3,4$, that is,     
\[
\M^{\id}(\k) = \{\, Q \in \M(\k) \,:\, q(Q) \neq t_i(Q) 
\quad (i = 1,2,3,4) \, \}, 
\]  
where $t_i = t_i(Q)$ denotes the $i$-th regular singular 
point of $Q$.  
\par 
Next we proceed to the Wronskian construction that 
recast each stable parabolic connection in 
$\M^{\id}(\k)$ to a Fuchsian differential equation 
in $\E(\k)$. 
Given a stable parabolic connection 
$Q = (E,\nabla,\psi,l)  \in \M(\k)$, we consider the 
locally constant sheaf 
\[
\L' = 
\Ker\,\left[\, \nabla|_{\P^1-D_t} : E|_{\P^1-D_t} 
\rightarrow E|_{\P^1-D_t} \ot \Om^1_{\P^1-D_t} 
\,\right]    
\]
of $\nabla$-horizontal sections on $\P^1-D_t$. 
By the stability of $Q$ we can show that the canonical 
projection $\pi$ in (\ref{eqn:proj}) induces an 
isomorphism of locally constant sheaves on $\P^1-D_t$, 
\[
\pi : \L' \rightarrow \pi(\L') \subset L|_{\P^1-D_t} 
\]
On the other hand, since the line bundle $L^{-1}$ is 
of degree one, there exists a unique connection  
$\delta : L^{-1} \to L^{-1} \ot \Om^1_{\P^1}(D_t)$ 
whose residue at each singular point is given by
\[
\Res_{t_i}(\delta) = \left\{
\begin{array}{cl}
-\l_i \qquad & (i = 1,2,3), \\[2mm]
\l_1+\l_2+\l_3-1 \qquad & (i=4),  
\end{array}
\right.  
\]
where $\l_i$ is given by (\ref{eqn:lambda}) in 
terms of $\k_i$. 
Let $\L''$ be the locally constant sheaf  
\[
\L'' = 
\Ker\,\left[\, \delta|_{\P^1-D_t} : L^{-1}|_{\P^1-D_t} 
\rightarrow L^{-1}|_{\P^1-D_t} \ot \Om^1_{\P^1-D_t} 
\,\right] \subset L^{-1}|_{\P^1-D_t}.      
\]
of $\delta$-horizontal sections on $\P^1-D_t$. 
Tensoring $\pi(\L')$ with $\L''$, we have a 
locally constant sheaf 
\begin{equation} \label{eqn:locconst} 
\L_Q = \pi(\L') \ot \L'' \subset \O_{\P^1-D_t},  
\end{equation}
canonically associated to the stable parabolic 
connection $Q \in \M(\k)$ (see Remark \ref{rem:shift} 
for the meaning of the tensoring with $\L''$). 
The construction so far is valid for any $Q \in \M(\k)$, 
but we have to put $Q$ in $\M^{\id}(\k)$ to obtain the 
following theorem.     
\begin{theorem}[Wronskian Isomorphism] 
\label{thm:wronskian} 
For a stable parabolic connection $Q \in \M^{\id}(\k)$, 
with singular points at $t \in T$,  
we consider the second-order monic differential 
equation on $\P^1-D_t$ whose solution sheaf is given 
by the locally constant sheaf $\L_Q$ in 
$(\ref{eqn:locconst})$. 
Then it is exactly such a Fuchsian differential equation 
as is formulated in Definition $\ref{def:fuchsian}$ with 
apparent singular point at $q = q(Q)$ given by 
$(\ref{eqn:apparent})$. 
Therefore there exists a well-defined morphism   
\begin{equation} \label{eqn:Phik}
\Phi_{\k} : \M^{\id}(\k) \to \E(\k). 
\end{equation} 
This morphism becomes an isomorphism. 
\end{theorem} 
\par
Combined with the identification $\E(\k) \simeq U$ in 
Definition \ref{def:CC}, where the set $U$ is defined 
by (\ref{eqn:U}), the isomorphism (\ref{eqn:Phik}) 
yields a local coordinate mapping  
\begin{equation} \label{eqn:coordMid}
\Psi_{\k} : \M^{\id}(\k) \to U, \quad Q \mapsto (q,p,t).  
\end{equation}
\par 
At the end of this subsection we emphasize that stability 
has been used many times in the Wronskian construction. 
In addition the following technical remark may be 
helpful. 
\begin{remark}[Shift of Exponents] \label{rem:shift} 
The essential factor in (\ref{eqn:locconst}) is the 
rank-two local system $\pi(\L')$, which is tensored with 
the rank-one local system $\L''$ just for shifting the 
exponents. 
By the tensoring with $\L''$, the exponents in 
Table \ref{tab:riemann} are shifted to those in 
Table \ref{tab:riemann2} by the vector   
$(\l_1,\l_2,\l_3,1-\l_1-\l_2-\l_3)$ at    
$t = (t_1,t_2,t_3,t_4)$. 
This process is necessary to obtain polynomial 
Hamiltonians as in (\ref{eqn:Hamiltonian}) 
(see Remark \ref{rem:polynh}).   
\end{remark}  
\subsection{Canonical Coordinate System}  
\label{subsec:coord} 
Combined with B\"acklund transformations, 
Theorem \ref{thm:wronskian} produces a canonical 
coordinate system on the moduli space $\M(\k)$. 
To see this, for each $\si \in \W$, consider the 
open subset 
\[
\M^{\si}(\k) = s_{\si}^{-1}(\M^{\id}(\si(\k))) \subset 
\M(\k),  
\]
where $s_{\si} : \M(\k) \to \M(\si(\k))$ is the B\"acklund 
transformation corresponding to $\si$ 
(see Figure \ref{fig:backlund}). 
Then there exists an open covering of the moduli 
space $\M(\k)$, 
\[ 
\M(\k) = \bigcup_{\si \in \W} \M^{\si}(\k)    
\] 
On each open subset $\M^{\si}(\k)$ we have an isomorphism 
\begin{equation} \label{eqn:poisson2} 
\Phi_{\k}^{\si} : \M^{\si}(\k) \rightarrow \E(\si(\k)),     
\end{equation}
defined as the composite of the sequence of isomorphisms
\[
\begin{CD} 
\M^{\si}(\k) @> s_{\si} >> \M^{\id}(\si(\k)) 
@> \Phi_{\si(\k)} >> \E(\si(\k)).  
\end{CD}
\] 
\par 
To see that (\ref{eqn:poisson2}) is a Poisson isomorphism, 
we make use of the following theorem. 
\begin{theorem}[Pull-Back Principle] \label{thm:PBP} 
Let $\k \in \K$ and put $a = \rh(\k) \in A$. 
We define the local Riemann-Hilbert correspondence 
$\RH_{\k}^{\si} : \E(\si(\k)) \to \R(a)$ as the 
composite of the sequence
\[
\begin{CD}
\E(\si(\k)) @> (\Phi_{\k}^{\si})^{-1} >> 
\M^{\si}(\k) \hookrightarrow \M(\k) @> \RH_{\k} >> \R(a). 
\end{CD}
\]
Then the fundamental $2$-form $\Om_{\E(\si(\k))}$ on 
$\E(\si(\k))$ is the pull-back of $\Om_{\R(a)}$ by 
$\RH_{\k}^{\si}$, 
\[
\Om_{\E(\si(\k))} = (\RH_{\k}^{\si})^* \Om_{\R(a)}. 
\]
\end{theorem}
\par
This theorem is due to Iwasaki \cite{Iwasaki2}, where 
the map $\RH_{\k}^{\si}$ is formulated\footnote{To be more 
precise, it is formulated only for $\si = \id$, but the 
modification for a general $\si$ is obvious.} directly 
without passing through the moduli space $\M(\k)$. 
Hence we have the commutative diagram
\[
\begin{CD}
\M^{\si}(\k) @> \sc{\mathrm{inclusion}} >> \M(\k) \\
@V \Phi_{\k}^{\si} VV  @VV \RH_{\k} V \\
\E(\si(\k)) @>> \RH_{\k}^{\si} > \R(a), 
\end{CD}
\]
where $\RH_{\k}$ is Poisson by Theorem \ref{thm:PF} 
while $\RH_{\k}^{\si}$ is also Poisson by Theorem 
\ref{thm:PBP} respectively. 
Therefore $\Phi_{\k}^{\si}$ becomes a Poisson isomorphism 
as desired. 
\begin{definition}[Canonical Coordinate System]   
\label{def:CCS} 
By the same procedure as in (\ref{eqn:coordMid}) 
the Poisson isomorphisms (\ref{eqn:poisson2}) induce local 
coordinate mappings 
\begin{equation} \label{eqn:CC2}
\Psi^{\si}_{\k} : \M^{\si}(\k) \to U_{\si}, 
\quad Q \mapsto (q_{\si},p_{\si},t) \qquad (\si \in \W),   
\end{equation} 
where $U_{\si}$ is a copy of $U$ endowed with the 
coordinates $(q_{\si},p_{\si},t)$. 
Note that we have $\Psi^{\si}_{\k} = 
\Psi_{\si(\k)} \ci s_{\si}$ with $\Psi_{\k}$ given 
by (\ref{eqn:coordMid}). 
The collection of maps (\ref{eqn:CC2}) is referred to as 
the {\it canonical coordinate system} on the moduli 
space $\M(\k)$.   
\end{definition}
\par 
We are now in a position to derive a Hamiltonian system of 
differential equations for $\PVI(\k)$ on each local chart 
$\M^{\si}(\k) \simeq U_{\si}$ based on the idea in 
Remark \ref{rem:hamiltonian}. 
\begin{theorem}[Hamiltonian System] \label{thm:HVI} 
In terms of the canonical coordinates $(q_{\si}, p_{\si},t)$ 
on $\M^{\si}(\k) \simeq U_{\si}$, the Painlev\'e flow 
$\PVI(\k)$ is expressed as the Hamiltonian system 
$\HVI^{\si}(\k)$, 
\begin{equation} \label{eqn:HVI} 
\dfrac{\partial q_{\si}}{\partial t_i} = 
\dfrac{\partial H_i(\si(\k))}{\partial p_{\si}}, 
\qquad 
\dfrac{\partial p_{\si}}{\partial t_i} = 
-\dfrac{\partial H_i(\si(\k))}{\partial q_{\si}},  
\end{equation}
with Hamiltonians 
$H_i(\si(\k)) = H_i(q_{\si},p_{\si},t;\si(\k))$ where  
$H_i(q,p,t;\k)$ is given by $(\ref{eqn:Hamiltonian})$.  
\end{theorem}
\par 
Indeed the Painlev\'e flow $\PVI(\k)$ is characterized by 
the condition that $\iota_{v} \Om_{\M(\k)} = 0$ for every 
$\F_{\PVI(\k)}$-horizontal vector field $v$. 
Since (\ref{eqn:poisson2}) is a Poisson isomorphism, this 
condition is equivalent to 
$\iota_{v} \Om_{\E(\si(\k))} = 0$, from which system 
(\ref{eqn:HVI}) readily follows. 
\begin{remark}[Malmquist Expression]  \label{rem:malmquist}
Malmquist \cite{Malmquist} obtained a Hamiltonian expression 
for $\PVI$ as early as 1923.  
Our expression (\ref{eqn:HVI}) is just a symmetric form of 
his expression that can be reduced to his original by the 
symplectic reduction in Remark \ref{rem:reduction}.  
Malmquist's expression was rediscovered by 
Okamoto \cite{Okamoto2,Okamoto3} in the isomonodromic context. 
Deriving Hamiltonian systems as in (\ref{eqn:HVI}) by the 
pull-back principle in Theorem \ref{thm:PBP} is due to 
Iwasaki \cite{Iwasaki2}, where he works on a 
Riemann surface of arbitrary genus. 
\end{remark} 
\begin{theorem}[Analytic Painlev\'e Property] 
\label{thm:APPHVI} 
For any $\k \in \K$ and $\si \in \W$ the Hamiltonian 
system $\HVI^{\si}(\k)$ has analytic Painlev\'e property. 
\end{theorem} 
\par
This theorem immediately follows from the geometric 
Painlev\'e property of the Painlev\'e flow $\F_{\PVI(\k)}$ 
(see Theorem \ref{thm:pp2}) and the algebraicity of 
the phase space $\M(\k)$ (see Remark \ref{rem:relation}).  
\begin{theorem}[Basic B\"acklund Transformations] 
\label{thm:BVI}
In terms of canonical charts in $(\ref{eqn:coordMid})$, 
\[ 
\M^{\id}(\k) \simeq U \simeq \M^{\id}(\si_i(\k)) 
\qquad (i = 0,1,2,3,4),  
\]
the basic B\"acklund transformation $s_i$ is expressed as 
the birational canonical transformation   
\begin{equation} \label{eqn:BVI} 
s_i(\k_j) = \k_j - \k_i \, c_{ij}, \qquad
s_i(t_j)  = t_j, \qquad 
s_i(q_j)  = q_j + \dfrac{\k_i}{q_i} \, u_{ij}.   
\end{equation}  
where $C = (c_{ij})$ is the Cartan matrix of type $D_4^{(1)}$ 
$($see Figure $\ref{fig:dynkin}$$)$ and   
\[ 
q_i = \left\{
\begin{array}{ll}
p     \quad & (i = 0), \\[2mm]
q-t_i \quad & (i = 1,2,3,4), 
\end{array}
\right. 
\qquad 
u_{ij} = \{q_i,q_j\} = 
\dfrac{\partial q_i}{\partial p} \dfrac{\partial q_j}{\partial q}-
\dfrac{\partial q_i}{\partial q} \dfrac{\partial q_j}{\partial p}.
\]   
\end{theorem} 
\par
As is mentioned after Theorem \ref{thm:backlund}, it is not 
so easy to derive the formula (\ref{eqn:BVI}) for $s_0$ 
from our definition of B\"acklund transformations in 
Definition \ref{def:backlund}. 
Our strategy is the coalescence of regular singular 
points along isomonodromic flow; see Inaba, Iwasaki and 
Saito \cite{IIS1}.    
\par  
Traditionally, such coordinate expressions as in 
(\ref{eqn:HVI}) and (\ref{eqn:BVI}) have been a 
starting point of the story. 
In the other way round, we end up with coordinate expressions 
as concrete realizations of the abstract dynamical 
system $\PVI$ that is defined conceptually. 
\begin{remark}[Gluing by B\"acklund Transformations] 
\label{rem:glue} 
The moduli space $\M(\k)$ is made up of local charts glued 
by B\"acklund transformations. 
Indeed it is clear from Definition \ref{def:CCS} that for 
$\si$, $\si' \in \W$ the transition 
function from $\M^{\si}(\k) \simeq U_{\si}$ to 
$\M^{\si'}(\k) \simeq U_{\si'}$ is just the B\"acklund 
transformation $s_{\si'\si^{-1}} = s_{\si'}s_{\si}^{-1}$.  
Noumi, Takano and Yamada \cite{NTY} showed that 
``the manifold of Painlev\'e system'' can be 
constructed in this way.   
Their empirical observation is trivial from our point of view, 
or even from the meta-physics: the phase space of a dynamical 
system should be made up of inertial coordinates glued together 
by symmetries of the system.  
\end{remark} 
\par 
The construction of moduli spaces and that of 
canonical coordinates tell us the following:  
\begin{remark}[Systems or Single Equations?] 
\label{rem:sysingle} 
In doing isomonodromic deformations, \linebreak 
some people work with first-order linear systems as 
in (\ref{eqn:schlesinger}),   
while others work with second-order single 
equations as in (\ref{eqn:fuchsian}). 
We may ask which approach is better. 
The answer is that both are important and necessary. 
Systems are sophisticated to stable parabolic 
connections and are used to construct the phase space of 
Painlev\'e dynamics; while single equations are used to 
construct canonical coordinates of the phase space that 
make it possible to get a concrete realization of the 
dynamics. 
Therefore both are important and necessary.   
\end{remark}  
\section{Summary} \label{sec:summary}
In this article we have observed the natural manner in 
which the {\it continuous} Hamiltonian system $\PVI$ 
induces two {\it discrete} Hamiltonian systems: 
\begin{enumerate}
\item {\bf B\"acklund transformations} as convering 
transformations of the Riemann-Hilbert correspondence. 
They describe the symmetries of $\PVI$.  
\item {\bf Poincar\'e return maps} (or the {\bf nonlinear 
monodromy}).    
Through the Riemann-Hilbert correspondence, they are realized 
as an area-preserving action of the modular group on smooth 
affine cubic surfaces, or on the minimal desingularizations 
of singular affine cubic surfaces. 
They describe the global structures, especially the 
multi-valuedness, of trajectories of $\PVI$.   
\end{enumerate} 
Here we recall that the geometric Painlev\'e property is 
a prerequisite for the well-definedness of 
Poincar\'e return maps. 
In this respect we have shown by the conjugacy method 
that 
\begin{enumerate}
\item[(3)] the {\bf geometric Painlev\'e property} of the 
Painlev\'e flow follows from that of the isomonodromic 
flow, which holds trivially, through the Riemann-Hilbert 
correspondence.   
\end{enumerate}
As to the Riccati component of $\PVI$ 
comprising the classical torajectories that can be 
linearized in terms of Gauss hypergeometric equations, 
we have given 
\begin{enumerate}
\item[(4)] a {\bf total picture of Riccati solutions}  
in terms of resolutions of singularities by the 
Riemann-Hilbert correspondence. 
\end{enumerate}
Concerning the concrete realization of the Painlev\'e 
dynamics, we have constructed   
\begin{enumerate}
\item[(5)] a {\bf canonical coordinate system} via the 
Wronskian construction, in terms of which Hamiltonian 
systems and B\"acklund transformations 
are written down explicitly. 
The nonlinear monodromy is also made explicit in terms 
of {\bf cubic surfaces} as in item (2).    
\end{enumerate}
\par
\begin{figure}[t]
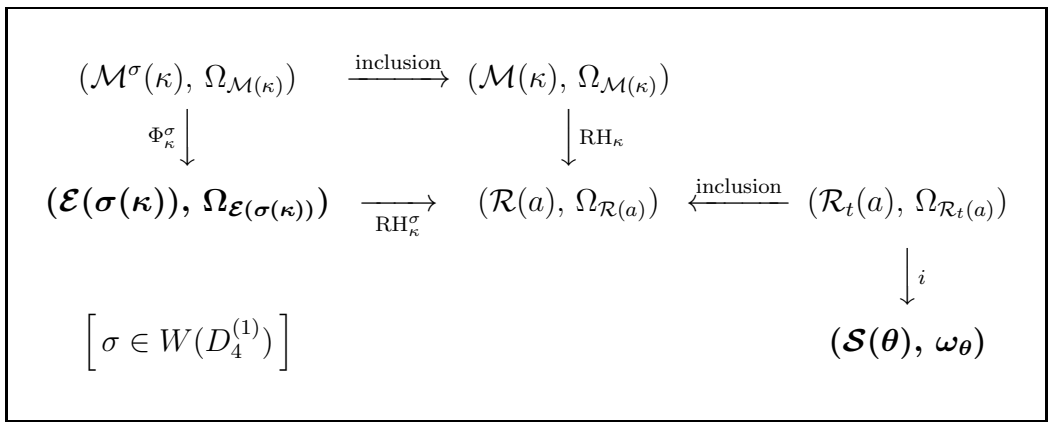

\begin{center}
\fbox{\begin{minipage}{.80\linewidth}
\[
\begin{CD}
(\M^{\si}(\k),\,\Om_{\M(\k)}) 
@> \incl >> (\M(\k),\,\Om_{\M(\k)}) @.  \\
@V \Phi_{\k}^{\si} VV       @VV \RH_{\k} V   @.  \\
\mbox{\bm$(\E(\si(\k)),\,\Om_{\E(\si(\k))})$}  
@>> \RH_{\k}^{\si} > (\R(a),\,\Om_{\R(a)}) 
@< \incl << (\R_t(a),\,\Om_{\R_t(a)}) \\
       @.              @.         @VV i V  \\
\left[\, \si \in \W \,\right] @. @. 
\mbox{\bm$(\mathcal{S}(\th),\,\omega_{\th})$}  
\end{CD}  
\]
\vspace{1mm}
\end{minipage}}
\end{center} 
\caption{Concluding Diagram: bold-faced are concrete objects} 
\label{fig:closing} 
\end{figure}

 
We started the main body of this article with the 
Guiding Diagram in Figure \ref{fig:guiding}. 
We wish to close the article with the Concluding Diagram 
in Figure \ref{fig:closing}. 
Located in the central position of the diagram, 
as well as of the development of our story, is 
the Riemann-Hilbert correspondence 
\[
\RH_{\k} : (\M(\k), \Om_{\M(\k)}) \to 
(\R(a), \Om_{\R(a)}) 
\]
in the precise moduli-theoretical setting.  
The Painlev\'e dynamics encoded in this abstract object 
is concretized on both sides of the diagram: 
Hamiltonian systems and B\"acklund transformations 
on the left-hand side, while nonlinear monodromy on the 
right-hand side, respectively.      
\par\bigskip 
{\bf Acknowledgment}. 
The authors are grateful to P. Boalch for his valuable 
comments on this article. 

\par\vspace{0.5cm}\noindent 
Michi-aki Inaba: Faculty of Mathematics, Kyushu University, 
Hakozaki, Higashi-ku, Fukuoka 812-8581 Japan; 
inaba@math.kyushu-u.ac.jp 
\par\bigskip\noindent
Katsunori Iwasaki: Faculty of Mathematics, Kyushu University, 
Hakozaki, Higashi-ku, Fukuoka 812-8581 Japan; 
iwasaki@math.kyushu-u.ac.jp 
\par\bigskip\noindent
Masa-Hiko Saito: Department of Mathematics, Faculty of 
Science, Kobe University, Rokko-dai, Nada-ku, Kobe 
657-8501 Japan; mhsaito@math.kobe-u.ac.jp 
\end{document}